\UseRawInputEncoding
\documentclass[a4paper,reqno,10pt,oneside]{amsart}
\usepackage{amsmath}
\usepackage[left=2.61cm,right=2.61cm,top=2.72cm,bottom=2.72cm]{geometry}
\usepackage[colorlinks=true,urlcolor=blue,
citecolor=red,linkcolor=blue,linktocpage,pdfpagelabels,
bookmarksnumbered,bookmarksopen]{hyperref}
\usepackage[english]{babel}
\usepackage{graphicx}
\usepackage[toc,page,header]{appendix}
\usepackage{mathrsfs}
\usepackage{amssymb,amsmath, amsfonts, amsthm}
\usepackage{latexsym}

\allowdisplaybreaks[1]

\numberwithin{equation}{section}

\renewcommand{\O}{\operatorname{O}}

\renewcommand{\(}{\left(}
\renewcommand{\)}{\right)}
\renewcommand{\[}{\left[}
\renewcommand{\]}{\right]}

\newtheorem{theorem}{Theorem}[section]
\newtheorem{proposition}[theorem]{Proposition}
\newtheorem{corollary}[theorem]{Corollary}
\newtheorem{lemma}[theorem]{Lemma}
\theoremstyle{definition}
\newtheorem{remark}[theorem]{Remark}
\theoremstyle{definition}
\newtheorem{definition}[theorem]{Definition}
\theoremstyle{definition}
\newtheorem{example}[theorem]{Example}

\renewcommand{\l }{\langle}
\renewcommand{\r }{\rangle}

\renewcommand{\O}{\mathcal{O}}

\newcommand{\T}{{\mathcal T}}

\newcommand{\N}{\mathbb{N}}
\renewcommand{\S}{\mathbb{S}}
\newcommand{\E}{\mathcal{E}}

\newcommand{\V}{\mathcal{V}}

\newcommand{\eps}{\varepsilon}

\renewcommand{\P}{\mathcal{P}}

\newcommand{\beq}{\begin{equation}}
\newcommand{\eeq}{\end{equation}}
\newcommand{\beqs}{\begin{equation*}}
\newcommand{\eeqs}{\end{equation*}}
\newcommand{\beqn}{\begin{eqnarray}}
\newcommand{\eeqn}{\end{eqnarray}}
\newcommand{\beqns}{\begin{eqnarray*}}
\newcommand{\eeqns}{\end{eqnarray*}}
\newcommand{\bdoc}{\begin{document}}
\newcommand{\edoc}{\end{document}}
\newcommand{\be}{\begin{enumerate}}
\newcommand{\ee}{\end{enumerate}}
\newcommand{\bdescr}{\begin{description}}
\newcommand{\edescr}{\end{description}}
\newcommand{\ba}{\begin{array}}
\newcommand{\ea}{\end{array}}
\newcommand{\intR}{\int_{\mathbb R^N}}
\newcommand{\R}{\mathbb R}
\newcommand{\RN}{\mathbb{R}^N}
\newcommand{\B}{\mathbb B}
\newcommand{\C}{\mathbb C}
\renewcommand{\H}{\mathcal H}
\renewcommand{\L}{\mathbb L}
\newcommand{\parallelsum}{\mathbin{\!/\mkern-5mu/\!}}
\newcommand{\e}{\varepsilon}
\newcommand{\SD}{\Sigma_D}
 \renewcommand{\(}{\left(}
\renewcommand{\)}{\right)}
\renewcommand{\[}{\left[}
\renewcommand{\]}{\right]}
\renewcommand{\appendixpagename}{\centering Appendix}

\begin{document}
\title[Existence of nonradial domains for overdetermined and isoperimetric problems]{Existence of nonradial domains for overdetermined and isoperimetric problems in nonconvex cones}

\author{ A\MakeLowercase{lessandro} Iacopetti, F\MakeLowercase{ilomena} Pacella, T\MakeLowercase{obias} Weth}

\thanks{\emph{Acknowledgements.} Research partially supported by Gruppo Nazionale per l'Analisi Matematica, la Pro\-ba\-bi\-li\-t\`a e le loro Applicazioni (GNAMPA) of the Istituto Nazionale di Alta Matematica (INdAM)}

\address[Alessandro Iacopetti]{Dipartimento di Matematica ``G. Peano", Universit\`a di Torino, Via Carlo Alberto 10, 10123 Torino, Italy}
\email{alessandro.iacopetti@unito.it}

\address[Filomena Pacella]{Dipartimento di Matematica, Sapienza Universit\`a di Roma,  P.le Aldo Moro 2, €"00185 Roma, Italy}
\email{pacella@mat.uniroma1.it}

\address[Tobias Weth]{Institut für Mathematik, Goethe-Universität Frankfurt, Robert-Mayer-Str. 10 D-60629 Frankfurt am Main, Germany}
\email{weth@math.uni-frankfurt.de}

\begin{abstract}
In this work we address the question of the existence of nonradial domains inside a nonconvex cone for which a mixed boundary overdetermined problem admits a solution. Our approach is variational, and consists in proving the existence of nonradial minimizers, under a volume constraint, of the associated torsional energy functional. In particular we give a condition on the domain $D$ on the sphere spanning the cone which ensures that the spherical sector is not a minimizer. Similar results are obtained for the relative isoperimetric problem in nonconvex cones. 
\end{abstract}

\maketitle

\section{Introduction}
In this paper we study an overdetermined problem for domains in a cone. This topic shares similarities with the question of characterising constant mean curvature hypersurfaces inside a cone (see \cite{PT,PT2}) and hence with the isoperimetric problem. Thus we will also show some results for it.\\

Let $D$ be a smooth domain on the unit sphere $\S^{N-1}$ and let $\Sigma_D$ be the cone spanned by $D$, namely
\beq\label{def:conespanned}
\Sigma_D:=\{x \in \R^{N}; \ x=s q,\ q\in D,\ s\in(0,+\infty)\}.
\eeq
For a domain $\Omega\subset \Sigma_D$ we set:
$$\Gamma_\Omega:=\partial \Omega \cap \Sigma_D, \ \Gamma_{1,\Omega}:=\partial\Omega\cap\partial\Sigma_D $$
and assume that $\mathcal{H}_{N-1}(\Gamma_{1,\Omega})>0$, where $\mathcal{H}_{N-1}(\cdot)$ denotes the $(N-1)$-dimensional Hausdorff measure. The set $\Gamma_\Omega$ is usually called the relative (to $\Sigma_D$) boundary of $\Omega$.

We consider the following overdetermined mixed boundary value problem:

\beq\label{eq:overdetermined}
\begin{cases}
-\Delta u = 1 & \text{in $\Omega$},\\
 u = 0 & \text{on $\Gamma_\Omega$},\\
\frac{\partial u}{\partial \nu} = 0 & \text{on $\Gamma_{1,\Omega}\setminus\{0\}$},\\
 \frac{\partial u}{\partial \nu} = -c<0 & \text{on $\Gamma_\Omega$},
\end{cases}
\eeq
for a constant $c>0$, where $\nu$ is the exterior unit normal. If $\Gamma_\Omega$ is not smooth then the constant normal derivative condition is understood to hold on the regular part of $\Gamma_\Omega$.

The overdetermined problem \eqref{eq:overdetermined} arises naturally in the study of critical points of a relative torsional energy of subdomains of the cone $\Sigma_D$ subject to a fixed volume contraint. Indeed, for any domain $\Omega$, as in \eqref{eq:overdetermined}, let us consider the torsion problem with mixed boundary conditions
\begin{equation}\label{eq:mixbvprobintro}
\begin{cases}
-\Delta u = 1 & \text{in $\Omega$,}\\
 u = 0 & \text{on $\Gamma_\Omega$,}\\
 \frac{\partial u}{\partial \nu} = 0 & \text{on $\Gamma_{1,\Omega}\setminus\{0\}$.}
\end{cases}
\end{equation}
It is easy to see that \eqref{eq:mixbvprobintro} has a unique weak solution $u_\Omega$ in the Sobolev space $H_0^1(\Omega; \Sigma_D)$ (see Sect. 6 or \cite{BV}), which is obtained by minimizing the functional
\beq\label{eq:defenergyJ}
J(v):= \frac{1}{2}\int_\Omega |\nabla v|^2 \ dx - \int_\Omega v \ dx.
\eeq
We then define the value 
\begin{equation}
  \label{eq:def-torsional-energy-intro}
\E(\Omega; \Sigma_D):=J(u_\Omega)=-\frac{1}{2} \int_\Omega |\nabla u_\Omega|^2 \ dx=-\frac{1}{2}\int_\Omega u_\Omega \ dx,
\end{equation}
and call it the {\em torsional energy of $\Omega$ in $\Sigma_D$}. Note here that the second and third equality in (\ref{eq:def-torsional-energy-intro}) hold since $u_\Omega$ is a weak solution of \eqref{eq:mixbvprobintro}. By definition, the domain-dependent functional $\Omega \mapsto \E(\Omega; \Sigma_D)$ represents a {\em relative} version of the classical torsional energy functional usually defined using the solution of the analogous Dirichlet problem.

Using domain derivative techniques, as for other similar problems in shape-optimization theory it can be proved that the critical points of the functional $\E(\Omega; \Sigma_D)$ with respect to volume-preserving deformations which leave the cone invariant, correspond to domains $\Omega$ for which $\frac{\partial u_\Omega}{\partial \nu}$ is constant on $\Gamma_\Omega$, i.e. $u_\Omega$ satisfies the overdetermined problem \eqref{eq:overdetermined} (see in \cite[Proposition 4.3]{PT2} if $\Gamma_\Omega$ is smooth and $u_\Omega$ has some Sobolev regularity, or Proposition \ref{prop:regmin} in the present paper in the nonsmooth case).

In this paper we intend to study the existence and the properties of domains for which a solution of \eqref{eq:overdetermined} exists. It is easy to see that for any spherical sector $\Omega_{D,R}:=B_R\cap \Sigma_D$, where $B_R=B_R(0)$ is the ball with radius $R>0$ centered at the origin (which is the vertex of the cone), the radial function 
\beq\label{eq:defradialsol}
u(x)=\frac{N^2 c^2-|x|^2}{2N}
\eeq
is a solution of \eqref{eq:overdetermined} for $\Omega=\Omega_{D,R}$. Therefore the question is whether the spherical sectors $\Omega_{D,R}$ are the only domains for which \eqref{eq:overdetermined} admits a solution. In the case of convex cones the answer was provided in \cite{PT}, obtaining the following result (see \cite[Theorem 1.1]{PT}):

\begin{theorem}\label{teo:pacellatralli}
If $\Sigma_D$ is convex, $\Gamma_\Omega$ is smooth and $u$ is a classical solution of \eqref{eq:overdetermined} such that $u \in W^{1,\infty}(\Omega)\cap W^{2,2}(\Omega)$ then:
$$\Omega=\Sigma_D\cap B_R(P_0),$$
where $B_R(P_0)$ is the ball centered at $P_0$ with radius $R=Nc$, and either $P_0=0$, i.e. $\Omega=\Omega_{D,R}$, or $P_0\in\partial \Sigma_D$ and $\Omega$ is a half-ball lying on flat part of $\partial\Sigma_D$.
\end{theorem}
Hence, if $\Sigma_D$ is a convex cone, not flat anywhere, then the radial domains $\Omega_{D,R}$ are the only domains admitting solutions of \eqref{eq:overdetermined}. Let us observe that the assumption $u \in W^{1,\infty}(\Omega)\cap W^{2,2}(\Omega)$ can be seen as a ``gluing condition". Indeed it is automatically satisfied whenever $\Gamma_\Omega$ and $\partial\Sigma_D$ intersect orthogonally (see \cite[Sect. 6]{PT}).

In the context of the variational formulation of problem \eqref{eq:overdetermined} described above, the result of Theorem \ref{teo:pacellatralli} gives a characterization of the smooth critical points of $\E(\Omega; \Sigma_D)$, restricted to the class of subdomains of fixed volume, in the case of convex cones.  In particular any local minimizer of $\E(\Omega; \Sigma_D)$ with a volume constraint is a spherical sector. Actually, using symmetrization methods in cones \cite{LPT, PTR} it can be proved (see \cite{PT2}) that this holds in a more general class of cones which are the ones having an isoperimetric property.

In contrast, the case of nonconvex cones is largely unexplored, which is the main motivation of the present paper. The variational formulation of the overdetermined problem suggests that to look for nonradial domains for which there exists a solution of  \eqref{eq:overdetermined} is equivalent to look for nonradial critical points of $\E(\Omega; \Sigma_D)$ under a volume constraint. In particular, if there are cones for which a minimizer of $\E(\Omega; \Sigma_D)$ (fixing the volume) exists and if we are able to show that it is not the spherical sector then we achieve our goal. This is the content of our first main result. 

Let us denote by $\lambda_1(D)$ the first nontrivial eigenvalue of the Laplace-Beltrami operator $-\Delta_{\S^{N-1}}$ on $D$ with zero Neumann boundary condition.
\begin{theorem}\label{mainteo}
If $D$ is a smooth domain of $\S^{N-1}$ such that:
\beq\label{eq:hpmainteo}
\lambda_1(D)<N-1 \ \ \hbox{and} \ \ \mathcal{H}_{N-1}(D)<\mathcal{H}_{N-1}(\S^{N-1}_+)
\eeq
where $\S^{N-1}_+$ is a half unit sphere, then there exists a bounded domain $\Omega^*$ which is a minimizer for $\E(\Omega; \Sigma_D)$ with a fixed volume, but $\Omega^*$ is not a spherical sector $\Omega_{D,R}$, for $R>0$.

Moreover there exists a critical dimension $d^*$ which can be either $5,6$ or $7$, such that for the relative boundary $\Gamma_{\Omega^*}$ it holds that:
 \begin{itemize}
\item[(i)] $\Gamma_{\Omega^*}$ is smooth if $N<d^*$;
\item[(ii)] $\Gamma_{\Omega^*}$ can have countable isolated singularities if $N=d^*$;
\item[(iii)] $\Gamma_{\Omega^*}$ can have a singular set of dimension $N-d^*$, if $N>d^*$.
\end{itemize}
In addition on the regular part of $\Gamma_{\Omega^*}$ the normal derivative $\frac{\partial u_{\Omega^*}}{\partial \nu}$ is constant, where $u_{\Omega^*}$ is the torsion function of $\Omega^*$.
\end{theorem}
The condition $\lambda_1(D)<N-1$ in \eqref{eq:hpmainteo} is the one which ensures that a spherical sector $\Omega_{D,R}$ cannot be a local minimizer for $\E(\Omega; \Sigma_D)$ among the class of smooth subdomains of $\Sigma_D$ with fixed volume, because it implies that it is not a stable critical point with respect to volume-preserving deformations (see Theorem \ref{prop:critfirsteigenminim}). To prove this, we restrict the torsional energy functional to the class of strictly star-shaped sets $\Omega$ in $\Sigma_D$ with fixed volume $c>0$, and we show the instability of the spherical sector $\Omega_{D,R}$ with $|\Omega_{D,R}|=c$ within this class. The reason to consider strictly star-shaped domains is that the relative boundary $\Gamma_\Omega$ of a strictly star-shaped set is a radial graph of a function $\varphi$ on $D$. This allows to study $\E(\Omega; \Sigma_D)$ as a functional on $\varphi \in C^2(\overline{D})$.

On the other hand, the condition $\mathcal{H}_{N-1}(D)<\mathcal{H}_{N-1}(\S^{N-1}_+)$ is the one which allows to prove the existence of a minimizer for $\E(\Omega; \Sigma_D)$ (see Theorem \ref{teo:mainteoexistmin} and Corollary \ref{cor:mainteoexistmin}). In the Appendix we give examples of domains $D$ on $\S^{N-1}$ satisfying both conditions in \eqref{eq:hpmainteo}.

Let us observe that, since $\Sigma_D$ is not bounded, the existence of a minimizer for $\E(\Omega; \Sigma_D)$ is not obvious. To prove Theorem \ref{mainteo} we use the concentration-compactness principle of P. L. Lions (see \cite{LIO}). It was first used in shape-optimization Dirichlet problems in \cite{BU}. Having mixed boundary conditions, we cannot make use of the same proof as in \cite{BU}. We also stress that, as the cone $\Sigma_D$ is not convex and since we do not have any information on the contact angle between $\Sigma_D$ and $\Gamma^*$, some care is needed to prove that the normal derivative $\frac{\partial u_{\Omega^*}}{\partial \nu}$ of the torsion function $u_{\Omega^*}$ is constant on the regular part of $\Gamma^*$ (see Proposition \ref{prop:regmin}). Finally, the regularity statements follow from the results of \cite{DJ, JS} and \cite{W}.\\

As announced we also consider the isoperimetric problem in the cone to get a analogous nonradiality result using the same strategy.

The isoperimetric problem in the cone consists in minimizing the relative perimeter $\P(E; \Sigma_D)$ among all possibile finite relative perimeter sets $E$ contained in the cone $\Sigma_D$, with a fixed volume. It was proved in \cite{LP}, and later in \cite{RR, FI, CRS}, that if $\Sigma_D$ is a convex cone then the only minimizer of $\P(E; \Sigma_D)$ with a fixed volume are the spherical sectors $\Omega_{D,R}$. This holds also in ``almost" convex cones as shown in \cite{BF} (see also \cite{PT2}). If the cone is not convex, a counterexample is given in \cite{LP}.

Here we show that under the same conditions \eqref{eq:hpmainteo} a minimizer of $\P(E;\Sigma_D)$ exists but is not the spherical sector $\Omega_{D,R}$. So we have:
\begin{theorem}\label{mainteo2}
If $D$ is a smooth domain of $\S^{N-1}$ such that \eqref{eq:hpmainteo} holds then there exists a bounded set of finite perimeter $E^*$ inside $\Sigma_D$ which minimizes $\P(E;\Sigma_D)$ for any fixed volume and $E^*$ is not a spherical sector $\Omega_{D,R}$, $R>0$. Moreover for the relative boundary $\Gamma_{E^*}$ it holds:
\begin{itemize}
\item[(i)] $\Gamma_{E^*}$ can have a closed singular set $\widetilde{\Gamma}_{E^*}$ of Hausdorff dimension less than or equal to $N-7$;\vspace{2pt}
\item[(ii)] $\Gamma_{E^*}\setminus \widetilde{\Gamma}_{E^*}$ is a smooth embedded hypersurface with constant mean curvature;\vspace{3pt}
\item[(iii)] if $x\in \overline{\Gamma_{E^*}\setminus \widetilde{\Gamma}_{E^*}} \cap \partial \Sigma_D$ then $\Gamma_{E^*}\setminus \widetilde{\Gamma}_{E^*}$ is a smooth CMC embedded hypersurface with boundary in a neighborhood of $x$ and meets $\partial\Sigma_D$ orthogonally.
\end{itemize}
\end{theorem}
As for Theorem \ref{mainteo}, the condition $\lambda_1(D)<N-1$ is the one which ensures that $\Omega_{D,R}$ cannot be a local minimizer (see Theorem \ref{prop:nonisopercrit}) and to prove this we again work in the class of smooth star-shaped sets. Instead the existence follows by results obtained in \cite{RR}, while the regularity of minimizers derives from classical results for isoperimetric problems. 

As a consequence of Theorem \ref{mainteo2} we get that whenever  \eqref{eq:hpmainteo} holds there exists a CMC hypersurface in the cone, namely $\Gamma_{E^*}$, intersecting $\partial\Sigma_D$ orthogonally, which is not a spherical cap centered at the vertex of the cone. It is important to notice that $\Gamma_{E^*}$ cannot be a smooth radial graph. Indeed, by \cite[Theorem 1.3]{PT} and \cite[Theorem 1.1]{PT2}, we know that if $\Gamma_{E^*}$ was a CMC radial graph intersecting $\partial\Sigma_D$ orthogonally then $E^*$ would be a spherical sector $\Omega_{D,R}$, and this holds in any cone without requiring convexity hypotheses. It would be very interesting to understand what kind of CMC hypersurface $\Gamma_{E^*}$ could be.

Finally we observe that, from our results and \cite[Theorem 1.1]{LP} (or \cite[Theorem 1.1]{PT}), we easily recover the inequality $\lambda_1(D)\geq N-1$ whenever $D$ is convex. This was proved in \cite[Theorem 4.3]{ES} (see also \cite[Theorem 4.1]{AM}).\\

The paper is organized as follows. In Section 2 we provide some geometric preliminaries. In Section 3 we study the torsional energy functional $\E(\Omega; \Sigma_D)$ on strictly star-shaped domains in the cone; while in Section 4 we derive the formulas for the first and second variations of $\E(\Omega; \Sigma_D)$ when the volume is fixed. In Section 5 we prove that the first condition in \eqref{eq:hpmainteo} allows to prove that the spherical sector is not a local minimizer for $\E(\Omega; \Sigma_D)$. The long Section 6 is devoted to study the question of the existence of minimizers of $\E(\Omega; \Sigma_D)$ with a volume constraint. Their properties are described in Section 7 where the proof of Theorem \ref{mainteo} is deduced. Finally in Section 8 we study the isoperimetric problem and prove Theorem \ref{mainteo2}. In the Appendix we give examples of nonconvex domains satisying the condition \eqref{eq:hpmainteo}.

\section{Some preliminaries}
In this section we fix some notation and we collect, for the reader's convenience, some definitions and known facts from Riemannian Geometry that will be used throughout the paper.
 
Given a smooth manifold $M$, we denote by $T_pM$ the tangent space at $p\in M$, by $\mathcal{T}(M)$ the space of tangent vector fields on $M$ and by $TM$ the tangent bundle.

We denote by $\l\cdot,\cdot\r$ or $\boldsymbol{\cdot}$ the standard scalar product in $\R^{N}$, by $|\cdot|$ the Euclidean norm, and by $\nabla^0$ the flat connection of $\R^{N}$. In the special case $M=D$, where $D\subset\S^{N-1}$ is a domain of the unit sphere in $\RN$, we denote by $\nabla$ the induced Levi-Civita connection on $D$, namely $$\nabla_X Y:=(\nabla_{X}^0Y)^\top,\ \text{for any $X,Y \in \T(D)$,}$$ 
 where $\top:T\R^{N}\to TD$ is the orthogonal projection. 
If we further assume that $D$ is a proper and smooth domain of $\S^{N-1}$ it will be always understood that $D$ is considered as a submanifold with boundary, equipped with the induced Riemannian metric. 
 
  If $\varphi:D\to \R$ is a smooth function, we adopt, respectively, the notations $d\varphi$, $\nabla \varphi$, to indicate the differential and the gradient of $\varphi$, which is the only vector field on $D$ such that $$d\varphi [X] = \l X, \nabla \varphi\r,  \ \ \hbox{for any} \ X \in \mathcal{T}(D).$$
 We will also use sometimes the notation $\nabla_{\S^{N-1}} \varphi$ instead of $\nabla \varphi$ to make a distinction with respect to the usual gradient of real valued functions defined in open subsets of $\RN$. The second covariant derivative of $\varphi$ is defined as
\beq\label{def:seccivder}
 \nabla_{X,Y} \varphi:= \nabla_X \nabla_Y \varphi - \nabla_{{\nabla_X Y}} \varphi, \ \ \hbox{for any} \ X,Y \in \mathcal{T}(D),
\eeq
and the Hessian of $\varphi$, denoted by $\nabla^2 \varphi$ or by $D^2 \varphi$, is the symmetric 2-tensor given by
$$\nabla^2 \varphi\ (X,Y):=\nabla_{X,Y} \varphi, \ \ \hbox{for any} \ X,Y \in \mathcal{T}(D).$$
The Laplacian of $\varphi$, denoted by $\Delta \varphi$, is the trace of the Hessian. Again, when there is a chance of confusion with the standard Laplacian we will use the notation $\Delta_{\S^{N-1}} \varphi$ instead of $\Delta \varphi$.\\

Let $\{e_1,\ldots,e_{N-1}\}$ be a local orthonormal frame field for $D$. For any $i,j\in\{1,\ldots,N-1\}$ we define the connection form $\omega_{ij}$ as

\begin{equation}\label{defconnform}
\omega_{ij}(X):=\l \nabla_X e_j , e_i\r,\ X \in \T(D).
\end{equation}
We recall that the connection forms are skew symmetric and in terms of the $\omega_{ij}$'s we can write 
\begin{equation}\label{dercoveiej}
 \nabla_{e_i}e_j = \sum_{k=1}^{N-1} \omega_{kj}(e_i) e_k.
\end{equation}
We denote by $\varphi_i$ the covariant derivative $\nabla_{e_i}\varphi$, and we recall that, by definition, $\nabla_{e_i} \varphi=d\varphi[e_i]$. It is easy to check that the gradient of $\varphi$ can be written as
$$\nabla \varphi= \sum_{i=1}^{N-1} \varphi_i e_i.$$
Finally, taking $X=e_i$, $Y=e_j$ in \eqref{def:seccivder} and using \eqref{defconnform} we have 
\begin{equation}\label{covarderij}
\nabla_{e_i,e_j} \varphi= \nabla_{e_i}\varphi_j -\sum_{k=1}^{N-1}\omega_{kj}(e_i) \varphi_k.
\end{equation}

From now on we will use the notation $\varphi_{ij}$ to denote $\nabla_{e_i,e_j} \varphi$. In particular the Laplacian of $\varphi$ can be written as $ \Delta\varphi= \sum_{i=1}^{N-1} \varphi_{ii}$.\\



Now we consider the special case of radial graphs.
\begin{definition}\label{def:radialgraphassociated}
Let $D \subset \S^{N-1}$ be a domain and let $\varphi \in C^2(D)$. We denote by $\Gamma_\varphi$
 the associated radial graph to $\varphi$, namely
$$
 \Gamma_\varphi:=\{x\in \R^N; \ x=e^{\varphi(q)}q, \ q\in D\}.
$$
\end{definition}
Clearly $\Gamma_\varphi$ is a $(N-1)$-dimensional manifold (of class $C^2$). We consider the map $\mathcal{Y}:D \to \Gamma_\varphi$ defined by 
\beq\label{eq:defparstand}
\mathcal{Y}(q):=e^{\varphi(q)}q, \ \ \ q\in D.
\eeq
For any fixed $q\in D$, let $\gamma_i:(-\delta,\delta) \to D$, $\gamma_i=\gamma_i(t)$ be a curve contained in $D$ and such that $\gamma_i(0)=q$, $\gamma^\prime_i(0)=e_i(q)$, for $i=1,\ldots,N-1$. Since
\begin{equation}\label{eq1covder}
\left.\frac{d(\mathcal{Y}\circ \gamma_i)}{dt}\right|_{t=0} = e^\varphi (\varphi_i q + e_i)
\end{equation}
then a local basis for $T_{\mathcal{Y}(q)}\Gamma_\varphi$ is given by
$$ E_i(q)= e^\varphi (e_i + \varphi_i q), \ \ i=1,\ldots,N-1,$$
and the components of the induced metric are 
$$g_{ij}=\l E_i, E_j\r=e^{2\varphi} (\l e_i, e_j \r + \varphi_i \varphi_j \l q, q\r )=e^{2\varphi} (\delta_{ij} + \varphi_i\varphi_j).$$
We denote by $\nu(\mathcal{Y}(q))$ the exterior unit normal at $\mathcal{Y}(q)\in \Gamma_\varphi$. It is easy to check that
\beq\label{eq:Gauss}
\nu(\mathcal{Y}(q))= \frac{q- \sum_{i=1}^{N-1} \varphi_i e_i }{(1+|\nabla \varphi|^2)^{1/2}}=  \frac{q- \nabla \varphi }{(1+|\nabla \varphi|^2)^{1/2}}.
\eeq
In addition by direct computation we see that the coefficients of the second fundamental form are
$$\mathbf{II}_{ij}
=\frac{e^{\varphi}\left(\delta_{ij} + \varphi_i\varphi_j -\varphi_{ij}\right)}{(1+|\nabla \varphi|^2)^{1/2}},$$
for any $i,j=1,\ldots,N-1$ (see \cite{Lo03} or \cite{BI2} for more details).
Finally, since the mean curvature at $\mathcal{Y}(q) \in \Gamma_\varphi$ is given by 
$$N H(\mathcal{Y}(q))=  \sum_{i,j=1}^{N-1} g^{ij} \mathbf{II}_{ij},$$
where $(g^{ij})$ is the inverse matrix of $(g_{ij})$, namely
\begin{equation}\label{invmatrixmetric}
g^{ij}=e^{-2\varphi}\left(\delta_{ij}-\frac{\varphi_i\varphi_j}{1+|\nabla \varphi|^2}\right),
\end{equation}
then, by a straightforward computation we see that $\varphi$ must satisfy the following equation
\beq\label{eq:meancurvaturerg}
\sum_{i,j=1}^{N-1} \left((1+|\nabla \varphi|^2)\delta_{ij} - \varphi_i \varphi_j)\right)\varphi_{ij}= (N-1)(1+|\nabla \varphi|^2) - (N-1) e^\varphi (1+|\nabla \varphi|^2)^{3/2}H(\mathcal{Y}(q)).
\eeq
Writing \eqref{eq:meancurvaturerg} in divergence form we obtain the well known equation for radial graphs of prescribed mean curvature (see \cite{Lo03} or \cite{TW})

\begin{equation}\label{eq:divformt}
\displaystyle -\mathrm{div}_{\S^{N-1}}\left(\frac{\nabla \varphi}{\sqrt{1+|\nabla \varphi|^{2}}}\right)+\frac{N-1}{\sqrt{1+|\nabla \varphi|^{2}}}=(N-1) e^{\varphi}H(e^{\varphi}q)\ \ \text{in }D.
\end{equation}

\section{Torsional energy for domains in cones}
In this section we define and study the torsional energy for smooth domains in cones and then we focus on the class of strictly star-shaped domains.

Let $D$ be a smooth proper domain of $\S^{N-1}$ and let $\Sigma_D$ be the cone spanned by $D$. For a bounded domain $\Omega\subset\Sigma_D$ we set:
$$\Gamma_\Omega:=\partial\Omega \subset \Sigma_D, \ \ \ \Gamma_{1,\Omega}:=\partial\Omega\cap\partial\Sigma_D,$$
and assume that $\H_{N-1}(\Gamma_{1,\Omega})>0$ and that $\Gamma_\Omega$ is a smooth $(N-1)$-dimensional manifold whose boundary $\partial\Gamma_\Omega=\partial\Gamma_{1,\Omega}\subset \partial\Sigma_D\setminus\{0\}$ is a smooth $(N-2)$-dimensional manifold. The set $\Gamma_\Omega$ is often called the relative (to $\Sigma_D$) boundary of $\Omega$.

We consider the following  mixed boundary value problem:
\begin{equation}\label{eq:mixbvprob}
\begin{cases}
-\Delta u = 1 & \text{in $\Omega$,}\\
 u = 0 & \text{on $\Gamma_\Omega$,}\\
 \frac{\partial u}{\partial \nu} = 0 & \text{on $\Gamma_{1,\Omega}\setminus\{0\}$.}
\end{cases}
\end{equation} 
It is easy to see that \eqref{eq:mixbvprob} admits a unique weak solution $u_\Omega$ in the space $H_0^{1}(\Omega\cup \Gamma_{1,\Omega})$ which is the Sobolev space of functions in $H^{1}(\Omega)$ whose trace vanishes on $\Gamma_\Omega$. Indeed $u_\Omega$ is the only minimizer of the functional
$$J(v):= \frac{1}{2}\int_\Omega |\nabla v|^2 \ dx - \int_\Omega v \ dx$$
in the space $H_0^{1}(\Omega\cup \Gamma_{1,\Omega})$ and we remark that $u_\Omega>0$ a.e. in $\Omega$, by the maximum principle (we refer to \cite{PT, PT2} for more details). 

Usually, the function $u_\Omega$ is called \emph{torsion function} of $\Omega$ and its energy $J(u_\Omega)$ represents the \emph{torsional energy} of the domain $\Omega$. This allows to consider the functional 
$$\E(\Omega; \Sigma_D)=J(u_\Omega)$$
which is defined on the domains contained in $\Sigma_D$. 

From the weak formulation of \eqref{eq:mixbvprob} we have
$$\int_\Omega |\nabla u_\Omega|^2 \ dx = \int_\Omega u_\Omega \ dx,$$
 which implies that
 \beq\label{eq:torsionminimizer}
 \E(\Omega;\Sigma_D)=-\frac{1}{2} \int_\Omega |\nabla u_\Omega|^2 \ dx = -\frac{1}{2}\int_\Omega u_\Omega \ dx. 
 \eeq
Now we focus on the special case when $\Omega$ is strictly star-shaped with respect to the origin which is the vertex of the cone $\Sigma_D$. Thus we consider the relative boundary $\Gamma_\Omega$ as the radial graph in $\Sigma_D$ of a function $\varphi \in C^2(\overline D, \R)$ as defined in Sect. 2. Therefore we denote $\Omega$ by $\Omega_\varphi$ which can be described as:
\begin{equation}\label{def:omegaphi}
\Omega_\varphi:=\{x\in \Sigma_D; \ x=s q, \ 0<s<e^{\varphi(q)}, \ q\in D\}.
\end{equation}
We restrict the torsional energy functional $\E$ to this class of domains and we denote it by $\mathscr{E}$, i.e. we set $$\mathscr{E}(\varphi):=\E(\Omega_\varphi; \Sigma_D).$$ 
We observe that $\mathscr{E}$ is a functional defined on $C^2(\overline{D},\R)$ and we compute its first and second derivatives. To this aim we point out that taking variations of $\varphi$ in $C^2(\overline{D},\R)$ corresponds to taking variations of $\Omega_\varphi$ in the class of strictly star-shaped domains (of class $C^2$).

Let us set for simplicity $$\Gamma_\varphi:=\Gamma_{\Omega_\varphi}, \ \ \ \Gamma_{1,\varphi}:=\Gamma_{1,\Omega_\varphi}.$$

If $v \in C^2(\overline D, \R)$ and $t \in (-\delta, \delta)$, where $\delta>0$ is a fixed number, we consider the domain variations $\Omega_{\varphi+tv} \subset \Sigma_D$, ${t\in(-\delta,\delta)}$. Let $\xi:(-\delta,\delta)\times \Sigma_D \to \Sigma_D$ be the map defined by
$$
\xi(t,x)=e^{tv\left(\frac{x}{|x|}\right)}x.
$$
It is elementary to check that, for a fixed $t\in (-\delta,\delta)$ the restriction 
\beq\label{eq:defxidiff}
\xi|_{\Omega_\varphi}(t,\cdot):\Omega_\varphi \to \Omega_{\varphi+tv}
\eeq
 is a diffeomorphism whose inverse $\left(\xi|_{\Omega_\varphi}\right)^{-1}:\Omega_{\varphi+tv} \to \Omega_\varphi$ is given by  $$\left(\xi|_{\Omega_\varphi}\right)^{-1}(x)=e^{-tv(\frac{x}{|x|})}x=\xi(-t,x).$$ Moreover by definition we have $\xi(t,x) \in \partial \Sigma_D\setminus\{0\}$ for all $(t,x) \in (-\delta,\delta) \times  \partial \Sigma_D\setminus\{0\}$. In particular $\xi$ is the flow associated to the vector field $V$ on $\Sigma_D$ given by
\begin{equation}\label{defvfV}
V(x):=v\left(\frac{x}{|x|}\right) x,
\end{equation}
since $\xi(0,x)=x$ and $\frac{d\xi}{dt}(t,x)=e^{tv\left(\frac{x}{|x|}\right)} v\left(\frac{x}{|x|}\right)x=V(\xi(t,x))$, and $(\Omega_{\varphi+tv})_{t\in(-\delta,\delta)}$ is a deformation of $\Omega_\varphi$ associated to the vector field $V$ (see \cite[Definition 1.1]{LS}). We now compute the derivative of $\mathscr{E}$ with respect to a variation $v \in C^2(\overline D,\R)$.

\begin{lemma}\label{lem:firstvar}
Let $\varphi \in C^2(\overline D,\R)$ and assume that $u_{\Omega_\varphi}\in W^{1,\infty}(\Omega_\varphi)\cap W^{2,2}(\Omega_\varphi)$. Then for any $v \in C^2(\overline D,\R)$ it holds
$$\mathscr{E}^\prime(\varphi)[v]=-\frac{1}{2} \int_D e^{N\varphi}\ v \left(\frac{\partial u_{\Omega_\varphi} }{\partial \nu}\left(e^{\varphi}q\right)\right)^2 \ d\sigma, $$
where $d\sigma$ is the $(N-1)$-dimensional area element of $\S^{N-1}$.
\end{lemma}
\begin{proof}
Let $\varphi \in C^2(\overline D,\R)$ as in the statement and let $v \in C^2(\overline D,\R)$.
By definition we have
\begin{equation}\label{eq:torsionalvart}
\mathscr{E}(\varphi+tv)=\E(\Omega_{\varphi+tv};\Sigma_D)= -\frac{1}{2}\int_{\Omega_{\varphi+tv}} u_{\Omega_{\varphi+tv}} \ dx,
\end{equation}
where $u_{\Omega_{\varphi+tv}}$ is the only (positive) weak solution to
\begin{equation}\label{eq:mixbvprob2}
\begin{cases}
-\Delta u = 1 & \text{in ${\Omega_{\varphi+tv}}$,}\\
 u = 0 & \text{on $\Gamma_{\varphi+tv}$,}\\
 \frac{\partial u}{\partial \nu} = 0 & \text{on $\Gamma_{1,\varphi+tv}\setminus\{0\}$.}
\end{cases}
\end{equation} 
Writing \eqref{eq:torsionalvart} in polar coordinates we obtain that
$$\mathscr{E}(\varphi+tv)= -\frac{1}{2}\int_D\int_0^{e^{\varphi+tv}} \rho^{N-1} u_{\Omega_{\varphi+tv}}(\rho q) \ d\rho d\sigma.$$
Let  $\hat \Phi:(-\delta,\delta)\to H^{1}_0(\Omega_\varphi\cup \Gamma_{1,\varphi})$ be the map defined by $$\hat \Phi(t):=\hat u_t,$$ 
where $\hat u_t:=u_{\Omega_{\varphi+tv}}\circ \xi(t,\cdot)|_{\Omega_\varphi} \in H^{1}_0(\Omega_\varphi\cup \Gamma_{1,\varphi})$, $\xi|_{\Omega_\varphi}(t,\cdot):\Omega_\varphi \to \Omega_{\varphi+tv}$ is the diffeomorphism given by \eqref{eq:defxidiff}. From the proof of \cite[Proposition 4.3]{PT2} we know that  $\hat \Phi$ is differentiable and thus we infer that $u_{\Omega_{\varphi+tv}}$ is differentiable with respect to $t$. Hence, by the Leibniz integral rule for differentiation of integral functions we get that
$$
\begin{array}{lll}
\displaystyle \frac{d}{dt}\left(\mathscr{E}(\varphi+tv)\right)&=&\displaystyle -\frac{1}{2}\int_D e^{(N-1)(\varphi+tv)} e^{\varphi+tv} v\ u_{\Omega_{\varphi+tv}}(e^{\varphi+tv} q) \ d\sigma\\[12pt]
&& \displaystyle  -\frac{1}{2}\int_D \int_0^{e^{\varphi+tv}} \rho^{N-1} \frac{d}{dt}\left(u_{\Omega_{\varphi+tv}}\right)(\rho q) \ d\rho d\sigma
\end{array}
$$
In view of \eqref{eq:mixbvprob2} we have $u_{\Omega_{\varphi+tv}}(e^{\varphi+tv} q)=0$ on $D$ for any $t\in(-\delta, \delta)$. In particular computing at $t=0$ we have
\begin{equation}\label{eq:firstvar1}
\begin{array}{lll}
\displaystyle \mathscr{E}^\prime(\varphi)[v]=  \frac{d}{dt}\left.\left(\mathscr{E}(\varphi+tv)\right)\right|_{t=0}&=&\displaystyle  -\frac{1}{2}\int_D \int_0^{e^{\varphi(q)}} \rho^{N-1} \frac{d}{dt}\left.\left(u_{\Omega_{\varphi+tv}}\right)\right|_{t=0}(\rho q) \ d\rho d\sigma\\[12pt]
&=&\displaystyle  -\frac{1}{2}\int_{\Omega_\varphi} \frac{d}{dt}\left.\left(u_{\Omega_{\varphi+tv}}\right)\right|_{t=0} \ dx.
\end{array}
\end{equation}
Setting $u^\prime:=\left.\frac{d}{dt}\left(u_{\Omega_{\varphi+tv}}\right)\right|_{t=0}$ and arguing as in the proof of \cite[Proposition 4.3]{PT2}, where the assumption $u_{\Omega_\varphi}\in W^{1,\infty}(\Omega_\varphi)\cap W^{2,2}(\Omega_\varphi)$ is used, we infer that $u^\prime \in H^{1}_0(\Omega_\varphi\cup\Gamma_{1,\varphi})$ satisfies
\begin{equation}\label{eq:mixbvprobuprime}
\begin{cases}
-\Delta u^\prime = 0 & \text{in $\Omega_{\varphi}$,}\\
 u^\prime = -\frac{\partial u_{\Omega_\varphi}}{\partial\nu} \langle V,\nu\rangle & \text{on $\Gamma_{\varphi}$,}\\
 \frac{\partial u^\prime}{\partial \nu} = 0 & \text{on $\Gamma_{1,\varphi}\setminus\{0\}$.}
\end{cases}
\end{equation} 
In particular, in view of \eqref{defvfV} and since $\Gamma_{\varphi}$ is a radial graph we have $\nu(x)=\frac{\frac{x}{|x|}-\nabla_{\S^{N-1}}\varphi \left(\frac{x}{|x|}\right)}{\sqrt{1+|\nabla_{\S^{N-1}} \varphi\left(\frac{x}{|x|}\right)|^2}}$, for any $x\in \Gamma_\varphi$ (see \eqref{eq:Gauss}), and thus 
\begin{equation}\label{eq:Vnuradgraph}
\langle V,\nu\rangle=\frac{|x|}{\sqrt{1+\left|\nabla_{\S^{N-1}} \varphi\left(\frac{x}{|x|}\right)\right|^2}} v\left(\frac{x}{|x|}\right)\ \ \text{on $\Gamma_\varphi$}.
\end{equation}
Rewriting \eqref{eq:firstvar1} in terms of $u^\prime$, applying Green's second identity (which holds also in conic domains, since it is a consequence of the divergence theorem, see e.g. \cite[Lemma 2.1]{PT}) and taking into account \eqref{eq:mixbvprob} (with $\Omega=\Omega_\varphi$), \eqref{eq:mixbvprobuprime} and \eqref{eq:Vnuradgraph} we get that
\begin{equation}\label{eq:firstafinal}
\begin{array}{lll}
\displaystyle  \mathscr{E}^\prime(\varphi)[v]&=&\displaystyle  -\frac{1}{2}\int_{\Omega_\varphi} u^\prime \ dx.\\[12pt]
&=&\displaystyle  \frac{1}{2}\int_{\Omega_\varphi} u^\prime {\Delta u_{\Omega_\varphi}} \ dx -\frac{1}{2}\int_{\Omega_\varphi} \underbrace{\Delta u^\prime}_{=0} u_{\Omega_\varphi} \ dx\\[22pt]
&=&\displaystyle  \frac{1}{2}\int_{\Gamma_\varphi} u^\prime \frac{\partial u_{\Omega_\varphi}}{\partial\nu} \ d\sigma_{\Gamma_\varphi} + \frac{1}{2}\int_{\Gamma_{1,\varphi}\setminus\{0\}} u^\prime \underbrace{\frac{\partial u_{\Omega_\varphi}}{\partial\nu}}_{=0}  \ d\sigma_{\Gamma_{1,\varphi}\setminus\{0\}}\\[24pt]
&=&\displaystyle - \frac{1}{2}\int_{\Gamma_\varphi} \left(\frac{\partial u_{\Omega_\varphi}}{\partial\nu}(x)\right)^2 \frac{|x|}{\sqrt{1+\left|\nabla_{\S^{N-1}} \varphi\left(\frac{x}{|x|}\right)\right|^2}} v\left(\frac{x}{|x|}\right)  \ d\sigma_{\Gamma_\varphi},
\end{array}
\end{equation}
where $d\sigma_{\Gamma_\varphi}$, $d\sigma_{\Gamma_{1,\varphi}\setminus\{0\}}$ are the $(N-1)$-dimensional area elements of $\Gamma_\varphi$, $\Gamma_{1,\varphi}\setminus\{0\}$, respectively.
Finally, writing $x=e^{\varphi(q)}q$, $q\in D$, observing that $d\sigma_{\Gamma_\varphi}=e^{(N-1)\varphi}\sqrt{1+|\nabla_{\S^{N-1}}\varphi|^2} d\sigma$ and $\frac{x}{|x|}=q$, then from \eqref{eq:firstafinal} we obtain that
$$ \mathscr{E}^\prime(\varphi)[v] = - \frac{1}{2}\int_{D} e^{N\varphi} v\left(\frac{\partial u_{\Omega_\varphi}}{\partial\nu}(e^{\varphi}q)\right)^2 \ d\sigma,$$
and this completes the proof.
\end{proof}
For the second variation of the functional $\mathscr{E}$ we have
\begin{lemma}\label{lem:secvarT}
Let $\varphi$ be as in Lemma \ref{lem:firstvar}. Then for any $v,w \in C^2(\overline D, \R)$ it holds
\begin{equation}\label{eq:secvarthesisgeneral}
\begin{array}{lll}
\displaystyle\mathscr{E}^{\prime\prime}(\varphi)[v,w]&=&\displaystyle - \frac{N}{2}\int_{D} e^{N\varphi}\,  v\, w\left(\frac{\partial u_{\Omega_\varphi}}{\partial\nu}(e^{\varphi}q)\right)^2 d\sigma\\[12pt]
&& \displaystyle - \int_{D} e^{N\varphi}\,  v\, \frac{\partial u_{\Omega_\varphi}}{\partial\nu} (e^{\varphi}q)\, \frac{\partial u^\prime_w}{\partial\nu} (e^{\varphi}q) \ d\sigma\\[12pt]
&&\displaystyle - \int_{D} e^{N\varphi}  v\, w\, \frac{\partial u_{\Omega_\varphi}}{\partial\nu}(e^{\varphi}q)\,  D^2u_{\Omega_\varphi}(e^{\varphi(q)}q) e^{\varphi}q  \boldsymbol{\cdot} \nu \ d\sigma\\[12pt]
&&\displaystyle + \int_{D} e^{N\varphi}  v\, \frac{\partial u_{\Omega_\varphi}}{\partial\nu}(e^{\varphi}q)\, \frac{\nabla u_{\Omega_\varphi}(e^{\varphi}q)  \boldsymbol{\cdot} \nabla_{\S^{N-1}} w}{\sqrt{1+|\nabla_{\S^{N-1}}\varphi|^2}} \ d\sigma\\[12pt]
&&\displaystyle +\int_{D} e^{N\varphi}  v\left(\frac{\partial u_{\Omega_\varphi}}{\partial\nu}(e^{\varphi}q)\right)^2 \frac{\nabla_{\S^{N-1}}\varphi  \boldsymbol{\cdot} \nabla_{\S^{N-1}} w}{(1+|\nabla_{\S^{N-1}}\varphi|^2)} \ d\sigma,
\end{array}
\end{equation}
where $u^\prime_w=\left.\frac{d}{ds}\left(u_{\Omega_{\varphi+sw}}\right)\right|_{s=0}$ is the solution to \eqref{eq:mixbvprobuprime} with $V$ given by $V(x)=w\left(\frac{x}{|x|}\right)x$.
\end{lemma}
\begin{proof}
Let us fix $v, w \in C^2(\overline D, \R)$, by definition and by Lemma \ref{lem:firstvar} we have
\begin{equation}\label{eq:secvardef}
\begin{array}{lll}
\displaystyle \mathscr{E}^{\prime\prime}(\varphi)[v,w] &=& \displaystyle\left.\frac{d}{ds}\left( - \frac{1}{2}\int_{D} e^{N(\varphi+sw)} v\left(\frac{\partial u_{\Omega_{\varphi+sw}}}{\partial\nu}(e^{\varphi+sw}q)\right)^2 \ d\sigma \right)\right|_{s=0},\
\end{array}
\end{equation}
and thus 
\begin{equation}\label{eq:secvarfinal}
\begin{array}{lll}
\displaystyle \mathscr{E}^{\prime\prime}(\varphi)[v,w] &=&\displaystyle - \frac{1}{2}\int_{D} e^{N\varphi}N  v\, w\, \frac{\partial u_{\Omega_\varphi}}{\partial\nu}(e^{\varphi}q)\, \ d\sigma\\[12pt]
 &&\displaystyle - \int_{D} e^{N\varphi} v\, \frac{\partial u_{\Omega_{\varphi}}}{\partial\nu}(e^{\varphi}q)\left. \frac{d}{ds}\left(\frac{\partial u_{\Omega_{\varphi+sw}}}{\partial\nu}(e^{\varphi+sw}q)\right)\right|_{s=0} \ d\sigma.
\end{array}
\end{equation}
Since $\Gamma_{\varphi+sw}$ is a radial graph, then, in view of \eqref{eq:Gauss}, we have
\begin{equation}\label{eq:secvarpartialnormal}
\begin{array}{lll}
\displaystyle\frac{\partial u_{\Omega_{\varphi+sw}}}{\partial\nu}(e^{\varphi+sw}q)
&=&\displaystyle\nabla u_{\Omega_{\varphi+sw}}(e^{\varphi+sw}q)  \boldsymbol{\cdot} \frac{q-\nabla_{\S^{N-1}}(\varphi+sw)}{\sqrt{1+|\nabla_{\S^{N-1}} (\varphi+sw)|^2}}.
\end{array}
\end{equation}
As in the proof of Lemma \ref{lem:firstvar} we consider the map $\hat \Phi:(-\delta,\delta)\to H^{1}_0(\Omega_\varphi\cup \Gamma_{1,\varphi})$, defined by $$\hat \Phi(s)= \hat u_s:=u_{\Omega_{\varphi+sw}}\circ \xi(s,\cdot)|_{\Omega_\varphi}.$$
Moroever, let $G:H^{1}_0(\Omega_\varphi\cup \Gamma_{1,\varphi}) \to L^2(\Omega_\varphi\cup \Gamma_{1,\varphi}, \R^N)$, given by $G(f):=\nabla f$. Since $G$ is a bounded linear operator, then $G$ is differentiable, $G^\prime(f)[g]=\nabla g$ for any $g \in H^{1}_0(\Omega_\varphi\cup \Gamma_{1,\varphi})$. In addition, as $\hat \Phi$ is differentiable (see the proof of \cite[Proposition 4.3]{PT2} for the details), then the composition $G\circ \hat \Phi: (-\delta,\delta)\to L^2(\Omega_\varphi\cup \Gamma_{1,\varphi}, \R^N)$ is differentiable and $$(G\circ \hat \Phi)^\prime(s)=G^\prime(\hat \Phi(s))[\hat \Phi^\prime(s)]=\nabla\hat \Phi^\prime(s) \ \ \forall s\in(-\delta,\delta).$$ 
In terms of $\hat u_s$ this means that
\beq\label{eq:eqdsnablahatus}
\frac{d}{ds}\left(\nabla \hat u_s\right)= \nabla \left(\frac{d\hat u_s}{ds}\right).
\eeq
In addition, since $u_{\Omega_{\varphi+sw}}=\hat u_s \circ \xi(-s,\cdot)|_{\Omega_{\varphi+sw}}$ it follows that also $s\mapsto \nabla u_{\Omega_{\varphi+sw}}$ is differentiable. We claim that 
\beq\label{eq:claimSect3}
\frac{d}{ds}\left(\nabla u_{\Omega_{\varphi+sw}}\right)=\nabla\left(\frac{d}{ds}  u_{\Omega_{\varphi+sw}}\right). 
\eeq
Indeed, setting $\xi_s:=\xi(-s,\cdot)|_{\Omega_{\varphi+sw}}$, since $$\frac{\partial}{\partial x_i} u_{\Omega_{\varphi+sw}} =\frac{\partial}{\partial x_i}\left(\hat u_s \circ \xi_s\right)=\nabla \hat u_s(\xi_s) \boldsymbol{\cdot} \frac{\partial \xi_s}{\partial x_i},$$
then, using \eqref{eq:eqdsnablahatus}, we get that  
 $$\frac{d}{ds}\left(\frac{\partial}{\partial x_i}u_{\Omega_{\varphi+sw}}\right)=\left(\nabla\left( \frac{d \hat u_s}{ds}\right) (\xi_s)+ D^2\hat u_s(\xi_s) \frac{d\xi_s}{ds}\right) \boldsymbol{\cdot} \frac{\partial \xi_s}{\partial x_i}+\nabla \hat u_s(\xi_s) \boldsymbol{\cdot} \frac{d}{ds}\left(\frac{\partial \xi_s}{\partial x_i}\right),$$ 
for any $i=1,\ldots,N$. On the other hand, since $\xi:(-\delta,\delta)\times\Sigma_D\to \Sigma_D$ is smooth we deduce that $$\frac{\partial}{\partial x_i}\frac{d\xi_s}{ds}= \frac{d}{ds}\frac{\partial \xi_s}{\partial x_i}$$ 
 and thus by a straightforward computation we obtain 
$$\frac{d}{ds}\left(\frac{\partial}{\partial x_i}u_{\Omega_{\varphi+sw}}\right)= \left( \nabla\left( \frac{d \hat u_s}{ds}\right) (\xi_s)+ D^2\hat u_s(\xi_s) \frac{d\xi_s}{ds}\right) \boldsymbol{\cdot} \frac{\partial \xi_s}{\partial x_i}+\nabla \hat u_s(\xi_s) \boldsymbol{\cdot}\frac{\partial}{\partial x_i}\frac{d\xi_s}{ds}=\frac{\partial}{\partial x_i}\left(\frac{d}{ds}  u_{\Omega_{\varphi+sw}}\right),$$
for any $i=1,\ldots,N$, which proves Claim \eqref{eq:claimSect3}.\\

Thanks to \eqref{eq:secvarpartialnormal} and \eqref{eq:claimSect3} we have

\begin{equation}\label{eq:secvar1}
\begin{array}{lll}
&&\displaystyle \frac{d}{ds}\left.\left(\frac{\partial u_{\Omega_{\varphi+sw}}}{\partial\nu}(e^{\varphi+sw}q)\right)\right|_{s=0}\\[12pt]
&=&\displaystyle \left(\nabla u^\prime_w(e^{\varphi}q)+D^2 u_{\Omega_\varphi}(e^{\varphi}q)e^{\varphi}w q\right)  \boldsymbol{\cdot} \frac{q-\nabla_{\S^{N-1}}\varphi}{\sqrt{1+|\nabla_{\S^{N-1}} \varphi|^2}}\\[12pt]
&&\displaystyle + \nabla u_{\Omega_{\varphi}}(e^{\varphi}q)  \boldsymbol{\cdot} \left( - \frac{\nabla_{\S^{N-1}}w}{\sqrt{1+|\nabla_{\S^{N-1}} \varphi|^2}} - \frac{(q-\nabla_{\S^{N-1}}\varphi) (\nabla_{\S^{N-1}} \varphi   \boldsymbol{\cdot} \nabla_{\S^{N-1}}w)}{(1+|\nabla_{\S^{N-1}} \varphi|^2)^{3/2}}\right)\\[16pt]
&=& \left(\nabla u^\prime_w(e^{\varphi}q)+D^2 u_{\Omega_\varphi}(e^{\varphi}q)e^{\varphi}w q\right)  \boldsymbol{\cdot} \nu\\[12pt]
&&\displaystyle + \nabla u_{\Omega_{\varphi}}(e^{\varphi}q)  \boldsymbol{\cdot} \left(-\frac{\nabla_{\S^{N-1}}w}{\sqrt{1+|\nabla_{\S^{N-1}} \varphi|^2}} -\frac{\nabla_{\S^{N-1}} \varphi   \boldsymbol{\cdot} \nabla_{\S^{N-1}}w}{(1+|\nabla_{\S^{N-1}} \varphi|^2)}\, \nu\right)\\[18pt]
&=&\displaystyle \frac{\partial u^\prime_w}{\partial \nu}(e^{\varphi}q)+ w\, D^2 u_{\Omega_\varphi}(e^{\varphi}q)e^{\varphi}q  \boldsymbol{\cdot} \nu\\[12pt]
&&\displaystyle  - \frac{ \nabla u_{\Omega_{\varphi}}(e^{\varphi}q) \boldsymbol{\cdot} \nabla_{\S^{N-1}}w}{\sqrt{1+|\nabla_{\S^{N-1}} \varphi|^2}} - \frac{\partial u_{\Omega_\varphi}}{\partial \nu}(e^{\varphi}q)\frac{\nabla_{\S^{N-1}} \varphi   \boldsymbol{\cdot} \nabla_{\S^{N-1}}w}{(1+|\nabla_{\S^{N-1}} \varphi|^2)}.
\end{array}
\end{equation}
Finally, combining \eqref{eq:secvarfinal} and \eqref{eq:secvar1}  we readily obtain \eqref{eq:secvarthesisgeneral}. The proof is complete.
\end{proof}
\section{Volume-constrained critical points for the torsional energy of star-shaped domains}
For any $\varphi \in C^2(\overline D,\R)$, the volume of the associated star-shaped domain $\Omega_\varphi$ (see \eqref{def:omegaphi}) is given by  
\beq\label{eq:defVolOmegaphi}
\mathcal{V}(\varphi)=|\Omega_\varphi|=\frac{1}{N}\int_D e^{N\varphi} \ d\sigma,
\eeq
where $d\sigma$ is the $(N-1)$-dimensional area element of $\S^{N-1}$.
It is easy to check that $\V$ is of class $C^2$ and for any $v,w \in C^2(\overline D,\R)$ it holds
\begin{equation}\label{eq:firstvarvol}
\V^\prime(\varphi)[v]= \int_D e^{N\varphi}\, v \ d\sigma,
\end{equation}
and
\begin{equation}\label{eq:secvarvol}
\V^{\prime\prime}(\varphi)[v,w]= N \int_D e^{N\varphi}\,v\, w \ d\sigma.
\end{equation}
For a number $c>0$ we define:
\beq\label{def:M}
M=\{\varphi \in C^2(\overline D, \R); \ \V(\varphi)=c\}.
\eeq
Clearly $M$ is a smooth manifold and for any $\varphi\in M$ it holds
\beq\label{def:TvarphiM} 
T_\varphi M=\left\{v \in C^2(\overline D, \R); \ \int_D e^{N\varphi}\, v \ d\sigma = 0\right\}.
\eeq
 We consider the restriction of the torsional energy to the domains corresponding to functions $\varphi\in M$, namely  the functional  defined by
\beq\label{eq:defI}
 I(\varphi):=\left.\mathscr{E}(\varphi)\right|_{\varphi \in M}=\left.\E(\Omega_\varphi; \Sigma_D)\right|_{\varphi \in M}. 
 \eeq
 If $\varphi\in M$ is critical point  of $I$ then there exists $\lambda \in \R$ such that
\begin{equation}\label{eq:LagMult}
\mathscr{E}^\prime(\varphi) = \lambda \V^\prime(\varphi).
\end{equation}
As a straightforward consequence of Lemma \ref{lem:firstvar} and \eqref{eq:firstvarvol} we have
\begin{lemma}\label{lem:lagmult}
Let $\varphi \in M$ be a critical point for $I$ and assume that $u_{\Omega_\varphi}\in W^{1,\infty}(\Omega_\varphi)\cap W^{2,2}(\Omega_\varphi)$. Then the Lagrange multiplier $\lambda$ is negative and $$\frac{\partial u_{\Omega_\varphi}}{\partial \nu}\equiv- \sqrt{-2\lambda}\ \ \text{on $\Gamma_\varphi$}.$$
\end{lemma}
\begin{proof}
Let $\varphi \in M$ be a critical point for $I$ and assume that $u_{\Omega_\varphi}\in W^{1,\infty}(\Omega_\varphi)\cap W^{2,2}(\Omega_\varphi)$, then, from \eqref{eq:LagMult} and exploiting Lemma \ref{lem:firstvar} and \eqref{eq:firstvarvol}, we have
$$-\frac{1}{2} \int_D e^{N\varphi}\, v \left(\frac{\partial u_{\Omega_\varphi} }{\partial \nu}\left(e^{\varphi}q\right)\right)^2 \ d\sigma = \lambda \int_D e^{N\varphi}\, v \ d\sigma,$$
for any $v \in C^2(\overline D, \R)$. Hence we readily obtain that
$$\int_D e^{N\varphi}\, v\left[ \left(\frac{\partial u_{\Omega_\varphi} }{\partial \nu}\left(e^{\varphi}q\right)\right)^2 +2 \lambda \right]\ d\sigma=0,$$
and from the arbitrariness of $v \in C^2(\overline D, \R)$ we easily deduce that $\lambda<0$ and 
 \begin{equation}\label{eq:lem1vccp}
 \left(\frac{\partial u_{\Omega_\varphi} }{\partial \nu}\right)^2 =-2 \lambda \ \text{on $\Gamma_\varphi$.}
 \end{equation}
 Now, recalling that $u_{\Omega_\varphi}$ is the only (positive) weak solution to
\begin{equation}\label{eq:mixbvprobvarphi}
\begin{cases}
-\Delta u = 1 & \text{in $\Omega_\varphi$,}\\
 u = 0 & \text{on $\Gamma_\varphi$,}\\
 \frac{\partial u}{\partial \nu} = 0 & \text{on $\Gamma_{1,\varphi}\setminus\{0\}$,}
\end{cases}
\end{equation} 
then from standard regularity estimates we infer that $u_{\Omega_\varphi}$ is smooth in $\Omega_\varphi$, and from Hopf's lemma we get that $\frac{\partial u_{\Omega_\varphi} }{\partial \nu}<0$ on $\Gamma_\varphi$. 
Hence, in view of \eqref{eq:lem1vccp}, we obtain
$$\frac{\partial u_{\Omega_\varphi}}{\partial \nu} =-\sqrt{-2 \lambda} \ \text{on $\Gamma_\varphi$}.$$
\end{proof}
\begin{remark}
From Lemma \ref{lem:lagmult} we deduce that each critical point of $I$ produces a star-shaped domain $\Omega_\varphi$ for which the overdetermined problem \eqref{eq:overdetermined} has a solution. We recall that, as shown in \cite[Proposition 4.3]{PT2}, each critical point of the functional $\E(\Omega; \Sigma_D)$ on the whole family of domains in $\Sigma_D$, with a volume constraint, is a domain for which \eqref{eq:overdetermined} has a solution. Hence Lemma \ref{lem:lagmult} shows that the same statement holds even if the variations are taken only in the class of star-shaped domains.
\end{remark}
In the next result we compute the second derivative of $I$ at critical point along variations in $T_\varphi M$.
\begin{lemma}\label{lem:secvarI}
Let $\varphi \in M$ be a critical point for $I$ and let $v,w \in T_\varphi M$. Then
$$I^{\prime\prime}(\varphi)[v,w]=\mathscr{E}^{\prime\prime}(\varphi)[v,w] - \lambda \V^{\prime\prime}(\varphi)[v,w],$$
where $\lambda$ is the Lagrange multiplier.
\end{lemma} 
\begin{proof}
By definition if $\varphi \in M$ is a critical point for $I$, the second variation $I^{\prime\prime}(\varphi)[v,w]$ along the variations $v,w \in T_\varphi M$ is given by
$$I^{\prime\prime}(\varphi)[v,w]=\left.\frac{\partial^2 I(\Psi(t,s))}{\partial t \partial s}\right|_{(t,s)=(0,0)},$$
where $\Psi:(-\epsilon,\epsilon)\times (-\epsilon,\epsilon)\to M$ is a smooth surface in $M$ such that
$$\Psi(0,0)=\varphi, \ \frac{\partial \Psi}{\partial t}(0,0)=v,  \ \frac{\partial \Psi}{\partial s}(0,0)=w.$$
We recall that by definition it holds $I(\Psi(t,s))= \mathscr{E}(\Psi(t,s))$. Since
$$\frac{\partial}{\partial s} \left(\mathscr{E}(\Psi(t,s))\right)=\mathscr{E}^\prime(\Psi(t,s))\left[\frac{\partial\Psi}{\partial s}(t,s)\right],$$
we have
\begin{equation}\label{eq:secvartorsionPsi}
\frac{\partial}{\partial t} \frac{\partial}{\partial s}\left(\mathscr{E}(\Psi(t,s))\right)=\mathscr{E}^{\prime\prime}(\Psi(t,s))\left[\frac{\partial\Psi}{\partial s}(t,s), \frac{\partial\Psi}{\partial t}(t,s)\right] + \mathscr{E}^\prime(\Psi(t,s))\left[\frac{\partial^2\Psi}{\partial t\partial s}(t,s)\right].
\end{equation}
On the other hand, since $\Psi(t,s) \in M$ we have $\V(\Psi(t,s))=c$ for any $(t,s)\in(-\epsilon,\epsilon)\times (-\epsilon,\epsilon)$, and thus differentiating with respect to $t$ we infer that $\V^\prime(\Psi(t,s))[\frac{\partial \Psi}{\partial s}(t,s)]=0$. Differentiating again with respect to $s$ we obtain
\begin{equation}\label{eq:secvarvolpsi}
\V^{\prime\prime}(\Psi(t,s))\left[\frac{\partial \Psi}{\partial s}(t,s), \frac{\partial \Psi}{\partial t}(t,s)\right]+\V^\prime(\Psi(t,s))\left[\frac{\partial}{\partial t}\frac{\partial \Psi}{\partial s}(t,s)\right]=0
\end{equation}
Hence, computing \eqref{eq:secvartorsionPsi} at $(t,s)=(0,0)$, since $\varphi=\Psi(0,0)$ is a critical point of $I$ and taking into account that \eqref{eq:LagMult}, \eqref{eq:secvarvolpsi}, we get that
\begin{eqnarray*}
\frac{\partial}{\partial t}\frac{\partial}{\partial s}\left. \left(\mathscr{E}(\Psi(t,s))\right)\right|_{(t,s)=(0,0)}&=&\mathscr{E}^{\prime\prime}(\Psi(t,s))\left[v,w\right] +\mathscr{E}^{\prime}(\Psi(t,s))\left[\frac{\partial}{\partial t}\frac{\partial\Psi}{\partial s}(0,0)\right]\\[6pt]
 &=&\mathscr{E}^{\prime\prime}(\Psi(t,s))\left[v,w\right] +\lambda \V^\prime(\varphi) \left[\frac{\partial}{\partial t}\frac{\partial\Psi}{\partial s}(0,0)\right]\\[6pt]
 &=&\mathscr{E}^{\prime\prime}(\Psi(t,s))\left[v,w\right] -\lambda \V^{\prime\prime}(\varphi) \left[v,w\right],
\end{eqnarray*}
which proves the desired relation.
\end{proof}

\begin{remark}\label{rem:varphi0}
When $\varphi\equiv 0$ then $\Omega_\varphi$ is the unit spherical sector $\Omega_{D,1}=\Sigma_D\cap B_1$, where $B_1=B_1(0)$ is the unit ball in $\RN$ centered at the origin. We denote it by $\Omega_0$, while $\Gamma_0$ will be its relative boundary. In this case the torsion function $u_{\Omega_0}$ is known to be the radial function $u_{\Omega_0}(x)=\frac{1-|x|^2}{2N}$. Then we can choose $c=|\Omega_0|$ in the definition of $M$ and the tangent space to $M$ at $\varphi\equiv 0$ is
 $T_0M=\{v \in C^2(\overline D, \R);\ \int_D v \ d\sigma=0 \}$. It is easy to check that $\nabla u_{\Omega_0}=-\frac{1}{N}x$, for $x \in \Sigma_D\cap B_1$, and  $\frac{\partial u_{\Omega_0}}{\partial \nu}=-\frac{1}{N}$ on $\Gamma_0$, so that $\Omega_0$ is a critical point for $I$ with  $\lambda=-\frac{1}{2N^2}$. Finally $D^2 u_{\Omega_0}(x)=-\frac{1}{N} \mathbb{I}_N$, for $x \in \Sigma_D\cap B_1$, where $\mathbb{I}_N$ is the identity matrix of order $N$, and thus we readily have that $u_{\Omega_0}\in W^{1,\infty}(\Omega_0)\cap W^{2,2}(\Omega_0)$. 
\end{remark}

For the second variation we have

\begin{proposition}\label{prop:secvartorzero}
For any $v \in T_0M$ it holds
\begin{equation}\label{eq:secvarthesis}
\displaystyle I^{\prime\prime}(0)[v,v]=\displaystyle - \frac{1}{N^2}\int_{D}   v^2 \ d\sigma \displaystyle +\frac{1}{N} \int_{D}  v \frac{\partial u^\prime}{\partial\nu} \ d\sigma,
\end{equation}
where  $u^\prime=\left.\frac{d}{dt}\left(u_{\Omega_{0+tv}}\right)\right|_{t=0}$ (see \eqref{eq:mixbvprobuprime}) and $\frac{\partial u^\prime}{\partial\nu}$ is the normal derivative of $u^\prime$ on $\Gamma_0=D$.
\end{proposition}
\begin{proof}
First we observe that, taking $\varphi\equiv 0$, from \eqref{eq:mixbvprobuprime} and \eqref{defvfV} we have that $u^\prime$ satisfies
\begin{equation}\label{eq:mixboundvalue2}
\begin{cases}
-\Delta u^\prime = 0 & \text{in $\Omega_0$},\\
 u^\prime = \frac{1}{N} v & \text{on $\Gamma_0$},\\
  \frac{\partial u^\prime}{\partial \nu} = 0 & \text{on $\Gamma_{1,0}$}.
\end{cases}
\end{equation}

Then, taking $v \in C^2(\overline D, \R)$ such that $\int_D v \ d\sigma=0$, from Lemma \ref{lem:secvarT}, Lemma \ref{lem:secvarI}, Remark \ref{rem:varphi0} and \eqref{eq:secvarvol},   \eqref{eq:mixboundvalue2} we obtain
\begin{equation}
\begin{array}{lll}
\displaystyle I^{\prime\prime}[v,v]&=&\displaystyle - \frac{N}{2}\int_{D}  v^2 \left(-\frac{1}{N}\right)^2 \ d\sigma - \int_{D} v \left(-\frac{1}{N}\right) \frac{\partial u^\prime}{\partial\nu} \ d\sigma\\[12pt]
&&\displaystyle  - \int_{D} v^2 \left(-\frac{1}{N}\right) \left( -\frac{1}{N}q \boldsymbol{\cdot} q \right) \ d\sigma + \int_{D} v \left(-\frac{1}{N}\right) \nabla u_{\Omega_0}  \boldsymbol{\cdot} \nabla_{\S^{N-1}} v\ d\sigma\\[12pt]
&&\displaystyle - \left(-\frac{1}{2N^2}\right) N \int_D v^2 \ d\sigma\\[12pt]
&=&\displaystyle \left(- \frac{1}{2N}-\frac{1}{N^2}+\frac{1}{2N}\right)\int_{D} v^2 \ d\sigma +\frac{1}{N} \int_{D}  v \frac{\partial u^\prime}{\partial\nu} \ d\sigma,
\end{array}
\end{equation}
since $\nabla u_{\Omega_0}  \boldsymbol{\cdot} \nabla_{\S^{N-1}} v\equiv0$ in $D$  because $\nabla u_{\Omega_0}$ is proportional to the radial direction.
\end{proof}
\begin{remark}
We observe that thanks to \eqref{eq:secvarthesis}, since $u^\prime=\frac{1}{N} v$ on $\Gamma_0$, by \eqref{eq:mixboundvalue2} and recalling that $\Gamma_0=D$, we can write
$$ I^{\prime\prime}(0)[v,v]=\displaystyle - \int_{D}   (u^\prime)^2 \ d\sigma + \int_{D}  u^\prime\frac{\partial u^\prime}{\partial\nu} \ d\sigma.$$
Then by Green's identity and \eqref{eq:mixboundvalue2} we infer that
$$ I^{\prime\prime}(0)[v,v]=\displaystyle - \int_{D}   (u^\prime)^2 \ d\sigma + \int_{\Omega_0}  |\nabla u^\prime|^2 \ dx.$$
\end{remark}
\section{A condition for instability}
In this section we provide conditions on the domain $D\subset\S^{N-1}$ such that the corresponding spherical sector (i.e. the domain $\Omega_0$ associated to the function $\varphi\equiv 0$, see \eqref{def:omegaphi}) is not a local minimizer for the torsional energy functional under a volume constraint. This is achieved by showing that $\Omega_0$ is an unstable critical point of $I$, i.e. its Morse index is positive.

More precisely, let $M$ be the manifold defined in \eqref{def:M}, with $c=|\Omega_0|$ and let $I$ be as in \eqref{eq:defI}. As observed in Remark \ref{rem:varphi0} the function $\varphi\equiv 0$ belongs to $M$ and $\Omega_0$ is a critical point for $I$. The main result of this section is the following.

\begin{theorem}\label{prop:critfirsteigenminim}
Let $D\subset\S^{N-1}$ be a smooth domain and let $\lambda_1(D)$ be the first non trivial eigenvalue of the Laplace-Beltrami operator $-\Delta_{\S^{N-1}}$, with zero Neumann condition on $\partial D$.
It holds:
\begin{itemize}
\item[(i)] if $\lambda_1(D) < N-1$, then $\Omega_0$ is not a local minimizer for $I$;
\item[(ii)] if $\lambda_1(D) > N-1$, then $\Omega_0$ is a local minimizer for $I$.
\end{itemize}
\end{theorem}
\begin{proof}
To prove (i), let $(w_j)_{j\in \N}$ be a $L^2(D)$-orthonormal basis of eigenfunctions of the eigenvalue problem
\begin{equation}\label{eq:eigenneumannsphere}
 \begin{cases} -\Delta_{\S^{N-1}} w_j=\lambda_j w_j & \hbox{in $D$},\\[8pt]\quad \quad \quad \frac{\partial w_j}{\partial \nu_{_{\partial D}}} =0 & \hbox{on $\partial D$},\end{cases}
 \end{equation}
 where $\nu_{_{\partial D}}$ is the exterior unit co-normal to $\partial D$, i.e. for any $q\in\partial D$, $\nu_{_{\partial D}}(q)$ is the only unit vector in $T_q\S^{N-1}$ such that $\nu_{_{\partial D}}(q)\perp T_q\partial D$ and $\nu_{_{\partial D}}(q)$ points outward $D$. We define the following extension of $w_j$ to the cone $\Sigma_D$
 \beq\label{eq:harmonicextension}
 \tilde w_j(rq):=\frac{1}{N}r^{\alpha_j} w_j(q) \ \ q\in D, r>0,
 \eeq
 where 
\begin{equation}\label{eq:defalphaj} 
\alpha_j:=-\frac{N-2}{2}+\sqrt{\left(\frac{N-2}{2}\right)^2 +\lambda_j}.
\end{equation}
 We claim that $w=\tilde w_j\big|_{\Omega_0}$ is the unique solution of
 \begin{equation}\label{eq:mixboundvalue3}
\begin{cases}
-\Delta w = 0 & \text{in $\Omega_0$},\\
 w = \frac{1}{N} w_j & \text{on $\Gamma_0$},\\
  \frac{\partial w}{\partial \nu} = 0 & \text{on $\Gamma_{1,0}$}.
\end{cases}
\end{equation}
Indeed, writing the Laplace operator in polar coordinates and exploiting \eqref{eq:eigenneumannsphere} we easily check that 
$$\Delta \tilde w_j=\frac{\partial^2 \tilde w_j}{\partial^2 r} + \frac{N-1}{r} \frac{\partial \tilde w_j}{\partial r} + \frac{1}{r^2} \Delta_{\S^{N-1}} \tilde w_j = \left( \alpha_j(\alpha_j-1)+ \alpha_j(N-1)-\lambda_j\right) \frac{r^{\alpha_j-2} w_j(q)}{N}=0,$$
because $\alpha_j$ satisfies $\alpha_j^2 + (N-2)\alpha_j -\lambda_j=0$. Moreover, by definition, we have $\tilde w_j\big|_D=\frac{1}{N}w_j$ and $\frac{\partial \tilde w_j}{\partial \nu} = 0$ on $\Gamma_{1,0}$.

Now, let us take $j=1$. It is well known that first eigenfunction $w_1$ is smooth and satisfies $\int_D w_1 \ d\sigma=0$, i.e. $w_1\in T_0M$. Computing $I^{\prime\prime}(0)[w_1,w_1]$, thanks to Proposition \ref{prop:secvartorzero} and taking into account  that $\tilde w_1\big|_{\Omega_0}$ is the solution of \eqref{eq:mixboundvalue3}, with $j=1$, we get that
\begin{equation}\label{eq:secvarthesis2}
\displaystyle I^{\prime\prime}(0)[w_1,w_1]=\displaystyle - \frac{1}{N^2}\int_{D}   w_1^2 \ d\sigma \displaystyle +\frac{1}{N} \int_{D}  w_1 \left(\frac{\partial \tilde w_1}{\partial\nu}\right) \ d\sigma.
\end{equation}
Then, since the $L^2(D)$-norm of $w_1$ is equal to 1, the exterior unit normal $\nu$ to $\Gamma_0$ is the radial direction, and
\beq\label{eq:dernormeigenw}
 \frac{\partial \tilde w_1}{\partial\nu}=\frac{1}{N}\alpha_1 r^{\alpha_1-1} w_1=\frac{\alpha_1}{N} w_1\ \ \hbox{on $D$},
\eeq
 from \eqref{eq:secvarthesis2} we obtain
$$
\displaystyle I^{\prime\prime}(0)[w_1,w_1]=\displaystyle - \frac{1}{N^2}  +\frac{\alpha_1}{N^2}.
$$
Thus we deduce
$$ I^{\prime\prime}(0)[w_1,w_1]<0 \ \ \hbox{if and only if} \ \ -1+\alpha_1<0. $$
Finally, from \eqref{eq:defalphaj} it is immediate to check that $\alpha_1<1$ is equivalent to $\lambda_1(D)<N-1$ and the proof of (i) is complete.\\

To prove (ii), let $v\in T_0M$ and assume, without loss of generality, that $\int_D v^2 \ d\sigma=1$. Taking $(w_j)_{j\in\N}$ as in the proof of (i), since $v \in T_0 M$ we can write $$v=\sum_{j=1}^\infty (v,w_j)_{L^2(D)} w_j.$$
Let $\tilde w_j$ be the harmonic extension of $w_j$ defined in \eqref{eq:harmonicextension}. Then, as $\tilde w_j\big|_{\Omega_0}$ is a solution to \eqref{eq:mixboundvalue3} for any $j\in\N$, we infer that
$\tilde v:=\sum_{j=1}^\infty (v,w_j)_{L^2(D)} \tilde w_j$ is a solution to
\begin{equation*}
\begin{cases}
-\Delta u = 0 & \text{in $\Omega_0$},\\
 u = \frac{1}{N} v & \text{on $\Gamma_0$},\\
  \frac{\partial u}{\partial \nu} = 0 & \text{on $\Gamma_{1,0}$}.
\end{cases}
\end{equation*}
Thus, by Proposition \ref{prop:secvartorzero}, we get
$$\displaystyle I^{\prime\prime}(0)[v,v]=\displaystyle - \frac{1}{N^2}\int_{D}   v^2\ d\sigma \displaystyle +\frac{1}{N} \int_{D}  v \ \frac{\partial \tilde v}{\partial\nu} \ d\sigma.$$
As in \eqref{eq:dernormeigenw} we have that $\frac{\partial \tilde w_j}{\partial\nu}=\frac{\alpha_j}{N} w_j$ on $D$, for any $j\in\N$. Hence, since $\int_D v^2 \ d\sigma=1$, we deduce
$$\displaystyle I^{\prime\prime}(0)[v,v]=\displaystyle - \frac{1}{N^2} +\frac{1}{N^2} \int_{D}  v \sum_{j=1}^\infty \alpha_j (v,w_j)_{L^2(D)} w_j \ d\sigma=- \frac{1}{N^2} +\frac{1}{N^2} \sum_{j=1}^\infty \alpha_j   (v,w_j)_{L^2(D)}^2.$$
Now, if $\lambda_1(D)>N-1$ it follows that $\alpha_1>1$, and, as $(\alpha_j)_{j\in\N}$ is a nondecreasing sequence, we obtain
\beq\label{eq:secformvv}
\displaystyle I^{\prime\prime}(0)[v,v]>- \frac{1}{N^2} +\frac{1}{N^2} \sum_{j=1}^\infty   (v,w_j)_{L^2(D)}^2=0,
\eeq
having used that $\sum_{j=1}^\infty   (v,w_j)_{L^2(D)}^2=1$, as $\int_D v^2(q) \ d\sigma=1$. Hence (ii) holds.

\end{proof}

We conclude this section with a useful criterion for checking the property $\lambda_1(D)<N-1$. To this end
let ${e} \in \S^{N-1}$ and let $u_{{e}}\in C^\infty(\R^N)$ be the function  defined by 
\beq\label{eq:defue}
u_{{e}}(x)=x \boldsymbol{\cdot} {e},
\eeq
which satisfies
\beq\label{eq:LaplaceBeltue}
-\Delta_{\S^{N-1}} u_{{e}}=(N-1) u_{{e}} \ \ \hbox{on $\S^{N-1}$.}
\eeq
We have
\begin{proposition}\label{lemma:unstablcrit}
Let $D$ be a smooth proper domain of $\S^{N-1}$ and let ${e}\in\S^{N-1}$ satisfy
$$ \int_{D} u_{{e}}\ d\sigma=0.$$
Assume that either one of the following holds:
\begin{itemize}
\item[(i)] $\displaystyle\int_{\partial D} u_{{e}} \frac{\partial u_{{e}}}{\partial \nu} \ d\hat\sigma<0$.\\[2pt]
\item[(ii)] $\displaystyle \int_{\partial D} u_{{e}} \frac{\partial u_{{e}}}{\partial \nu} \ d\hat\sigma=0$,
and $u_{{e}}$ is not an eigenfunction of
$
 \begin{cases} -\Delta_{\S^{N-1}} w=\lambda w& \hbox{in $D$},\\[8pt]\quad \quad \quad \frac{\partial w}{\partial \nu} =0 & \hbox{on $\partial D$},\end{cases}
$
\end{itemize}
where $d\hat\sigma$ is the $(N-2)$-dimensional area element of $\partial D$ and $\nu=\nu_{_{\partial D}}$ is the exterior unit co-normal to $\partial D$.
Then $\lambda_1(D)<N-1$.
\end{proposition}
\begin{proof}
Taking $u_{{e}}$ as test function in the variational characterization of the first non-trivial eigenvalue of $-\Delta_{\S^{N-1}}$ with zero Neumann condition on $\partial D$, applying Green's identity and exploiting \eqref{eq:LaplaceBeltue}, we have
$$\int_D |\nabla u_{{e}}|^2 \ d\sigma =\int_{\partial D} u_{{e}} \frac{\partial u_{{e}}}{\partial \nu} \ d\hat\sigma- \int_{D} u_{{e}}  \Delta u_{{e}} \ d\sigma=\int_{\partial D} u_{{e}} \frac{\partial u_{{e}}}{\partial \nu} \ d\hat\sigma+  (N-1)\int_{D} u^2_{{e}}  \ d\sigma.$$
Therefore, if (i) holds it follows that
\beq\label{eq:varchareigen}
\displaystyle\frac{\int_D |\nabla u_{{e}}|^2 \ d\sigma}{\int_{D} u_{{e}}^2  \ d\sigma}<N-1
\eeq
which implies that $\lambda_1(D)<N-1$. This completes the proof for the case (i). On the other hand, under the assumption (ii), the equality sign in \eqref{eq:varchareigen} holds, but as $u_{{e}}$ is not an eigenfunction it follows that $N-1$ cannot be the smallest non-trivial eigenvalue.
\end{proof}

\section{Existence of volume-constrained minimizers for the torsional energy}

Let $D\subset \S^{N-1}$ be a domain of the unit sphere and let $\Sigma_D$ be the cone generated by $D$. We will always assume that $D$ is smooth so that  $\Sigma_D$ is smooth exept at the vertex. In Sect. 3 we defined the torsional energy $\mathcal{E}(\Omega;\Sigma_D)$ for smooth domains  $\Omega\subset \Sigma_D$ strictly star-shaped with respect to the vertex of the cone. In this section we study the minimization problem for the torsional energy under a volume constraint in a larger class of sets. Thus we recall some definitions. 
\begin{definition}
We say that $\Omega\subset \RN$ is quasi-open, if for any $\varepsilon>0$, there exists an open set $\Lambda_\e$ such that $\mathrm{cap}(\Lambda_\e)\leq \e$ and $\Omega\cup\Lambda_\e$ is open, where $\mathrm{cap}(\Lambda_\e)$ denotes the capacity of $\Lambda_\e$.
\end{definition}
For any quasi-open set $\Omega\subset \Sigma_D$ we consider the Sobolev space:
$$H^1_0(\Omega; \Sigma_D):=\left\{u \in H^1(\Sigma_D); \ \ u=0 \ \ \hbox{q.e. on} \ \Sigma_D\setminus\Omega  \right\},$$
where q.e. means quasi-everywhere, i.e. up to sets of zero capacity.

\begin{definition}
We say that $u$ is a (weak) solution of the mixed boundary value problem
\begin{equation}\label{eq:mixbvprobquasiopen}
\begin{cases}
-\Delta u = 1 & \text{in $\Omega$,}\\
 u = 0 & \text{on $\partial\Omega \cap \Sigma_D$,}\\
 \frac{\partial u}{\partial \nu} = 0 & \text{on $\partial \Sigma_D$,}
\end{cases}
\end{equation}
if $u\in H^1_0(\Omega; \Sigma_D)$ and
$$\int_{\Sigma_D} \nabla u \boldsymbol{\cdot} \nabla v \ dx = \int_{\Sigma_D} v \ dx \ \ \ \ \forall v \in H^1_0(\Omega; \Sigma_D).$$
\end{definition}

\begin{remark}\label{rem:compacteigen}
As $\Sigma_D$ is connected and smooth (execpt at the vertex) then $\Sigma_D$ is uniformly Lipschitz. Thus if $|\Omega|<+\infty$ the inclusion $H^1_0(\Omega; \Sigma_D)\hookrightarrow L^2(\Sigma_D)$ is compact (see \cite[Proposition 2.3-(i)]{BV}). This implies that the functional
\beq\label{eq:deftorsion}
J(v)=\frac{1}{2}\int_{\Sigma_D} |\nabla v|^2 \ dx - \int_{\SD} v \ dx 
\eeq
has a unique minimizer $u_\Omega \in H^1_0(\Omega; \Sigma_D)$ which is the unique (weak) solution to \eqref{eq:mixbvprobquasiopen}, which is called \textit{energy function} or \textit{torsion function} of $\Omega$.  We also recall that $\Omega=\{u_\Omega>0\}$ up to a set of zero capacity (see \cite[Proposition 2.8-(e)]{BV}). Moreover we denote by $\lambda_1(\Omega; \Sigma_D)$ the first eigenvalue of the Laplacian in $H_0^1(\Omega; \Sigma_D)$, i.e. 
\beq\label{eq:deffirsteigen}
\lambda_1(\Omega; \Sigma_D)=\min_{v \in H_0^1(\Omega; \Sigma_D)\setminus\{0\}} \frac{\int_{\Sigma_D} |\nabla v|^2 \ dx}{\int_{\Sigma_D} v^2 \ dx}.
\eeq
Then, as before, we define the torsional energy of $\Omega$ (relative to $\Sigma_D$) as:
\beq\label{def:torsionquasiopen}
\mathcal{E}(\Omega; \Sigma_D)=J(u_\Omega)=-\frac{1}{2} \int_\Omega |\nabla u_\Omega|^2 \ dx=-\frac{1}{2} \int_\Omega u_\Omega \ dx. 
\eeq

\end{remark}
We want to study the problem of minimizing the functional $\E(\Omega;\Sigma_D)$ among quasi-open sets of uniformly bounded measure. Therefore, fixing $c>0$ we define
\beq\label{def:inftorsionvolumeconstraint}
\mathcal{O}_c(\Sigma_D):=\inf\{\E(\Omega;\Sigma_D); \ \ \Omega \ \hbox{quasi-open}, \ \Omega\subset \Sigma_D, \ |\Omega|\leq c\}.
\eeq
Our aim is to give a sufficient condition on the cone $\Sigma_D$ (hence on $D$) for the infimum in \eqref{def:inftorsionvolumeconstraint} to be achieved. We begin by recalling some known properties of the function $u_\Omega$ that will be used in this section.
 
 \begin{proposition}\label{prop1:sect8}
 Let $c>0$. There exists a positive constant $C$ depending only on $N$, $\Sigma_D$ and $c$ such that for any quasi-open subset $\Omega$ of $\Sigma_D$ with $|\Omega| \leq c$, it holds:
\begin{itemize}
\item[(i)] $u_\Omega$ is bounded and $\|u_\Omega\|_{L^\infty(\Sigma_D)}\leq C |\Omega|^{2/N}$;
\item[(ii)] $\int_{\Sigma_D} |\nabla u_\Omega|^2\  dx \leq C |\Omega|^{\frac{N+2}{N}}$;
\item[(iii)] $\int_{\Sigma_D} u_\Omega^2\  dx \leq C |\Omega|^{\frac{N+4}{N}}$.
\end{itemize}
 \end{proposition}  
\begin{proof}
Let us fix $c>0$. Since the cone $\Sigma_D$ is a uniformly Lipschitz connected open set of $\RN$ then we can apply \cite[Lemma 2.5]{BV}. Hence, for any quasi-open subset $\Omega\subset\Sigma_D$, with $|\Omega|\leq c$, fixing $p \in ]N/2, +\infty[$ and taking $f=\chi_\Omega$, where $\chi_\Omega$ denotes the characteristic function of $\Omega$,   we obtain from \cite[Lemma 2.5]{BV} that there exists a positive constant $\tilde C$ depending on $N$, $p$, $\Sigma_D$ and $c$ only such that
$$
\|u_\Omega\|_{L^\infty(\Sigma_D)} \leq \tilde C \|f\|_{L^p(\Sigma_D)} |\Omega|^{2/N - 1/p}=  \tilde C  |\Omega|^{1/p}  |\Omega|^{2/N - 1/p} = \tilde C|\Omega|^{2/N},
$$
which gives (i).

Next, taking $u_\Omega$ as test function in the weak formulation of \eqref{eq:mixbvprobquasiopen} we get $$\int_\Omega |\nabla u_\Omega|^2 \ dx = \int_\Omega  u_\Omega \ dx, $$
and, by (i), we obtain
$$\int_\Omega |\nabla u_\Omega|^2 \ dx \leq C |\Omega|^{2/N} |\Omega|= C |\Omega|^{\frac{N+2}{N}},$$
i.e. (ii). Finally (iii) is a trivial consequence of (i) since
$$ \int_{\Omega} u_\Omega^2 \ dx\leq \|u_\Omega\|_{L^\infty(\Sigma_D)}^2 |\Omega|\leq C^2 |\Omega|^{4/N}|\Omega|=C^2  |\Omega|^{\frac{N+4}{N}}.$$
\end{proof}

Notice that as a straightforward consequence of the previous result it holds that $\O_c(\Sigma_D)>-\infty$.

\begin{remark}\label{remark1:sect8}
As remarked in \cite[Remark 4.2]{PT2} there is a natural invariance by scaling in our problem, which, in particular, allows to claim that the infimum as in \eqref{def:inftorsionvolumeconstraint}, but with volume bounded by another constant $\lambda>0$, can be easily computed from $\O_c(\Sigma_D)$. Namely we have
\beq\label{eq:relationscalinginfenergy}
\lambda^{-\frac{N+2}{N}}\O_\lambda(\Sigma_D)=c^{-\frac{N+2}{N}}\O_c(\Sigma_D)=\O_1(\Sigma_D).
\eeq
 Indeed, for any quasi-open $\Omega \subset \Sigma_D$, for any $t>0$ it holds that $t\Omega\subset \Sigma_D$, $|t\Omega|=t^N|\Omega|$, and it is easy to check that $u_{t\Omega}(x)=t^2 u_\Omega\left(\frac{x}{t}\right)$ and
\beq\label{eq:scalinvariance}
\E(t\Omega; \Sigma_D)=t^{N+2} \E(\Omega; \Sigma_D).
\eeq
In particular $\O_c(\Sigma_D)$ can be defined by taking $|\Omega|=c$ in \eqref{def:inftorsionvolumeconstraint} and either a minimizer exists for any fixed volume or there are no minimizers whatever bound for the volume is chosen.
\end{remark}

Among the quasi-open sets in $\Sigma_D$ we can consider the spherical sectors 
\beq\label{def:sphericalsector}
\Omega_{D,R}:=\Sigma_D\cap B_R(0). 
\eeq
In this case the solution of \eqref{eq:mixbvprobquasiopen} is radial and explicitly given by
\beq\label{def:radialfundamentalsoltorsion}
u_{\Omega_{D,R}}(x)=\begin{cases} \frac{R^2-|x|^2}{2N} & \hbox{if $x\in \Omega_{D,R}$},\\
0 & \hbox{if $x\in \Sigma_D\setminus{\Omega_{D,R}}$},
\end{cases}
\eeq
and its energy is
\beq\label{eq:torsionfundamental}
 \E(\Omega_{D,R}; \Sigma_D)=-\frac{1}{2N^2(N+2)} R^{N+2} \mathcal{H}_{N-1}(D).
 \eeq
 Therefore, by \eqref{def:inftorsionvolumeconstraint}, we have
 \beq\label{eq:negativitytenergy}
\O_c(\Sigma_D) \leq \E(\Omega_{D,R_c}; \Sigma_D)<0,
 \eeq
 where $R_c=R_c(D)>0$ is such that $|\Omega_{D,R_c}|=R_c^{N} \mathcal{H}_{N-1}(D)=c$, namely $R_c(D)=\left(\frac{c}{ \mathcal{H}_{N-1}(D)}\right)^{\frac{1}{N}}$. 
 \begin{remark}\label{rem:monotonenergyODRc}
 Notice that for any $c>0$ it holds
 \beq\label{eq:energysphericalsector}
 \E(\Omega_{D,R_c(D)}; \Sigma_D)=-\frac{1}{2N^2(N+2)} c^{\frac{N+2}{N}} \left[\mathcal{H}_{N-1}(D)\right]^{-\frac{2}{N}}
 \eeq
which means that $\E(\Omega_{D,R_c}; \Sigma_D)$ is monotone increasing with respect to $\mathcal{H}_{N-1}(D)$.
\end{remark}
 \begin{remark}\label{rem:energyhalfball}
When $D$ is a hemisphere, let us say for convenience the upper hemisphere, denoted by $\S^{N-1}_+=\S^{N-1}\cap\{(x_1,\ldots,x_N)\in \RN; \ x_N>0\}$, then the cone $\Sigma_{\S^{N-1}_+}$ coincides with the half-space $\R^N_+$. In this case it is well known, for example by symmetrization, that $\O_c\left(\Sigma_{\S^{N-1}_+}\right)$ is achieved by any half-ball of measure $c$ and 
\beq\label{eq:exprtorsionhemisphereR}   
\O_c\left(\Sigma_{\S^{N-1}_+}\right)= \E\left(\Omega_{\S^{N-1}_+,R_c(\S^{N-1}_+)}; \Sigma_{\S^{N-1}_+}\right) =-\frac{{\omega}_N}{4N(N+2)} \left(\frac{2c}{ N\omega_N}\right)^{\frac{N+2}{N}}.
\eeq
\end{remark}
In the general case, using the smoothness of the cone, we prove in Proposition \ref{prop:ineqtesi} that it always holds
\beq\label{eq:ineqtesi}
\O_c\left(\Sigma_D\right) \leq \O_c\left(\Sigma_{\S^{N-1}_+}\right).
\eeq
The main result of this section is to show that if the strict inequality holds in \eqref{eq:ineqtesi} then the infimum is achieved. Indeed we have:
\begin{theorem}\label{teo:mainteoexistmin}
Let $c>0$ and assume that
\beq\label{eq:conditionD}
 \O_c\left(\Sigma_D\right)< \O_c\left(\Sigma_{\S^{N-1}_+}\right),
 \eeq
 then $\O_c(\Sigma_D)$ is achieved.
\end{theorem}

\begin{proof}
Let $(\Omega_n)_n \subset \Sigma_D$ be a minimizing sequence for $\O_c(\Sigma_D)$ and consider the  corresponding energy functions $u_{\Omega_n}\in H^1_0(\Omega_n;\Sigma_D)$ for any $n\in \N$. By definition we have $$\E(\Omega_n; \Sigma_D)=-\frac{1}{2}\int_{\Omega_n} u_{\Omega_n} \ dx \to \mathcal{O}_c(\Sigma_D), \ \hbox{ as $n\to +\infty$}.$$
Setting $u_n:=u_{\Omega_n}$, since $|\Omega_n|\leq c$, for any $n\in \N$, by Proposition \ref{prop1:sect8} we find a positive constant $C_1$ independent of $n$ such that 
\beq\label{eq:boundH1teosect8}
\|u_n\|_{H^1(\Sigma_D)} \leq C_1 \ \ \forall n \in \N.
\eeq
In particular, up to a subsequence (still denoted by $(u_n)_n$), we have $\|u_n\|_{L^2(\Sigma_D)}^2\to \lambda$, for some $\lambda\geq 0$. We first observe that $\lambda>0$. Otherwise, if $\lambda=0$, by H\"older's inequality and exploiting the uniform bound $|\Omega_n|\leq c$, we would have
\beq\label{eq2:teo1sect8}
\|u_n\|_{L^1(\Sigma_D)}\to 0, \ \hbox{as $n\to +\infty$},
\eeq
which implies that $\E(\Omega_n; \Sigma_D)\to 0$, as $n\to +\infty$, contradicting  $\O_c(\Sigma_D)<0$ (see \eqref{eq:negativitytenergy}). Therefore, as $(u_n)_n$ is bounded in $H^1(\Sigma_D)$ and
\beq\label{eq0:teo1sect8}
\|u_n\|^2_{L^2(\Sigma_D)} \to \lambda,\ \ \hbox{for $n\to +\infty$},
\eeq
for some $\lambda>0$, we can apply, with small modifications in the proof, the concentration-compactness principle of P. L. Lions (see \cite[Lemma III.1]{LIO}). Hence, there exists a subsequence $(u_{n_k})_k$ satisfying one of three following possibilities:
\begin{itemize}
\item[(i)] there exists $(y_{n_k})_k \subset \overline{\Sigma_D}$ satisfying $$ \forall \e>0 \ \exists R>0 \ \ \hbox{such that} \  \ \int_{B_R(y_{n_k})\cap \Sigma_D} u_{n_k}^2 \ dx \geq \lambda - \eps \ \ \forall k\in \N;\vspace{-4pt}$$
\item[(ii)] $\displaystyle \lim_{k\to +\infty} \sup_{y \in \Sigma_D} \int_{B_R(y)\cap \Sigma_D} u_{n_k}^2 \ dx=0$, for all $R>0$;\\[2pt]
\item[(iii)] there exists $\alpha\in]0,\lambda[$ such that for all $\e>0$, there exist $k_0\geq 1$ and two sequences $(u_{1,k})_k$, $(u_{2,k})_k$ bounded in $H^1(\Sigma_D)$ satisfying, for $k\geq k_0$
$$\|u_{n_k} - u_{1,k}-u_{2,k}\|_{L^2(\Sigma_D)} \leq 4\e,  \ \left| \int_{\Sigma_D} u_{1,k}^2 \ dx -\alpha \right| \leq \e, \  \ \left| \int_{\Sigma_D} u_{2,k}^2 \ dx -(\lambda-\alpha) \right| \leq \e,$$
$$\textrm{dist}(\textrm{supp}(u_{1,k}), \textrm{supp}(u_{2,k})) \to +\infty, \ \hbox{as $k\to+\infty$},$$
$$\liminf_{k\to +\infty} \int_{\Sigma_D} |\nabla u_{n_k}|^2 - |\nabla u_{1,k}|^2 -  |\nabla u_{2,k}|^2 \ dx \geq 0.$$
\end{itemize}
  We now divide the proof in some steps. We begin by showing that the ``vanishing'' case (ii) cannot occur.\\
 
 \noindent\textbf{Step 1:} (ii) cannot happen.\\
 
 Assume by contradiction that (ii) holds. The idea is to show that $u_{n_k}\to 0$ strongly in  $L^2(\Sigma_D)$, as $k\to +\infty$, contradicting \eqref{eq0:teo1sect8}. To prove this we invoke \cite[Lemma 1.21]{WI} (whose proof can be easily adapted for functions in $H^1(\Sigma_D)$), which claims that (ii) and \eqref{eq:boundH1teosect8} imply $u_{n_k}\to 0$ in  $L^p(\Sigma_D)$, for any $2<p<2^*$, as $k\to +\infty$, where $2^*=\frac{2N}{N-2}$ is the critical Sobolev exponent. Then, exploiting that $u_{n_k}\in H_0^1(\Omega_{n_k}; \Sigma_D)$ and $|\Omega_{n_k}|\leq c$, by H\"older's inequality we readily conclude that $u_{n_k}\to 0$ in $L^2(\Sigma_D)$, as $k\to +\infty$.\\

In the next step we prove that the ``dichotomy'' case (iii) cannot occur.\\

\noindent\textbf{Step 2:} (iii) cannot happen.\\

Assume by contradiction that (iii) holds.  
We claim that, up to a further subsequence, there exists  another minimizing sequence $(\widetilde \Omega_{n_k})_k\subset \Sigma_D$, with $\widetilde \Omega_{n_k} \subset \Omega_{n_k}$, for any $k$, satisfying:
\begin{itemize}
\item $\widetilde\Omega_{n_k} = \Omega_{1,k} \cup \Omega_{2,k}$, for some quasi-open subsets $\Omega_{1,k}$, $\Omega_{2,k}$ of $\Omega_{n_k}$;
\item  $\mathrm{dist}( \Omega_{1,k},  \Omega_{2,k})\to +\infty$, as $k\to +\infty$;
\item $c_i:=\liminf_{k\to +\infty} |\Omega_{i,k}|>0$, for $i=1,2$.
\end{itemize}
Indeed, by (iii) and a diagonal argument, we find bounded subsequences $(u_{1,k})_k$, $(u_{2,k})_k$ in $H^1(\Sigma_D)$ (still indexed by $k$) satisfying
\beq\label{eq:dichotomyLions}
\begin{array}{lll}
&&\displaystyle \|u_{n_k} - u_{1,k}-u_{2,k}\|_{L^2(\Sigma_D)} \to 0, \ \ \int_{\Sigma_D} u_{1,k}^2 \ dx \to \alpha, \ \ \int_{\Sigma_D} u_{2,k}^2 \ dx  \to (\lambda-\alpha), \ \hbox{as $k\to +\infty$,}\\[10pt]
&&\displaystyle\textrm{dist}(\textrm{supp}(u_{1,k}), \textrm{supp}(u_{2,k})) \to +\infty, \ \hbox{as $k\to+\infty$},\\[6pt]
&&\displaystyle\liminf_{k\to +\infty} \int_{\Sigma_D} |\nabla u_{n_k}|^2 - |\nabla u_{1,k}|^2 -  |\nabla u_{2,k}|^2 \ dx \geq 0.
\end{array}
\eeq
By the proof of \cite[Lemma III.1]{LIO} we see that $u_{1,k}$, $u_{2,k}$ can be chosen to be non-negative and in addition, since $u_{n_k} \in H^1_0(\Omega_{n_k};\Sigma_D)$, we also have that $u_{1,k}, u_{2,k} \in H^1_0(\Omega_{n_k};\Sigma_D)$ for any $k$. In particular, as $u_{1,k}, u_{2,k} \in H^1(\Sigma_D)$, setting $\Omega_{1,k}:=\{u_{1,k}>0\}$, $\Omega_{2,k}:=\{u_{2,k}>0\}$ it follows that $\Omega_{1,k}$, $\Omega_{2,k}$ are quasi-open subsets of $\Sigma_D$. Therefore, $\widetilde \Omega_{k}:= \Omega_{1,k} \cup \Omega_{2,k}$ is a quasi-open set contained in $\Omega_{n_k}$ and denoting by $\tilde u_{n_k}:=u_{\widetilde \Omega_{n_k}}$ the torsion function of $\widetilde \Omega_{n_k}$ and arguing as in \cite[Sect. 3.3]{BU} (with obvious small modifications), we infer that
\beq\label{eq3:teo1sect8}
 \|u_{n_k} - \tilde u_{n_k}\|_{H^1(\Sigma_D)}\to 0, \ \ \hbox{as $k\to+\infty$}.
\eeq
From \eqref{eq3:teo1sect8} it follows that $(\widetilde \Omega_{n_k})_k$ is a minimizing sequence for $\mathcal{O}_c(D)$. 
Moreover, by construction and \eqref{eq:dichotomyLions} we readily deduce that $\mathrm{dist}( \Omega_{1,k},  \Omega_{2,k})\to +\infty$, as $k\to +\infty$. Finally, setting
\beq\label{eq:decomeganc1c2}
c_i:=\liminf_{n\to+\infty} |\Omega_{i,k}|, \ \ i=1,2,
\eeq
it holds that
$c_i>0$ for $i=1,2$.  Indeed, assuming by contradiction, for instance, that $c_1=0$, by H\"older's inequality and Sobolev's inequality (note that $\Sigma_D$ satisfies the cone condition) we get
$$\int_{\Sigma_D} u_{1,k}^2 \ dx \leq \left(\int_{\Sigma_D} |u_{1,k}|^{2^*} \ dx \right)^{\frac{2}{2^*}} |\Omega_{1,k}|^{\frac{2}{N}} \leq C(N,\Sigma_D) \int_{\Sigma_D} |\nabla u_{1,k}|^2 \ dx \ |\Omega_{1,k}|^{\frac{2}{N}}.$$
Now, recalling that $(u_{1,k})_k$ is a bounded sequence in $H^1(\Sigma_D)$, from the previous inequality and since we are assuming $c_1=0$  we deduce that
$$\liminf_{k\to+\infty} \int_{\Sigma_D} u_{1,k}^2 \ dx=0,$$
which contradicts \eqref{eq:dichotomyLions}. Hence $c_1>0$, and by the same argument we infer that $c_2>0$. The proof of the claim is complete.\\

In order to conclude the proof of \textbf{Step 2} we show that the previous claim leads to a contradiction. To this end we begin observing that by invariance by dilatation (see Remark \ref{remark1:sect8}) it follows that our minimization problem is equivalent to $$\mathcal{M}(\Sigma_D):=\inf\left\{\frac{\E(\Omega; \Sigma_D)}{|\Omega|^{\frac{N+2}{N}}}; \ \  \Omega \ \hbox{quasi-open}, \ \Omega\subset \Sigma_D,\ |\Omega|>0\right\}.$$
In particular by \eqref{eq:relationscalinginfenergy} and a straightforward computation we check that 
$$\O_c(\Sigma_D)=c^{\frac{N+2}{N}}\O_1(\Sigma_D)=c^{\frac{N+2}{N}} \mathcal{M}(\Sigma_D),$$
and a minimizing sequence for $\mathcal{M}(\Sigma_D)$ is given by $\Lambda_k:=|\widetilde \Omega_{n_k} |^{-\frac{1}{N}}\widetilde \Omega_{n_k}$. Then, setting $\Lambda_{i,k}:=|\widetilde \Omega_{n_k} |^{-\frac{1}{N}}\Omega_{i,k}$, $i=1,2$, and in view of the previous claim, up to a subsequence, we have $\Lambda_k= \Lambda_{1,k} \cup  \Lambda_{2,k}$, where $|\Lambda_{i,k}|\to \frac{c_i}{c_1+c_2}>0$, as $k\to +\infty$, $i=1,2$, and $\Lambda_{1,k} \cap \Lambda_{2,k}=\emptyset$ for all sufficiently large $k$. 
Now, as $\Lambda_{i,k} \subset \Sigma_D$ are quasi-open subsets with positive measure, then by definition of $\mathcal{M}(\Sigma_D)$, for any $i=1,2$, we have 
\beq\label{eq23:teo1sect8}
 \frac{\E(\Lambda_{i,k}; \Sigma_D)}{|\Lambda_{i,k}|^{\frac{N+2}{N}}}\geq \mathcal{M}(\Sigma_D).
\eeq
In addition, assuming without loss of generality that  $|\Lambda_{1,k}|\geq |\Lambda_{2,k}|$, by an elementary computation, we deduce that
\beq\label{eq24:teo1sect8}
\begin{array}{lll}
\displaystyle \left(|\Lambda_{1,k}|+|\Lambda_{2,k}|\right)^{\frac{N+2}{N}}&=& \displaystyle|\Lambda_{1,k}|^{\frac{N+2}{N}}+\frac{N+2}{N} |\Lambda_{1,k}|^{\frac{2}{N}}|\Lambda_{2,k}|+\frac{N+2}{N^2} (|\Lambda_{1,k}|+\xi_k)^{\frac{2-N}{N}} |\Lambda_{2,k}|^2\\[6pt]
&\geq&\displaystyle|\Lambda_{1,k}|^{\frac{N+2}{N}}+|\Lambda_{2,k}|^{\frac{N+2}{N}}+\frac{N+2}{N^2} (|\Lambda_{1,k}|+\xi_k)^{\frac{2-N}{N}} |\Lambda_{2,k}|^2,
\end{array}
\eeq
where $\xi_k \in [0,|\Lambda_{2,k}|]$.
Then, setting $\mathcal{M}_k(\Sigma_D):=\frac{\E(\Lambda_k; \Sigma_D)}{|\Lambda_k|^{\frac{N+2}{N}}}<0$, recalling that $(\Lambda_k)_k$ is minimizing for $\mathcal{M}(\Sigma_D)$, exploiting the properties of $\Lambda_k$ and taking into account \eqref{eq24:teo1sect8}, \eqref{eq23:teo1sect8}, we infer that for all sufficiently large $k$ it holds
\beq\label{eq25:teo1sect8}
\begin{array}{lll}
\displaystyle \E(\Lambda_{k}; \Sigma_D)&=& \displaystyle \mathcal{M}_k(\Sigma_D)|\Lambda_k|^{\frac{N+2}{N}}\\[6pt]
&=& \displaystyle \left(\mathcal{M}(\Sigma_D)+o(1)\right) \left(|\Lambda_{1,k}|+|\Lambda_{2,k}|\right)^{\frac{N+2}{N}}\\[6pt]
&\leq& \displaystyle  \left(\mathcal{M}(\Sigma_D)+o(1)\right) \left(|\Lambda_{1,k}|^{\frac{N+2}{N}}+|\Lambda_{2,k}|^{\frac{N+2}{N}}+\frac{N+2}{N^2} \xi_k^{\frac{2-N}{N}} |\Lambda_{2,k}|^2\right)\\[6pt]
&=& \displaystyle \mathcal{M}(\Sigma_D) |\Lambda_{1,k}|^{\frac{N+2}{N}}+ \mathcal{M}(D)|\Lambda_{2,k}|^{\frac{N+2}{N}}+\mathcal{M}(\Sigma_D)\frac{N+2}{N^2}(|\Lambda_{1,k}|+\xi_k)^{\frac{2-N}{N}}  |\Lambda_{2,k}|^2 + o(1)\\[6pt]
&\leq& \displaystyle \E(\Lambda_{1,k}; \Sigma_D)+ \E(\Lambda_{2,k}; \Sigma_D)+\mathcal{M}(\Sigma_D)\frac{N+2}{N^2}(|\Lambda_{1,k}|+\xi_k)^{\frac{2-N}{N}}  |\Lambda_{2,k}|^2 + o(1)\\[6pt]
&=& \displaystyle \E(\Lambda_{k}; \Sigma_D)+\mathcal{M}(\Sigma_D)\frac{N+2}{N^2}(|\Lambda_{1,k}|+\xi_k)^{\frac{2-N}{N}}  |\Lambda_{2,k}|^2 + o(1)\\[6pt]
&<&\displaystyle \E(\Lambda_{k}; \Sigma_D)
\end{array}
\eeq
where the last inequality is strict because $|\Lambda_{i,k}|\to \frac{c_i}{c_1+c_2}>0$, for $i=1,2$, $k\to +\infty$, and $\mathcal{M}(\Sigma_D)<0$. Clearly \eqref{eq25:teo1sect8} is contradictory and this concludes the proof of \textbf{Step 2}.\\

From \textbf{Step 1} and \textbf{Step 2} we know that the only admissible case is (i). 
%
Roughly speaking (i) states that, there exists a sequence $(y_{n_k})_k \subset \overline{\Sigma_D}$ such that for a sufficiently large ball $B_R(y_{n_k})$, the mass of $u_{n_k}$ is concentrated in  $B_R(y_{n_k})\cap \Omega_{n_k}$, while the part in the complement $B_R^\complement(y_k)\cap \Omega_{n_k}$ is negligible. With (i) at hand we can show that the same happens for the energy, in particular the possible tails of $\Omega_{n_k}$ do not play a role. This is the content of the next technical step:\\ 
 
\noindent\textbf{Step 3:} For any fixed $\eps>0$ there exist $\bar R>1$ and $\bar k \in \N$, both depending only on $\eps$, such that
\beq\label{eq:claim2Teosect8}
\E(\Omega_{n_k}; \Sigma_D)\geq \E(B_{2 R}(y_{n_k})\cap \Omega_{n_k}; \Sigma_D) - 2c\sqrt{2\eps}, \ \ \forall k\geq \bar k \ \forall R\geq \bar R.
\eeq

Let us fix $\eps>0$ and let $R>0$ be tha radius given by (i). Let $\varphi \in C^\infty_c(\RN)$ such that $0\leq \varphi \leq 1$, $\varphi\equiv 1$ in $B_{R}(0)$, $\varphi\equiv 0$ in $B_{2R}^\complement(0)$ and $|\nabla \varphi|\leq \frac{C_0}{R}$ in $\RN$, where $C_0>0$ is a constant independent of $R$. We set $\varphi_k(x):=\varphi(x-y_{n_k})$ and observe that 
\beq\label{eq16:teo1sect8}
\begin{array}{lll}
\displaystyle \int_{\Sigma_D} |\nabla u_{n_k}|^2 \ dx &\geq&\displaystyle  \int_{\Sigma_D} |\nabla u_{n_k}|^2\varphi_k^2 \ dx\\[12pt]
&=&\displaystyle \int_{\Sigma_D} |\nabla (u_{n_k}\varphi_k)|^2 \ dx- \underbrace{2 \int_{\Sigma_D}  u_{n_k} \varphi_k \nabla u_{n_k} \boldsymbol{\cdot} \nabla \varphi_k \ dx}_{\textbf{(I)}} -  \underbrace{\int_{\Sigma_D} u_{n_k}^2 |\nabla  \varphi_k |^2 \ dx}_{\textbf{(II)}}.
\end{array}
\eeq
Then, exploiting the properties of $\varphi_k$, H\"older's inequality and \eqref{eq:boundH1teosect8}  we infer that
$$|\textbf{(I)}| \leq  \frac{C_0}{R} \|\nabla u_{n_k}\|_{L^2(\Sigma_D)} \|u_{n_k}\|_{L^2(\Sigma_D)} \leq \frac{C_0C_1^2}{R},$$
where $C_0$, $C_1$ are both independent of $k$ and $R$. Similarly, for $|\textbf{(II)}|$ we have
$$|\textbf{(II)}| \leq  \frac{C_0^2C_1}{R^2},$$
and thus by \eqref{eq16:teo1sect8} and assuming without loss of generality that $R>1$ we obtain that
\beq\label{eq17:teo1sect8}
 \int_{\Sigma_D} |\nabla u_{n_k}|^2 \ dx \geq \int_{\Sigma_D} |\nabla (u_{n_k}\varphi_k)|^2 \ dx- \frac{C_2}{R},
\eeq
where $C_2>0$ is independent of $k$ and $R$. In addition, let us write
$$ -\int_{\Sigma_D}  u_{n_k} \ dx = -\int_{\Sigma_D}  u_{n_k}\varphi_k \ dx - \underbrace{\int_{\Sigma_D\cap B_R(y_k)}  u_{n_k}(1-\varphi_k) \ dx}_{\textbf{(III)}} - \underbrace{\int_{\Sigma_D\cap B_R^\complement(y_k)}  u_{n_k}(1-\varphi_k) \ dx}_{\textbf{(IV)}}.$$
We first observe that $\textbf{(III)}=0$, because $\varphi_k\equiv1$ in $B_R(y_k)$, while for \textbf{(IV)}, applying H\"older's inequality, taking into account that $u_{n_k}=0$ q.e. on $\Sigma\setminus\Omega_{n_k}$ and the properties of $\varphi_k$, we get that 
$$|\textbf{(IV)}| \leq \left( \int_{\Sigma_D\cap B_R^\complement(y_k)} u_{n_k}^2 \ dx\right)^{\frac{1}{2}} |B_{2R}(y_k)\cap\Omega_{n_k}|.$$
Now, thanks to (i) and \eqref{eq0:teo1sect8} it follows that $\|u_{n_k}\|_{L^2(\Sigma_D\cap B_R^\complement(y_k))} \leq \sqrt{2\eps}$ for all sufficiently large $k$, and thus, as $|\Omega_{n_k}|\leq c$, we deduce that 
$$|\textbf{(IV)}| \leq \sqrt{2\eps} c.$$
Summing up, we have proved that
\beq\label{eq18:teo1sect8}
-\int_{\Sigma_D}  u_{n_k} \ dx \geq -\int_{\Sigma_D}  u_{n_k}\varphi_k \ dx -  \sqrt{2\eps} c.
\eeq
Hence, combining \eqref{eq17:teo1sect8}, \eqref{eq17:teo1sect8} and recalling the definition of the functional $J$ (see \eqref{eq:deftorsion}) we obtain
$$J(u_{n_k})\geq J(u_{n_k}\varphi_k) - \frac{C_2}{2R}-\sqrt{2\eps} c.$$
Since $\eps$ is fixed and $C_2$ is independent of $R$ and $k$, up to taking a larger $R$, we can assume that $\frac{C_2}{2R}<\sqrt{2\eps} c$. Then, observing that $u_{n_k}\varphi_k \in H_0^1(B_{2R}(y_{n_k})\cap \Omega_{n_k}; \Sigma_D)$ we finally get
$$\E(\Omega_{n_k}; \Sigma_D)=J(u_{n_k})\geq J(u_{n_k}\varphi_k) - 2\sqrt{2\eps} c \geq \E(B_{2R}(y_{n_k})\cap \Omega_{n_k}; \Sigma_D) - 2\sqrt{2\eps} c,$$
so that \textbf{Step 3} is proved.\\
 
In the next step we prove that the sequence of points $(y_{n_k})_k\subset \overline\Sigma_D$ provided by (i) is bounded.\\

 \noindent\textbf{Step 4:} The sequence $(y_{n_k})_k\subset \overline\Sigma_D$ existing by (i) is bounded.\\

Assume by contradiction that there exists a subsequence (still indexed by $k$) such that $$\lim_{k\to+\infty}|y_{n_k}|= +\infty.$$ 

Thanks to the assumption \eqref{eq:conditionD} we can fix $\eps>0$ sufficiently small so that 
\beq\label{eq20:teo1sect8}
\mathcal{O}_c(\Sigma_D)+2c\sqrt{2\eps}<\mathcal{O}_c\left(\Sigma_{\S^{N-1}_+}\right),
\eeq
and by \textbf{Step 3} we find $R$ sufficiently large depending only on $\eps$, such that  for all sufficiently large $k$
\beq\label{eq:thesisStep3}
\E(\Omega_{n_k}; \Sigma_D)\geq \E(B_{2R}(y_{n_k})\cap \Omega_{n_k}; \Sigma_D) - 2c\sqrt{2\eps}.
\eeq
We observe that $B_{2R}(y_{n_k})\cap \Omega_{n_k}$ intersects the boundary of $\Sigma_D$. More precisely, for all sufficiently large $k$, it holds that
\beq\label{eq:Hintersect}
\mathcal{H}_{N-1}\left(\overline{(B_{2R}(y_{n_k})\cap \Omega_{n_k})}\cap\partial\Sigma_D\right)>0.
\eeq
Indeed, on the contrary, setting for convenience $\Theta_{R,k}:=B_{2R}(y_{n_k})\cap \Omega_{n_k}$ there exists a subsequence (still indexed by $k$) such that $\mathcal{H}_{N-1}(\overline{\Theta_{R,k}}\cap\partial\Sigma_D)=0$ for all $k\in \N$, and by the same argument of \cite[Remark 4.3]{BV} we conclude that $H_0^1(\Theta_{R,k};\Sigma_D)= H_0^1(\Theta_{R,k})$ and thus 
\beq\label{eq19:teo1sect8bis}
\E( \Theta_{R,k}; \Sigma_D)=\E(\Theta_{R,k}; \RN),
\eeq
where $\E(\Theta_{R,k}; \RN)$ denotes the ``free'' energy of $\Theta_{R,k}$, under  a homogeneous Dirichlet condition, namely, $\E(\Theta_{R,k}; \RN)$ is the minimizer in $H_0^1(\Theta_{R,k})$ of the functional $J(v)=\frac{1}{2}\int_{\RN} |\nabla v|^2 \ dx - \int_{\RN} v \ dx$.
 Then, by considering the Schwartz symmetrization  $u_{\Theta_{R,k}}^*$ of the energy function $u_{\Theta_{R,k}}$ associated to $\Theta_{R,k}$, and thanks to the classical P\'olya-Szeg\"o inequality, we infer that  
\beq\label{eq19:teo1sect8}
\begin{array}{lll}
\displaystyle\E( \Theta_{R,k}; \RN)&=&\displaystyle \frac{1}{2}\int_{\Theta_{R,k}} |\nabla u_{\Theta_{R,k}}|^2 \ dx - \int_{ \Theta_{R,k}} u_{\Theta_{R,k}} \ dx\\[12pt]
& \geq&\displaystyle  \frac{1}{2}\int_{\Theta_{R,k}^*} |\nabla u_{\Theta_{R,k}}^*|^2 \ dx - \int_{ \Theta_{R,k}^*} u_{\Theta_{R,k}}^* \ dx\\[12pt]
 & \geq&\displaystyle \E(\Theta_{R,k}^*; \RN),
\end{array}
\eeq
Hence, as $\Theta_{R,k}^*$ is a ball, with $c_k:=|\Theta_{R,k}^*|=|\Theta_{R,k}|\leq c$, then from \eqref{eq19:teo1sect8bis},\eqref{eq19:teo1sect8}, taking into account Remark \ref{rem:monotonenergyODRc} and \eqref{eq:exprtorsionhemisphereR} (noticing that $\E(\Theta_{R,k}^*; \RN)=\E(\Omega_{\S^{N-1},R_{c_k}(\S^{N-1})}; \Sigma_{\S^{N-1}})$, where $\Omega_{\S^{N-1},R_{c_k}(\S^{N-1})}$ is the ball centred at the origin of radius $R_{c_k}(\S^{N-1})$ with $|\Omega_{\S^{N-1},R_{c_k}(\S^{N-1})}|=c_k$, see \eqref{def:sphericalsector}, \eqref{eq:energysphericalsector}), we deduce that
\beq\label{eq21:teo1sect8}
\E(\Theta_{R,k}; \Sigma_D)  \geq \E(\Theta_{R,k}^*; \RN)> \O_{c_k}\left(\Sigma_{\S^{N-1}_+}\right) \geq \mathcal{O}_c\left(\Sigma_{\S^{N-1}_+}\right).
\eeq
Finally, recalling that $\Theta_{R,k}=B_{2R}(y_{n_k})\cap \Omega_{n_k}$, then from \eqref{eq:claim2Teosect8} and \eqref{eq21:teo1sect8} we have, for large $k$,
$$\E(\Omega_{n_k}; \Sigma_D)\geq\mathcal{O}_c\left(\Sigma_{\S^{N-1}_+}\right)- 2c\sqrt{2\eps}.$$
Hence passing to the limit as $k\to +\infty$ we conclude that
$$\O_c(\Sigma_D)\geq\mathcal{O}_c\left(\Sigma_{\S^{N-1}_+}\right)- 2c\sqrt{2\eps},$$
which contradicts \eqref{eq20:teo1sect8}.\\

Then, by \eqref{eq:Hintersect}, there exists $k_0\in \N$ such that $\mathrm{dist}(y_{n_k}, \partial\Sigma_D)\leq 2R$ for all $k\geq k_0$ and we can find a sequence of points $(z_k)_k \subset\partial\Sigma_D\setminus\{0\}$ such that $z_k \in (\overline{B_{2R}(y_{n_k})\cap \Omega_{n_k}})\cap\partial\Sigma_D$ and $B_{2R}(y_{n_k})\subset B_{4R}(z_k)$, for all $k\geq k_0$. Then, by monotonicity of the torsional energy $\E$ with respect to the set inclusion, noticing that $H_0^1(B_{2R}(y_{n_k})\cap \Omega_{n_k}; \Sigma_D)\subset H_0^1(B_{4R}(z_k)\cap \Omega_{n_k}; \Sigma_D)$, we have
\beq\label{eq:boundenergymonton}
\E(B_{2R}(y_{n_k})\cap \Omega_{n_k}; \Sigma_D)\geq \E(B_{4R}(z_k)\cap \Omega_{n_k}; \Sigma_D),
\eeq
for all $k\geq k_0$. Clearly, by construction, $|B_{4R}(z_k)\cap \Omega_{n_k}| \leq c$ for all $k\geq k_0$ and $|z_k|\to +\infty$, as $k\to +\infty$. We claim that, up to a further subsequence (still indexed by $k$) it holds
\beq\label{eq:claimStep4}
\E(B_{4R}(z_k)\cap \Omega_{n_k}; \Sigma_D) \geq \O_c\left(\Sigma_{\S_+^{N-1}}\right)+o(1),
\eeq
for all sufficiently large $k$, where $o(1)\to 0$ as $k\to +\infty$.\\ 

Notice that if Claim \eqref{eq:claimStep4} holds then the proof of \textbf{Step 4} is complete. Indeed combining \eqref{eq:thesisStep3}, \eqref{eq:boundenergymonton} and \eqref{eq:claimStep4}, up to a subsequence, we have
$$\E(\Omega_{n_k}; \Sigma_D)\geq  \O_c(\Sigma_{\S^{N-1}_+}) +o(1)-2c\sqrt{2\eps},$$
for all sufficiently large $k$. Then, passing to the limit as $k\to +\infty$ we get
$$\mathcal{O}_c(\Sigma_D)\geq \mathcal{O}_c\left(\Sigma_{\S^{N-1}_+}\right)-2c\sqrt{2\eps},$$
but this contradicts \eqref{eq20:teo1sect8}.\\

\noindent\textbf{Proof of Claim \eqref{eq:claimStep4}}:\\

We first observe that since $\partial\Sigma_D\setminus\{0\}$ is a smooth hypersurface, then, for any fixed $q\in\partial D \subset \partial\Sigma_D\setminus\{0\}$ there exitsts an open neighborhood $V$ of $q$ in $\partial\Sigma_D\setminus\{0\}$ such that $V-q$ is the graph over $T_q\partial\Sigma_D$ of a smooth function $g: U\to \R$, with $g(0)=0$, where $U$ is an open neighborhood of the origin in $T_q\partial\Sigma_D$ (without loss of generality we can assume that $U$ is a ball and $g$ is smooth in $\overline{U}$). Namely, fixing a orthonormal base $\mathcal{B}^\prime:=\{v_1,\ldots,v_{N-1}\}$ of $T_q\partial\Sigma_D$ and choosing $\mathcal{B}:=\{v_1,\ldots,v_{N-1},-\nu(q)\}$ as orthonormal base of $\RN$, where $-\nu(q)$ is the inner unit normal of $\partial\Sigma_D$ at $q$, denoting by $x^\prime=(x_1^\prime,\ldots,x^\prime_{N-1})$ the coordinates of the points in $T_q\partial\Sigma_D$ with respect to $\mathcal{B}^\prime$ and by $x=(x^\prime,x_N)$ the coordinates in $\RN$ with respect to $\mathcal{B}$, we can identify
$$V-q=\{(x^\prime,x_N) \in \RN; \ x^\prime=(x_1^\prime,\ldots,x_{N-1}^\prime)\in U, \ x_N=g(x_1^\prime,\ldots, x_{N-1}^\prime)\},$$
where $U$ is the orthogonal projection of $V-q$ onto $T_q\partial\Sigma_D$. To be precise, if $\varphi$ is a local parametrization centered at $q$, i.e. $\varphi(0)=q$, by writing $\varphi-q=\sum_{i=1}^{N-1} [(\varphi-q) \boldsymbol{\cdot} v_i] v_i +(\varphi-q) \boldsymbol{\cdot} (-\nu(q))$, and since $\sum_{i=1}^{N-1} [(\varphi-q) \boldsymbol{\cdot} v_i] v_i$ is a locally invertible map from an open neighborhood of the origin in $\R^{N-1}$ to an open neighborhood of the origin in $T_q\partial\Sigma_D\cong\R^{N-1}$, then, denoting by $G$ its local inverse and taking $g(x^\prime):=[(\varphi-q)\circ G(x^\prime)] \boldsymbol{\cdot} (-\nu(q))$ we obtain the desired map. In particular, notice that since $\frac{\partial[(\varphi-q)\circ G]}{\partial x^\prime_i}(0)\in T_q\partial\Sigma_D$ it follows that $\frac{\partial g}{\partial x_i^\prime}(0)=0$, for any $i=1,\ldots,N-1$.

Now, since $\partial\Sigma_D$ is a cone it follows that for any $t>0$, $T_{tq}\partial\Sigma_D=T_q\partial\Sigma_D$, $\nu(tq)=\nu(q)$ and $tV-tq$ is the graph over $T_q\partial\Sigma_D$ of the map $g_t:tU\to \R$ defined by 
\beq\label{eq:defgt}
g_t(x^\prime):=tg\left(\frac{x^\prime}{t}\right),\ x^\prime\in tU. 
\eeq
For any $x^\prime \in tU$, for any $i=1,\ldots,N-1$, we have
\beq\label{eq:svilgt}
 \frac{\partial g_t}{\partial x^\prime_i}(x^\prime)= \frac{\partial g}{\partial x^\prime_i}\left(\frac{x^\prime}{t}\right)=\frac{\partial g}{\partial x^\prime_i}\left(0\right) + \frac{1}{t} \nabla_{x^\prime} \left[\frac{\partial g}{\partial x^\prime_i}\right]\left(\frac{\xi}{t}\right) \boldsymbol{\cdot} x^\prime,
 \eeq
where $\xi=\xi(x^\prime,t)$ belongs to the segment joining $0$ and $\frac{x^\prime}{t}$, and $\nabla_{x^\prime}$ denotes the gradient with respect to $x_1^\prime,\ldots,x_{N-1}^\prime$. Hence, for any fixed ball $B_{R_1}(0)\subset T_q\partial\Sigma_D$, for all $t>0$ sufficiently large such that $B_{R_1}(0) \subset tU$, recalling that $\frac{\partial g}{\partial x^\prime_i}\left(0\right)=0$, we have 
\beq\label{eq:estimatenablag}
 \max_{x^\prime \in \overline{B_{R_1}(0)}}{|\nabla_{x^\prime} g_t(x^\prime)|}\leq \frac{1}{t}  \max_{x^\prime \in \overline{U}} \sqrt{\sum_{i=1}^{N-1}\left|\nabla_{x^\prime} \left[\frac{\partial g}{\partial x^\prime_i}\right](x^\prime)\right|^2}R_1\leq \frac{C}{t},
\eeq
where $C$ is independent of $t$, and in particular
\beq\label{eq:limgt}
 \lim_{t\to+\infty}  \max_{x^\prime \in \overline{B_{R_1}(0)}}{|\nabla_{x^\prime} g_t(x^\prime)|}=0.
\eeq
Let $C^+_{B_{R_1}(0)}$, $E^+\left(g_t\big|_{B_{R_1}(0)}\right)$ be, respectively, the upper cylinder generated by $B_{R_1}(0)$ and the epigraph associated to $g_t\big|_{B_{R_1}(0)}$, namely
\beq\label{eq:defepigraph}
\begin{array}{lll}
\displaystyle C^+_{B_{R_1}(0)}&:=&\displaystyle \{(x^\prime,x_N)\in \RN; \ x^\prime=(x_1^\prime,\ldots,x_{N-1}^\prime)\in B_{R_1}(0),\ x_N>0\},\\[6pt]
\displaystyle E^+\left(g_t\big|_{B_{R_1}(0)}\right)&:=&\displaystyle \{(x^\prime,x_N)\in \RN; \ x^\prime=(x_1^\prime,\ldots,x_{N-1}^\prime)\in B_{R_1}(0),\ x_N>g_t(x_1^\prime,\ldots,x_{N-1}^\prime)\}.
\end{array}
\eeq
Then, the map $F_t: \overline{C^+_{B_{R_1}(0)}} \to \overline{E^+\left(g_t\big|_{B_{R_1}(0)}\right)}$ defined by 
\beq\label{eq:defdiffFt}
F_t(x^\prime,x_N)=(x^\prime, x_N+g_t(x^\prime)), \ (x^\prime,x_N)\in  \overline{C^+_{B_{R_1}(0)}},
\eeq
is a diffeomorphism whose Jacobian matrix is of the form
\beq\label{eq:JacFt}
\mathrm{Jac}(F_t)(x^\prime,x_N)=\left[\begin{array}{cc}\mathbb{I}_{N-1}&0_{N-1}^T\\ \nabla_{x^\prime} g_t(x^\prime)&1 \end{array}\right],
\eeq
where $\mathbb{I}_{N-1}$ is the identity matrix of order $N-1$, $0_{N-1}$ is the null vector in $\R^{N-1}$, $T$ is the transposition. Notice also that $\mathrm{Jac}(F_t)$ is independent of $x_N$ and in view of \eqref{eq:limgt} it holds that
\beq\label{eq:limJacFt}
\lim_{t\to +\infty}\|\mathrm{Jac}(F_t)-\mathbb{I}_{N} \|_{C^0\left(\overline{C^+_{B_{R_1}(0)}}\right)}= 0.
\eeq
Now, let us consider the sequence $(q_k)_k\subset \partial D$, where $q_k:=\frac{z_k}{|z_k|}$, and $(z_k)_k\subset \partial\Sigma_D\setminus\{0\}$ is the sequence appearing in Claim \eqref{eq:claimStep4}. Since $\partial D$ is a compact subset of $\S^{N-1}$ we deduce that, up to a subsequence (still indexed by $k$) it holds that $\textrm{dist}_{\S^{N-1}}(q_k, \bar q)\to 0$, as $k\to +\infty$, for some $\bar q\in\partial D$, where $\mathrm{dist}_{\S^{N-1}}$ denotes the geodesic distance in $\S^{N-1}$. Then, from the previous discussion there exist an open neighborhood $V_1$ of $\bar q$ in $\partial\Sigma_D\setminus\{0\}$, a convex open neighborhood $U_1$ of the origin in $T_{\bar q}\partial\Sigma_D$ and a smooth function $g_1:\overline{U_1}\to \R$ such that $g_1(0)=0$, $\nabla_{x^\prime}g_1(0)=0$, and $V_1-\bar q$ is the graph over $T_{\bar q}\partial\Sigma_D$ associated to $g_1\big|_{U_1}$, where $x^\prime=(x_1^\prime,\ldots,x_{N-1}^\prime)$ are the coordinates with respect to fixed orthonormal base $\{\bar v_1,\ldots,\bar v_{N-1}\}$ of $T_{\bar q}\partial\Sigma_D$. Since $\textrm{dist}_{\S^{N-1}}(q_k, \bar q)\to 0$, as $k\to +\infty$, then definitely $q_k \in V_1$, $\Pi_{T_{\bar q}\partial\Sigma_D}(q_k-\bar q)\in U_1$, where $\Pi_{T_{\bar q}\partial\Sigma_D}:\R^N\to T_{\bar q}\partial\Sigma_D$ is the orthogonal projection onto $T_{\bar q}\partial\Sigma_D$ and $\Pi_{T_{\bar q}\partial\Sigma_D}(q_k-\bar q)\to 0$, as $k\to +\infty$. Let $\bar x_k^\prime:=(\bar x_{1,k}^\prime,\ldots,\bar x_{N-1,k}^\prime)$ be the coordinates of $\Pi_{T_{\bar q}\partial\Sigma_D}(q_k-\bar q)$ and set
\begin{eqnarray*}
 U_{1,k}&:=&U_1-\bar x_k^\prime,\\
 g_{1,k}(x^{\prime})&:=&g_1(x^{\prime}+\bar x^\prime_k)-g_1(\bar x_{k}^\prime), \ x^\prime \in U_{1,k}.
\end{eqnarray*} 
Then we readily check that $V_1-q_k$ is a cartesian graph over $T_{\bar q}\partial\Sigma_D$, associated to $g_{1,k}:U_{1,k}\to \R$. Notice that, since $\bar x_k^\prime\to 0$, as $k\to+\infty$, then there exists a ball $B_{\bar R}(0)$ in $T_{\bar q}\partial\Sigma_D$ such that $B_{\bar R}(0)\subset U_{1,k}$ for all sufficiently large $k$. In particular, setting $t_k:=|z_k|$, recalling that $|z_k|\to +\infty$, as $k\to +\infty$, then, $ t_kU_{1,k}$ invades $T_{\bar q}\partial\Sigma_D$. As in \eqref{eq:defgt} we consider the rescaled map $h_{t_k}:t_kU_{1,k}\to\R$ defined by
$$h_{t_k}(x^{\prime})=t_k g_{1,k}\left(\frac{x^{\prime}}{t_k}\right)=t_k\left[g_1\left(\frac{x^{\prime}}{t_k}+\bar x_k^\prime\right)-g_1(\bar x_{k}^\prime)\right], \ x^{\prime}\in t_k U_{1,k},$$
and arguing as in \eqref{eq:svilgt} we have the expansion
\beq\label{eq:svilgt2}
\frac{\partial h_{t_k}}{\partial x^{\prime}_i}(x^{\prime})=\frac{\partial g_1}{\partial x^\prime_i}\left(\frac{x^{\prime}}{t_k}+\bar x_k^\prime\right)=\frac{\partial g_1}{\partial x^\prime_i}\left(\bar x_k^\prime\right) + \frac{1}{t_k} \left[\nabla_{x^\prime} \frac{\partial g_1}{\partial x^\prime_i}\left(\frac{\xi_k}{t_k}+\bar x_k^\prime\right) \right]\boldsymbol{\cdot} x^{\prime},
\eeq
where $\xi_k$ belongs to the segment joining $\bar x_k^\prime$ and $\frac{x^{\prime}}{t_k}$. Let us fix a ball $B_{R_1}(0)$ in $T_{\bar q}\partial\Sigma_D$, with $R_1$ to be chosen later and independently on $k$, and observe that $B_{R_1}(0)\subset t_kU_{1,k}$ for all sufficiently large $k$. Since $\bar x_k^\prime \to 0$, as $k\to +\infty$, and $\nabla_{x^\prime} g_1(0)=0$ we get that the first term in the right-hand side of \eqref{eq:svilgt2} goes to zero as $k\to +\infty$, and arguing as in \eqref{eq:estimatenablag} for the second term, we conclude that
\beq\label{eq:limthtk}
\lim_{k\to +\infty} \max_{x^{\prime} \in \overline{B_{R_1}(0)}}{|\nabla_{x^{\prime}} h_{t_k}(x^{\prime})|} = 0.
\eeq
Let $F_{t_k}:  \overline{C^+_{B_{R_1}(0)}} \to \overline{E^+\left(h_{t_k}\big|_{B_{R_1}(0)}\right)}$ be the diffeomorphism  defined by $$F_{t_k}(x^\prime,x_N)=(x^\prime,x_N+h_{t_k}(x^\prime)), \ (x^\prime,x_N)\in  \overline{C^+_{B_{R_1}(0)}},$$ 
where $E^+\left(h_{t_k}\big|_{B_{R_1}(0)}\right)$ is the epigraph associated to $h_{t_k}\big|_{B_{R_1}(0)}$ (see \eqref{eq:defepigraph}) and where we recall that $x=(x^\prime,x_N)$ are the coordinates with respect to the orthogonal base $\{\bar v_1,\ldots,\bar v_{N-1},-\nu(\bar q)\}$.

Now, let us consider the set $B_{4R}(z_k)\cap \Omega_{n_k}$ appearing in \eqref{eq:claimStep4}. Recalling that $z_k=t_k q_k$, since $t_k\to +\infty$, as $k\to+\infty$, and since $R$ is independent of $k$, then, for all sufficiently large $k$ it holds $$\overline{(B_{4R}(z_k)\cap \Omega_{n_k})} \cap \partial\Sigma_D-z_k\subset t_k(V_1-q_k).$$ 
We observe that for any $k$ we have
\beq\label{eq:simpleremark}
\E(B_{4R}(z_k)\cap \Omega_{n_k}; \Sigma_D)=\E(B_{4R}(z_k)\cap \Omega_{n_k}-z_k; \Sigma_D-z_k)
\eeq
and we set for brevity 
\beq\label{eq:simpleremark2}
\widetilde\Omega_{R,k}:=(B_{4R}(z_k)\cap \Omega_{n_k})-z_k.
\eeq
Notice that since ${\widetilde\Omega_{R,k}}$ is a uniformly bounded subset of $\RN$ and $t_k(V_1-q_k)$ is the cartesian graph of $h_{t_k}:t_kU_{1,k}\to\R$ then we can choose $R_1>0$ (independent of $k$) in such a way that \beq\label{eq:inclusionBCil}
F_{t_k}^{-1}\left(\ \overline{\widetilde\Omega_{R,k}}\ \right)\subset B_{R_1}(0)\times [0,+\infty[,\eeq
for all large $k$. Let us denote by $\widetilde u_{k}:= u_{\widetilde\Omega_{R,k}} \in H_0^1(\widetilde\Omega_{R,k}; \Sigma_D-z_k)$ the energy function of $\widetilde\Omega_{R,k}$ and let $\widetilde U_k:=\widetilde u_k\circ F_{t_k}$. Notice that  by construction and thanks to \eqref{eq:inclusionBCil}, it follows that $\tilde U_k$ extends to a function $\tilde U_k \in H_0^1(F_{t_k}^{-1}(\tilde\Omega_{R,k}); \R^N_+)$.

 Thanks to $\eqref{eq:JacFt}$ we have that $\textrm{det}(\textrm{Jac}(F_{t_k}))\equiv 1$ and thus
\beq\label{eq:transfL1norm}
\int_{\widetilde\Omega_{R,k}} \widetilde u_k \ dy=\int_{F_{t_k}^{-1}(\widetilde\Omega_{R,k})} \widetilde u_k\circ F_{t_k}\ \left|\textrm{det}(\textrm{Jac}(F_{t_k}))\right| \ dx =\int_{F_{t_k}^{-1}(\widetilde\Omega_{R,k})} \widetilde U_k \ dx,
\eeq
and 
\beq\label{eq:boundmeasFtk}
|F_{t_k}^{-1}(\widetilde\Omega_{R,k})|=|\widetilde\Omega_{R,k}|\leq c.
\eeq
Moreover, setting $$M_{1,k}:=\max_{(x^\prime,x_N)\in \overline{C^+_{B_{R_1}(0)}}} \max_{|\eta|\leq 1} |[\textrm{Jac}(F_{t_k})(x^\prime,x_N)] \eta |,$$
which is well defined (see \eqref{eq:JacFt}) and taking into account that
$$\nabla_x \widetilde U_k(x^\prime,x_N)=(\textrm{Jac}(F_{t_k})(x^\prime,x_N))^T \nabla_{y} \widetilde u_k(F_{t_k}(x^\prime,x_N)) \ \ \hbox{for a.e. $(x^\prime,x_N) \in F_{t_k}^{-1}(\widetilde\Omega_{R,k})$}, $$
we obtain that
\beq\label{eq:estimateF_tktildeU} 
\begin{array}{lll}
\displaystyle \int_{F_{t_k}^{-1}(\tilde\Omega_{R,k})} |\nabla_{x} \tilde U_k|^2 \ dx&\leq&\displaystyle M_{1,k}^2 \int_{F_{t_k}^{-1}(\tilde\Omega_{R,k})}|\nabla_{y} \widetilde u_k\circ F_{t_k}|^2 \ dx\\[12pt]
&=&\displaystyle M_{1,k}^2 \int_{\tilde\Omega_{R,k}}|\nabla_{y}  \widetilde u_k|^2 \ dy.
\end{array}
\eeq
In view of \eqref{eq:limthtk} and recalling \eqref{eq:JacFt}, \eqref{eq:limJacFt} we deduce that  $M_{1,k}\to 1$ as $k\to +\infty$. 
Therefore, recalling that $\widetilde u_k$ is the energy function of $\widetilde\Omega_{R,k}=B_{4R}(z_k)\cap \Omega_{n_k}-z_k$ and thanks to \eqref{eq:estimateF_tktildeU} we get that $ \int_{F_{t_k}^{-1}(\tilde\Omega_{R,k})} |\nabla_{x} \tilde U_k|^2 \ dx$ is bounded by a uniform positive constant. Moreover, by definition and thanks to \eqref{eq:transfL1norm}, \eqref{eq:estimateF_tktildeU} we have
\beq\label{eq:estimateenergyOmegaRk}
\begin{array}{lll}
\displaystyle\E(\tilde\Omega_{R,k}; \Sigma_D-z_k)&=&\displaystyle\frac{1}{2}\int_{\tilde\Omega_{R,k}}|\nabla_{y}   \widetilde u_k|^2 \ dy -\int_{\widetilde\Omega_{R,k}} \widetilde u_k \ dy\\[12pt]
&\geq&\displaystyle\frac{1}{2 } \frac{1}{M_{1,k}^2} \int_{F_{t_k}^{-1}(\tilde\Omega_{R,k})} |\nabla_{x} \tilde U_k|^2 \ dx  -\int_{F_{t_k}^{-1}(\tilde\Omega_{R,k})}  \tilde U_k \ dx\\[12pt]
&=&\displaystyle J( \tilde U_k) + \displaystyle\frac{1}{2}\left(\frac{1}{M_{1,k}^2}-1\right) \int_{F_{t_k}^{-1}(\tilde\Omega_{R,k})} |\nabla_{x} \tilde U_k|^2 \ dx.
\end{array}
\eeq
Hence, as $\tilde U_k \in H_0^1(F_{t_k}^{-1}(\tilde\Omega_{R,k}); \R^N_+)$, denoting by $W_k:=u_{F_{t_k}^{-1}(\tilde\Omega_{R,k})} \in H_0^1(F_{t_k}^{-1}(\tilde\Omega_{R,k}); \R^N_+)$ the energy function of $F_{t_k}^{-1}(\tilde\Omega_{R,k})$ then by reflection, symmetrization and taking into account \eqref{eq:boundmeasFtk} and \eqref{eq:exprtorsionhemisphereR} we infer that
\beq\label{eq:estimateenergyOmegaRk2}
\begin{array}{lll}
\E(\tilde\Omega_{R,k}; \Sigma_D-z_k)&\geq&\displaystyle J(W_k) + \displaystyle\frac{1}{2}\left(\frac{1}{M_{1,k}^2}-1\right) \int_{F_{t_k}^{-1}(\tilde\Omega_{R,k})} |\nabla_{x} \tilde U_k|^2 \ dx\\[12pt]
&\geq&\displaystyle \O_{c_k}(\Sigma_{\S^{N-1}_+}) + \displaystyle\frac{1}{2}\left(\frac{1}{M_{1,k}^2}-1\right) \int_{F_{t_k}^{-1}(\tilde\Omega_{R,k})} |\nabla_{x} \tilde U_k|^2 \ dx\\[12pt]
&\geq&\displaystyle \O_{c}(\Sigma_{\S^{N-1}_+}) + \displaystyle\frac{1}{2}\left(\frac{1}{M_{1,k}^2}-1\right) \int_{F_{t_k}^{-1}(\tilde\Omega_{R,k})} |\nabla_{x} \tilde U_k|^2 \ dx.
\end{array}
\eeq
Finally, since $\int_{F_{t_k}^{-1}(\tilde\Omega_{R,k})} |\nabla_{x} \tilde U_k|^2 \ dx$ is uniformly bonded and $M_{1,k}\to 1$, as $k\to +\infty$, then from \eqref{eq:simpleremark}, \eqref{eq:simpleremark2} and \eqref{eq:estimateenergyOmegaRk2}, we get 
$$\E(B_{4R}(z_k)\cap \Omega_{n_k}; \Sigma_D)\geq  \O_c(\Sigma_{\S^{N-1}_+}) +o(1),$$
for all sufficiently large $k$, and this proves Claim \eqref{eq:claimStep4}.\\


In the next step, we prove the pre-compactness of the sequence $(u_{n_k})_k$ in $L^2(\Sigma_D)$.\\

\noindent\textbf{Step 5:} The sequence $(u_{n_k})_k$ admits a subsequence which strongly converges in $L^2(\Sigma_D)$.\\

We first show that
\beq\label{eq7:teo1sect8}
\lim_{R\to+\infty} \sup_{k \in \N} \int_{B_R^\complement(0)\cap\Sigma_D} u_{n_k}^2 \ dx=0
\eeq
 Indeed, if \eqref{eq7:teo1sect8} is not true there exist $\eps^\prime>0$, a sequence $(R_m)_m\subset \R^+$ such that $R_m\to+\infty$, as $m\to +\infty$, and we find a subsequence $(n_{k_m})_m$ such that for all $m\in \N$
 \beq\label{eq10:teo1sect8}
  \int_{B_{R_m}^\complement(0)\cap\Sigma_D} u_{n_{k_m}}^2 \ dx\geq \frac{\eps^\prime}{2}.
 \eeq
On the other hand, taking $\eps=\frac{\eps^\prime}{4}$ in (i) we find $R^\prime>0$ depending only on $\eps^\prime$ such that for all $k \in \N$
\beq\label{eq8:teo1sect8}
 \int_{B_{R^\prime}(y_k)\cap\Sigma_D} u_{n_{k}}^2 \ dx\geq \lambda-\frac{\eps^\prime}{4}.
\eeq
Now, in view of \textbf{Step 4} we know that $(y_{n_k})_k$ is bounded, and thus there exists $R^{\prime\prime}>0$ independent of $k$ such that $B_{R^\prime}(y_k)\subset B_{R^{\prime\prime}}(0)$ for all $k$. Hence, from \eqref{eq8:teo1sect8} we get that
\beq\label{eq9:teo1sect8}
 \int_{B_{R^{\prime\prime}(0)}\cap\Sigma_D} u_{n_{k}}^2 \ dx \geq \int_{B_{R^\prime}(y_k)\cap\Sigma_D} u_{n_{k}}^2 \ dx\geq \lambda-\frac{\eps^\prime}{4},
\eeq
for all sufficiently large $k$. Finally, by writing
$$
\int_{\Sigma_D} u_{n_{k_m}}^2 \ dx = \int_{B_{R_m}(0)\cap\Sigma_D} u_{n_{k_m}}^2 \ dx+ \int_{B_{R_m}^\complement(0)\cap\Sigma_D} u_{n_{k_m}}^2 \ dx,
$$
and recalling that $R_m\to +\infty$, then, we have $R_m>R^{\prime\prime}$, for all sufficiently large $m$, and thus from \eqref{eq10:teo1sect8}, \eqref{eq9:teo1sect8} we deduce that
$$ \int_{\Sigma_D} u_{n_{k_m}}^2 \ dx \geq \lambda-\frac{\eps^\prime}{4} + \frac{\eps^\prime}{2} = \lambda+\frac{\eps^\prime}{4},$$
for all sufficiently large $m$, but this contradicts \eqref{eq0:teo1sect8} and \eqref{eq7:teo1sect8} is thus proved.\\

In order to prove the relative compactness of the sequence $(u_{n_k})_k$ in $L^2(\Sigma_D)$, it suffices to find, for any given $\eps>0$, a relative compact sequence $(v_k)_k$ in $L^2(\Sigma_D)$, depending on $\eps$, with the property that 
\begin{equation}
  \label{eq:sufficient-crit-compactness}
\|u_{n_k}-v_k\|_{L^2(\Sigma_D)}< \eps  \qquad \text{$\forall k \in \N$.}
\end{equation}
Indeed, the latter property readily implies that the set $\{u_{n_k};\ k \in \N \}$ is totally bounded in $L^2(\Sigma_D)$, and therefore it is relative compact since $L^2(\Sigma_D)$ is a Banach space. So let $\eps>0$. By \eqref{eq7:teo1sect8}, there exists $R>0$ with 
$$
\int_{B_R^\complement(0)\cap\Sigma_D} u_{n_k}^2 \ dx< \eps \qquad \text{$\forall k \in \N$.}
$$
Hence \eqref{eq:sufficient-crit-compactness} holds with $v_k: = \chi_{B_R(0)}u_{n_k}$ for $k \in \N$, where $\chi_{B_R(0)}$ denotes the characteristic function of the ball $B_R(0)$. Moreover, by \eqref{eq:boundH1teosect8} and the compactness of the embedding $H^1(B_R(0)\cap \Sigma_D) \hookrightarrow L^2(B_R(0)\cap \Sigma_D)$ the sequence of functions $u_{n_k}\big|_{B_R(0)}$, $k \in \N$ is relatively compact in $L^2(B_R(0)\cap \Sigma_D)$, which obviously implies that the sequence $(v_k)_k$ is relatively compact in $L^2(\Sigma_D)$, as required. We have thus established the relative compactness of the sequence $(u_{n_k})_k$ in $L^2(\Sigma_D)$, as claimed.\\

\noindent\textbf{Step 6:} Existence of a minimizer for $\mathcal{O}_c(\Sigma_D)$.\\

In the previous steps we proved that the sequence of energy functions $(u_n)_n$, associated to a minimizing sequence $(\Omega_n)_n \subset \Sigma_D$ for $\mathcal{O}_c(\Sigma_D)$, is bounded in $H^1(\Sigma_D)$ (see \eqref{eq:boundH1teosect8}) and possesses a subsequence which strongly converges in $L^2(\Sigma_D)$. Hence, up to a subsequence (still indexed by $n$ for convenience), we have $u_n \rightharpoonup \bar u$ in $H^1(\Sigma_D)$, for some $\bar u \in H^1(\Sigma_D)$, and $u_n \to \bar u$ in $L^2(\Sigma_D)$, as $n\to +\infty$.

We set $\Omega:=\{\bar u>0\} \subset \Sigma_D$. Since $\bar u\in H^1(\Sigma_D)$ then $\Omega$ is a quasi-open subset of $\Sigma_D$, in addition, arguing as in \cite[Proof of Lemma 5.2]{BV1}, namely using that $u_n \to \bar u$ in $L^2(\Sigma_D)$, as $n\to +\infty$, and applying Fatou's Lemma, we infer that 
$$|\Omega|=\int_{\Sigma_D} \chi_{\{\bar u>0\}} \ dx \leq \liminf_{n\to +\infty} \int_{\Sigma_D} \chi_{\{ u_n>0\}} \ dx =  \liminf_{n\to +\infty} |\Omega_n| \leq c. $$
We claim that $\Omega$ is a minimizer for $\mathcal{O}_c(\Sigma_D)$ and that $\bar u$ is the torsion function of $u_\Omega$. To prove this we first observe that as $u_n \to \bar u$ in $L^2(\Sigma_D)$ and since $|\Omega_n|\leq c$, $|\Omega|\leq c$ it follows that $u_n \to \bar u$ in $L^1(\Sigma_D)$. Indeed by construction we have $u_n\in H_0^1(\Omega_n;\Sigma_D)$, $\bar u \in H_0^1(\Omega;\Sigma_D)$, and by H\"older's inequality we deduce that
$$\int_{\Sigma_D}|u_n-\bar u|  \ dx \leq \left(\int_{\Omega_n \cup \Omega}|u_n-\bar u|^2  \ dx\right)^{\frac{1}{2}} |\Omega_n\cup \Omega|^{\frac{1}{2}}\leq \|u_n-\bar u\|_{L^2(\Sigma_D)} \sqrt{2c}.$$
Now, as $u_n \rightharpoonup \bar u$ in $H^1(\Sigma_D)$, we have
$$\|\bar u\|_{H^1(\Sigma_D)}^2 \leq \liminf_{n\to+\infty} \| u_n\|_{H^1(\Sigma_D)}^2,$$
and thus, since $u_n \to \bar u$ in $L^2(\Sigma_D)$, we readily get that
$$ \int_{\Sigma_D} |\nabla \bar u|^2 \ dx  \leq \liminf_{n\to+\infty} \int_{\Sigma_D} |\nabla u_n|^2 \ dx$$
Then, recalling the definition of the functional $J$ (see \eqref{eq:deftorsion}), exploiting that $u_n \to \bar u$ in $L^1(\Sigma_D)$ and since $u_n$ are the energy functions associated to $\Omega_n$, we obtain that
\beq\label{eq14:teo1sect8}
J(\bar u)\leq \liminf_{n\to+\infty} \left(\frac{1}{2} \int_{\Sigma_D} |\nabla u_n|^2 \ dx - \int_{\Sigma_D} u_n \ dx \right)=\liminf_{n\to+\infty} \E(\Omega_n; \Sigma_D) =\mathcal{O}_c(\Sigma_D).
\eeq
Finally, considering the energy function $u_{\Omega}$ associated to $\Omega$, i.e. the minimizer of $J$ in $H_0^1(\Omega;\Sigma_D)$, then, by the minimality of $u_{\Omega}$, since $\bar u \in H_0^1(\Omega;\Sigma_D)$ and thanks to \eqref{eq14:teo1sect8} we have
\beq\label{eq:energyineqeq}
\E(\Omega; \Sigma_D)=J(u_{\Omega})\leq J(\bar u) \leq \mathcal{O}_c(\Sigma_D).
\eeq
Therefore $\E(\Omega; \Sigma_D)=\mathcal{O}_c(\Sigma_D)$, and \eqref{eq:energyineqeq} implies that $J(u_\Omega)=J(\bar u)$. Hence $\Omega$ is a minimizer for $\mathcal{O}_c(\Sigma_D)$ and $\bar u=u_\Omega$ in $H^1(\Sigma_D)$.
\end{proof}

\begin{corollary}\label{cor:mainteoexistmin}
If $D\subset \S^{N-1}$ is a smooth domain such that 
\beq\label{eq:corollarythesis}
\mathcal{H}_{N-1}(D) < \mathcal{H}_{N-1}(\S^{N-1}_+)
\eeq
then $\O_c(\Sigma_D)$ is achieved, for any $c>0$.
\end{corollary}
\begin{proof}
By \eqref{def:radialfundamentalsoltorsion}--\eqref{eq:exprtorsionhemisphereR} we readily check that \eqref{eq:corollarythesis} implies the condition \eqref{eq:conditionD}, and by Theorem \ref{teo:mainteoexistmin} we conclude.
\end{proof}

We conclude this section with the following 
\begin{proposition}\label{prop:ineqtesi}
Let $D\subset \S^{N-1}$ be a smooth domain and let $c>0$. Then
\beq\label{eq:tesipropsect72}
\O_c\left(\Sigma_D\right) \leq \O_c\left(\Sigma_{\S^{N-1}_+}\right).
\eeq
\end{proposition}
\begin{proof}
Let us fix $q\in\partial D\subset\partial\Sigma_D\setminus\{0\}$ and let $\{v_1,\ldots,v_{N-1}\}$ be an orthonormal basis of $T_q\partial\Sigma_D$. We denote by $x=(x^\prime,x_N)$ the coordinates of points in $\RN$ with respect to $\{v_1,\ldots,v_{N-1},-\nu(q)\}$, where $-\nu(q)$ is the inner unit normal to $\partial\Sigma_D$ at $q$.
As seen in the proof of Claim \eqref{eq:claimStep4} there exist an open neighborhood $V$ of $q$ in $\partial\Sigma_D\setminus\{0\}$, an open neighborhood $U$ of the origin in $T_q\partial\Sigma_D$, and a smooth map $g:\overline{U}\to \R$, $g=g(x^\prime)$ such that $V-q$ is the graph over $T_q\partial\Sigma_D$ of $g\big|_{U}$.
 
Let  $B^+_{R}(0)\subset \R^N_+$ be a $N$-dimensional half-ball such that $|B^+_{R}(0)|=c$, i.e. $B^+_{R}(0)$ is a half-ball of volume $c$ contained in the upper half-space delimited by $T_q\partial\Sigma_D$. Let $u_{B^+_{R}(0)}\in H_0^1(B^+_{R}(0); \R^N_+)$ be the energy function of $B^+_{R}(0)$. Then, by definition and recalling Remark \ref{rem:energyhalfball}, we have 
\beq\label{eq:EnergyhalfR}
\O_c(\Sigma_{\S^{N-1}_+})=\E(B^+_{R}(0); \R^N_+)=J(u_{B^+_{R}(0)}).
\eeq
Let $B_{R_1}(0)$ be a ball in $T_q\partial\Sigma_D$, with $R_1>R$. Clearly
\beq\label{eq:contBHRBR1}
\overline{B^+_R(0)} \subset B_{R_1}(0)\times[0,+\infty[.
\eeq
Let $(t_k)_k\subset \R^+$ be a sequence such that $t_k\to +\infty$, as $k\to +\infty$, then, setting $z_k:=t_k q$ we obtain a diverging sequence of points on $\partial\Sigma_D\setminus\{0\}$. We consider the rescaled map $g_{t_k}:t_k U\to \R$ defined by \eqref{eq:defgt} and the associated diffeomorphism $F_{t_k}: \overline{C^+_{B_{R_1}(0)}} \to \overline{E^+\left(g_{t_k}\big|_{B_{R_1}(0)}\right)}$ given by \eqref{eq:defdiffFt}, where $E^+\left(g_{t_k}\big|_{B_{R_1}(0)}\right)$, $C^+_{B_{R_1}(0)}$ are defined by \eqref{eq:defepigraph}. The inverse diffeomorphism $F_{t_k}^{-1}$ is given by
\beq\label{eq:defdiffinvFt}
F_{t_k}^{-1}(x^\prime,x_N)=(x^\prime, x_N-g_{t_k}(x^\prime)), \ (x^\prime,x_N)\in  \overline{E^+\left(g_{t_k}\big|_{B_{R_1}(0)}\right)},
\eeq
and as done in \eqref{eq:JacFt}, \eqref{eq:limJacFt} we readily check that 
\beq\label{eq:detFinv1}
\mathrm{Det}(Jac(F_{t_k}^{-1}))\equiv 1,\ \lim_{k\to +\infty}\|\mathrm{Jac}(F^{-1}_{t_k})-\mathbb{I}_{N} \|_{C^0\left(\overline{E^+\left(g_{t_k}\big|_{B_{R_1}(0)}\right)}\right)}= 0.
\eeq
Moreover, setting $U_k:=u_{B^+_R(0)}\circ F_{t_k}^{-1}$ we notice that since $u_{B^+_R(0)}\in H_0^1(B^+_{R}(0); \R^N_+)$ (actually $u_{B^+_R(0)}= 0$ in $\R^N_+\setminus B^+_{R}(0)$, see \eqref{def:radialfundamentalsoltorsion} with $D=\S^{N-1}_+$), then, by construction, taking into account \eqref{eq:contBHRBR1}, it follows that $U_k$ extends to a function $U_k\in H_0^1(F_{t_k}\left(B^+_{R}(0)\right); \Sigma_D-z_k)$. Arguing as in \eqref{eq:transfL1norm}--\eqref{eq:estimateF_tktildeU}, taking into account \eqref{eq:detFinv1}, we infer that $|F_{t_k}\left(B^+_{R}(0)\right)|=|B^+_{R}(0)|=c$, 
\beq\label{eq:estimateF_tm1UkL1} 
\int_{F_{t_k}\left(B^+_{R}(0)\right)} U_k \ dx=\int_{B^+_{R}(0)} u_{B^+_{R}(0)} \ dy,
\eeq
\beq\label{eq:estimateF_tm1Uk} 
\displaystyle \int_{F_{t_k}\left(B^+_{R}(0)\right)} |\nabla_{x} U_k|^2 \ dx\leq \displaystyle M_{2,k}^2 \int_{B^+_{R}(0)}|\nabla_{y}  u_{B^+_R(0)}|^2 \ dy,
\eeq
where $$M_{2,k}:=\max_{(x^\prime,x_N)\in \overline{E^+\left(g_{t_k}\big|_{B_{R_1}(0)}\right)}} \max_{|\eta|\leq 1} |[\textrm{Jac}(F^{-1}_{t_k})(x^\prime,x_N)] \eta |,$$
and $M_{2,k}\to 1$, as $k\to +\infty$. Hence, combining \eqref{eq:EnergyhalfR}, \eqref{eq:estimateF_tm1UkL1} and \eqref{eq:estimateF_tm1Uk} we deduce

\beq\label{eq:stimasfinalpropsect72}
\begin{array}{lll}
\O_c(\Sigma_{\S^{N-1}_+})&=&\displaystyle\frac{1}{2}\int_{B^+_{R}(0)}|\nabla_{y}   u_{B^+_R(0)}|^2 \ dy -\int_{B^+_{R}(0)}  u_{B^+_R(0)} \ dy\\[12pt]
&\geq&\displaystyle \frac{1}{2}\frac{1}{M_{2,k}^2} \int_{F_{t_k}\left(B^+_{R}(0)\right)} |\nabla_{x} U_k|^2 \ dx  -\int_{F_{t_k}\left(B^+_{R}(0)\right)} U_k \ dx \\[12pt]
&=&\displaystyle J(U_k) + \displaystyle\frac{1}{2}\left(\frac{1}{M_{2,k}^2}-1\right)\int_{F_{t_k}\left(B^+_{R}(0)\right)} |\nabla_{x} U_k|^2 \ dx \\[12pt]
&\geq&\displaystyle \E\left(F_{t_k}\left(B^+_{R}(0)\right)\right) + o(1)
\end{array}
\eeq
where in the last inequality we used that $U_k\in H_0^1(F_{t_k}\left(B^+_{R}(0)\right); \Sigma_D-z_k)$, the definition of torsional energy of $F_{t_k}\left(B^+_{R}(0)\right)$, $M_{2,k}\to 1$, as $k\to +\infty$ and that $\int_{F_{t_k}\left(B^+_{R}(0)\right)} |\nabla_{x} U_k|^2 \ dx$ is uniformly bounded. Summing up, from \eqref{eq:stimasfinalpropsect72} and the definition of $\O_c(\Sigma_D)$ we finally have
$$\O_c(\Sigma_{\S^{N-1}_+}) \geq \E(F_{t_k}\left(B^+_{R}(0)\right); \Sigma_D-z_k) + o(1)= \E(F_{t_k}\left(B^+_{R}(0)\right)+z_k; \Sigma_D)+o(1)\geq \O_c(\Sigma_D) + o(1),$$
and passing to the limit as $k\to +\infty$ we obtain \eqref{eq:tesipropsect72}.
\end{proof}

 \section{Properties of minimizers and proof of Theorem \ref{mainteo}}

In this section we show some qualitative properties of the minimizers of the torsional energy functional with fixed volume (we refer to Sect. 6 for the notations). In view of the scaling invariance of our problem (see Remark \ref{remark1:sect8}) it suffices to focus on the case $\O_1(\Sigma_D)$. We begin by proving that any minimizer for $\O_1(\Sigma_D)$ is bounded.

\begin{proposition}
If $\Omega$ is a minimizer for $\O_1(\Sigma_D)$ then $\Omega$ is bounded.
\end{proposition}
\begin{proof}

We argue as in \cite[Sect. 2.1.2]{BHS} with slightly changes.
 Let $\Omega$ be a minimizer for $\O_1(\Sigma_D)$. In view of \cite[Theorem 1]{BU2}, in order to prove that $\Omega$ is bounded it is sufficient to show that $\Omega$ is a local shape subsolution for the energy $\T_D$, which means that there exist $\delta>0$ and $\Lambda>0$ such that for any quasi-open subset $\widetilde\Omega\subset \Omega$ with
$\|u_{\Omega}-u_{\widetilde\Omega}\|_{L^2{\Sigma_D}}< \delta$ it holds
$$ \E(\Omega; \Sigma_D)+\Lambda |\Omega| \leq  \E(\widetilde\Omega; \Sigma_D)+\Lambda |\widetilde\Omega|.$$
Let us assume, by contradiction, that there exist a sequence $(\Lambda_n)_n\subset \R^+$ with  $\Lambda_n\to 0$, as $n\to +\infty$, and an increasing sequence $(\widetilde \Omega_n)_n\subset \Omega$ of quasi-open subsets such that
\beq\label{eq1:sect9}
\E(\Omega; \Sigma_D)+\Lambda_n |\Omega| > \E(\widetilde\Omega_n; \Sigma_D)+\Lambda_n |\widetilde\Omega_n|,
\eeq
and $\|u_{\Omega}-u_{\widetilde \Omega_n}\|_{L^2{\Sigma_D}}\to 0$, as $n\to +\infty$. Then, let us fix $t_n>1$ such that
$$|t_n \widetilde\Omega_n|=t_n^N|\widetilde\Omega_n|=1= |\Omega|. $$
Obviously, $t_n\to 1^+$, as $n\to +\infty$, and by the minimality of $\Omega$ we have
$$  t_n^{N+2}\E(\widetilde\Omega_n; \Sigma_D)= \E(t_n\widetilde\Omega_n; \Sigma_D)\geq \E(\Omega; \Sigma_D).$$
Thus, from \eqref{eq1:sect9} we obtain
$$ \Lambda_n (|\Omega|-|\widetilde\Omega_n|) > \E(\widetilde\Omega_n; \Sigma_D) - \E(\Omega; \Sigma_D) \geq \frac{(t_n^{N+2}-1)}{t_n^{N+2}} (-\E(\Omega; \Sigma_D)),$$
and dividing by $t_n^N-1=\frac{|\Omega|}{|\widetilde\Omega_n|}-1$ we get
\beq\label{eq2:sect9}
\frac{-\E(\Omega; \Sigma_D)}{t_n^{N+2}} \frac{(t_n^{N+2}-1)}{t_n^{N}-1} <  \Lambda_n \frac{(|\Omega|-|\widetilde\Omega_n|)}{t_n^N-1}= \Lambda_n |\widetilde\Omega_n|.
\eeq
This gives a contradiction because, as $n\to+\infty$, the left-hand side of \eqref{eq2:sect9} converges to $-\frac{N+2}{N}\E(\Omega; \Sigma_D)$, while the right-hand side converges to zero.
\end{proof}

\begin{proposition}
If $\Omega$ is a minimizer for $\O_1(\Sigma_D)$ then the torsion function $u_\Omega\in H_0^1(\Omega; \Sigma_D)$ is Lipschitz continuous in any Lipschitz  domain $\omega \subset \Sigma_D$ such that $\overline\omega\subset\Sigma_D$  and $\Omega=\{u_\Omega>0\}$ is an open subset of $\Sigma_D$.
\end{proposition}
\begin{proof}
The result follows essentially as in the case of Dirichlet boundary conditions, which was addressed in the work \cite{BRP}. The extension to the case of mixed boundary conditions was done in \cite[Theorem 2.14]{MPV} for the problem of minimizing the first eigenvalue. In our situation the proof would be similar.
\end{proof}

Next, we prove that any minimizer is connected.
\begin{proposition}
If $\Omega$ is a minimizer for $\O_1(\Sigma_D)$ then $\Omega$ is a connected subset of $\Sigma_D$.
\end{proposition}
\begin{proof}
As seen in the proof of Theorem \ref{teo:mainteoexistmin}, \textbf{Step 2}, we notice that, as $|\Omega|=1$, then $\Omega$ is also a minimizer for 
$$\mathcal{M}(\Sigma_D):=\inf\left\{\frac{\E(\Omega; \Sigma_D)}{|\Omega|^{\frac{N+2}{N}}}; \ \  \Omega \ \hbox{quasi-open}, \ \Omega\subset \Sigma_D,\ |\Omega|>0\right\}.$$
Assume by contradiction that $\Omega$ is not connected. Then there exist two open subsets $\Omega_1$, $\Omega_2$ of $\Sigma_D$, with $\Omega_1, \Omega_2\neq\emptyset$, such that $\Omega_1\cap\Omega_2=\emptyset$ and
$$\Omega= \Omega_1\cup\Omega_2.$$
In addition, since $\Omega_i$ is an open nonempty subset of $\Sigma_D$ then $|\Omega_i|>0$, for $i=1,2$, and by construction we have
$$\frac{\E(\Omega_i; \Sigma_D)}{|\Omega_i|^{\frac{N+2}{N}}} \geq \mathcal{M}(\Sigma_D), \ \hbox{for $i=1,2$, and}\ \ \E(\Omega; \Sigma_D)=\E(\Omega_1; \Sigma_D) + \E(\Omega_2; \Sigma_D), \ \ |\Omega|=|\Omega_1|+|\Omega_2|.$$
Then, since $\frac{N+2}{N}>1$, by the convexity of $t\mapsto t^{\frac{N+2}{N}}$, taking into account that $|\Omega_1|>0$, $|\Omega_2|>0$ and $\mathcal{M}(\Sigma_D)<0$, we deduce that
\begin{eqnarray*}
\mathcal{M}(\Sigma_D) |\Omega|^{\frac{N+2}{N}} &=& \mathcal{M}(\Sigma_D) \left( |\Omega_1|+|\Omega_2|\right)^{\frac{N+2}{N}}\\
&<& \mathcal{M}(\Sigma_D) \left( |\Omega_1|^{\frac{N+2}{N}} + |\Omega_2|^{\frac{N+2}{N}}\right)\\
&\leq& \E(\Omega_1; \Sigma_D)+ \E(\Omega_2; \Sigma_D)=\E(\Omega; \Sigma_D),
\end{eqnarray*}
which is a contradiction.
\end{proof}

Concerning the regularity of the minimizing set, by the theory of free boundary problems we have:

\begin{proposition}\label{prop:regmin}
Let $\Omega \subset\Sigma_D$ be a minimizer for $\O_1(\Sigma_D)$ and let $\Gamma=\partial\Omega\cap\Sigma_D$ be its relative boundary. Then there exists a critical dimension $d^*$ which can be either 5,6 or 7, such that:
\begin{itemize}
\item[(i)] $\Gamma$ is smooth if $N<d^*$;
\item[(ii)] $\Gamma$ can have countable isolated singularities if $N=d^*$;
\item[(iii)] $\Gamma$ can have a singular set of dimension $N-d^*$, if $N>d^*$.
\end{itemize}
Moreover on the regular part of $\Gamma$ the normal derivative $\frac{\partial u_\Omega}{\partial \nu}$ is constant, namely $\frac{\partial u_\Omega}{\partial \nu}\equiv - \sqrt{\frac{2(N+2)}{N} |\mathcal{O}_1(\Sigma_D)|}$, where $u_\Omega$ is the torsion function of $\Omega$.
\end{proposition}
\begin{proof}
The points (i)--(iii) follows from the results of \cite{DJ, JS} and \cite{W}. Let us prove the last statement. 

Let $\Gamma_{reg}$ be the regular part of $\Gamma$ (which is a relative open set of $\Gamma$), let $x_0\in\Gamma_{reg}$, and let $B_r(x_0)$ be a small ball such that ${B_{2r}(x_0)}\subset \Sigma_D$ and $\Gamma\cap{B_{2r}(x_0)} \subset \Gamma_{reg}$. Moreover, let $\psi \in C^\infty_c(B_r(x_0))$ and consider the vector field $V$ given by $V=\psi \bar\nu$, where $\bar\nu$ is a smooth extension of the normal versor $\left.\nu\right|_{\Gamma\cap{B_{2r}(x_0)}}$ of $\Gamma\cap{B_{2r}(x_0)}$ to a smooth vector field defined in $\overline{B_r(x_0)}$. Hence, by construction, we have that $V:\R^N\to\R^N$ is a smooth vector field with compact support in $B_r(x_0)$, and in particular it holds that $V(0)=0$ and $V(x)=0\in T_x\partial\Sigma_D$ for all $x\in\partial\Sigma_D\setminus\{0\}$. This means that the associated flow $\xi:(-t_0,t_0)\times\overline{\Sigma_D}\to \overline{\Sigma_D}$, for some $t_0>0$, preserves the boundary $\partial\Sigma_D$ and we consider the induced deformation of $\Omega$, $(\Omega_t)_{t\in(-t_0,t_0)}$, where $\Omega_t:=\xi(t,\Omega)$. Actually, since $\mathrm{supp}(\psi) \subset B_r(x_0)$ we infer that $\xi(t,x)=x$ for $x \in B_{\frac{3}{2} r}^\complement(x_0)$, $t\in(-t_0,t_0)$.

Let $u_{\Omega_t}\in H^1_0(\Omega_t;\Sigma_D)$ be the torsion function of $\Omega_t$, for $t\in(-t_0,t_0)$ and let us set $u_t:=u_{\Omega_t}$. Arguing as in the proof of \cite[Proposition 4.3]{PT2} we can prove that the map from $(-t_0,t_0)$ to $H^1(\Sigma_D)$, $t\mapsto u_t$ is differentiable. In particular the function $f:(-t_0,t_0)\to L^1(\Sigma_D)$, given by $f(t)=|\nabla u_t|^2$ is differentiable. We also notice that, since $u_\Omega$ is a weak solution to \eqref{eq:mixbvprobquasiopen}, then by standard elliptic regularity theory $u_\Omega\in W^{2,2}(\Omega\cap B_r(x_0))$. In particular it holds that $V f(0)=\psi\bar\nu|\nabla u_\Omega|^2\in W^{1,1}(\Sigma_D,\R^N)$ and by easy modifications to the proof of \cite[Theorem 5.2.2]{HP} we infer that the function $t\mapsto \mathcal{E}(\Omega_t;\Sigma_D)= -\frac{1}{2}\int_{\Omega_t}|\nabla u_t|^2 \ dx$ is differentiable at $t=0$ and
$$
\frac{d}{dt}\left.\left( \mathcal{E}(\Omega_t;\Sigma_D)\right)\right|_{t=0}= -\int_\Omega \nabla u_\Omega \boldsymbol{\cdot} \nabla u^\prime \ dx - \frac{1}{2}\int_\Omega \mathrm{div}(V |\nabla u_\Omega|^2) \ dx,
$$
where $u^\prime=\frac{d}{dt}\left.\left(u_t\right)\right|_{t=0}$ is a solution to (see the proof of \cite[Proposition 4.3]{PT2})
\begin{equation}\label{eq:mixbvprobuprime2}
\begin{cases}
-\Delta u^\prime = 0 & \text{in $\Omega$,}\\
 u^\prime = -\frac{\partial u_{\Omega_\varphi}}{\partial\nu} \langle V,\nu\rangle & \text{on $\Gamma$,}\\
 \frac{\partial u^\prime}{\partial \nu} = 0 & \text{on $\Gamma_{1}\setminus\{0\}$.}
\end{cases}
\end{equation}
We point out that since the flow $\xi$ leaves invariant $B^\complement_{\frac{3}{2}r}(x_0)$ for all $t\in(-t_0,t_0)$, we have $u^\prime\equiv 0$ in $\Omega\cap B^\complement_{\frac{3}{2}r}(x_0)$ and thus
\beq\label{eq1:constnormderregpart}
\frac{d}{dt}\left.\left( \mathcal{E}(\Omega_t;\Sigma_D)\right)\right|_{t=0}= \underbrace{-\int_{\Omega\cap B_{\frac{3}{2}r}(x_0)} \nabla u_\Omega \boldsymbol{\cdot} \nabla u^\prime \ dx}_{\textbf{(I)}} - \underbrace{\frac{1}{2}\int_\Omega \mathrm{div}(V |\nabla u_\Omega|^2) \ dx}_{\textbf{(II)}}.
\eeq
Let us analize $\textbf{(I)}$. We first observe that as $\Gamma\cap B_{2r}(x_0)\subset \Gamma_{reg}$ and $u^\prime$ is a solution to \eqref{eq:mixbvprobuprime2}, then by standard elliptic regularity theory, it follows that $u^\prime\in W^{2,2}\left(\Omega\cap B_{\frac{3}{2}r}(x_0)\right)$ and it is smooth inside $\Omega$. Hence, applying the Green's formula, taking into account that $\Delta u^\prime=0$ in $\Omega\cap B_{\frac{3}{2}r}(x_0)$, $\frac{\partial u^\prime}{\partial \nu}=0$ on $\Omega\cap\partial B_{\frac{3}{2}r}(x_0)$ (because $u^\prime\equiv 0$ in $B_{\frac{3}{2} r}^\complement(x_0)$ and $u^\prime$ is smooth inside $\Omega$) and $u_\Omega=0$ on $\Gamma\cap \overline{B_{\frac{3}{2}r}(x_0)}$ we get that $\textbf{(I)}=0$. For $\textbf{(II)}$, applying the divergence theorem, and recalling the definition of $V$, at the end we obtain
\begin{equation}\label{eq:diffTorsionalOmegat}
\begin{array}{lll}
\displaystyle \frac{d}{dt}\left.\left( \mathcal{E}(\Omega_t;\Sigma_D)\right)\right|_{t=0} &=&\displaystyle - \frac{1}{2}\int_{\Gamma\cap \overline{B_r(x_0)}}  |\nabla u_\Omega|^2 \psi \ d\sigma.
\end{array}
\end{equation}
On the other hand, for the volume we have
\begin{equation}\label{eq:diffVolOmegat}
 \frac{d}{dt}\left.\left(|\Omega_t|\right)\right|_{t=0} =\int_\Gamma \langle V,\nu\rangle \ d\sigma=\int_{\Gamma\cap \overline{B_r(x_0)}} \psi \ d\sigma.
\end{equation}
Now, since $\Omega$ is a minimizer for $\O_1(\Sigma_D)$, then, recalling Remark \ref{remark1:sect8} and as observed in the proof of Step 2 of Theorem \ref{teo:mainteoexistmin} we get that $\Omega$ is also a minimizer for $\frac{\mathcal{E}(\Omega;\Sigma_D)}{|\Omega|^{\frac{N+2}{N}}}$. Thus from \eqref{eq:diffTorsionalOmegat}, \eqref{eq:diffVolOmegat}, since $|\Omega|=1$ and $\mathcal{E}(\Omega;\Sigma_D)= \mathcal{O}_1(\Sigma_D)<0$, we readily obtain that
$$\frac{d}{dt}\left.\left(\frac{\mathcal{E}(\Omega_t;\Sigma_D)}{|\Omega_t|^{\frac{N+2}{N}}} \right)\right|_{t=0}=- \frac{1}{2}\int_{\Gamma\cap \overline{B_r(x_0)}} \left( |\nabla u_\Omega|^2-\frac{2(N+2)}{N} |\mathcal{O}_1(\Sigma_D)|\right) \psi \ d\sigma=0.$$
Finally, from the arbitrariness of $\psi$ and $x_0$ we conclude that $|\nabla u_\Omega|^2\equiv \frac{2(N+2)}{N} |\mathcal{O}_1(\Sigma_D)|$ on $\Gamma_{reg}$, and as $u_\Omega=0$ on $\Gamma_{reg}$ then by Hopf's lemma it follows that $\frac{\partial u_\Omega}{\partial \nu}\equiv - \sqrt{\frac{2(N+2)}{N} |\mathcal{O}_1(\Sigma_D)|}$ on $\Gamma_{reg}$.

\end{proof}
We conclude this section with the following:
\begin{proof}[Proof of Theorem \ref{mainteo}]
It follows from Theorem \ref{prop:critfirsteigenminim}, Theorem \ref{teo:mainteoexistmin}, Corollary \ref{cor:mainteoexistmin} and Proposition \ref{prop:regmin}.
\end{proof}

\section{The isoperimetric problem and proof of Theorem \ref{mainteo2}}

In this section we study the isoperimetric problem in the class of strictly star-shaped domains in cones, i.e. domains in $\Sigma_D$ whose relative boundary is a radial graph.

Using the same notations of the previous sections, if $\varphi \in C^2(\overline D, \R)$ and $\Omega_\varphi$, $\Gamma_\varphi$ are, respectively, the associated star-shaped domain (see \eqref{def:omegaphi}), and the associated radial graph (see Definition \ref{def:radialgraphassociated}), the (relative) perimeter of $\Omega_\varphi$ in $\Sigma_D$ is given by: 
$$
\P(\Omega_\varphi; \Sigma_D)= \mathcal{H}_{N-1}(\Gamma_\varphi)= \int_D e^{(N-1)\varphi} \sqrt{1+|\nabla \varphi|^2} \ d\sigma,
$$
where $d\sigma$ is the $(N-1)$-dimensional area element of $\S^{N-1}$. Let us observe that $\P$ can be seen as a functional on the space $C^2(\overline D, \R)$ and, contrarily to the torsional energy functional of Section 3, its expression does not involve the associated domain $\Omega_\varphi$. Thus we set for brevity $\P(\varphi)=\P(\Omega_\varphi; \Sigma_D)$.

The derivative of $\mathcal{P}$ along a variation $v\in C^2(\overline D, \R)$ is given by
 \begin{eqnarray*}
\P^\prime(\varphi)[v]&=& \int_D e^{(N-1)\varphi} \left((N-1) \sqrt{1+|\nabla \varphi|^2} v + \frac{\nabla \varphi  \boldsymbol{\cdot} \nabla v}{ \sqrt{1+|\nabla \varphi|^2}}\right) \ d\sigma.
\end{eqnarray*}
Let $\V$ be the volume functional (see \eqref{eq:defVolOmegaphi}). We are concerned with critical points $\varphi$ of $\P$ subject to the volume constraint $\{\V=c\}$. Namely we consider the manifold $M$ defined by \eqref{def:M} and the restriction $\mathcal{I}:=\left.\P\right|_M$. A critical point $\varphi \in M$ for $\mathcal{I}$ satisfies
\beq\label{eq:PprimelambdaVprime}
\P^\prime(\varphi)=\lambda \V^\prime(\varphi),
\eeq
with a Lagrangian multiplier $\lambda \in \R$. In the next two propositions we prove that the radial graph $\Gamma_\varphi$ associated to a critical point $\varphi \in M$ is a CMC hypersurface which intersects orthogonally $\partial\Sigma_D\setminus\{0\}$.

This is well known if the variations of the domains are taken in the whole class of subsets of $\Sigma_D$ of finite relative perimeter (see \cite{RR}) but not obvious in our case.
 
\begin{proposition}\label{prop1sect5}
If $\varphi \in C^2(\overline D,\R)$ is a volume-constrained critical point for $\P$, then the associated radial graph $\Gamma_\varphi$ has constant mean curvature $H\equiv\frac{\lambda}{N-1}$.
\end{proposition}
\begin{proof}
Let $v \in C^{1}_c({D})$ be a variation with compact support. By definition (see \eqref{eq:PprimelambdaVprime}) there exists $\lambda \in \R$ such that
\beq\label{eq1:sect5}
\int_D e^{(N-1)\varphi} \left((N-1) \sqrt{1+|\nabla \varphi|^2} v+ \frac{\nabla \varphi  \boldsymbol{\cdot} \nabla v}{ \sqrt{1+|\nabla \varphi|^2}}\right) \ d\sigma= \lambda \int_D e^{N\varphi}v \ d\sigma.
\eeq
Let us observe that
$$\int_D e^{(N-1)\varphi} \frac{\nabla \varphi  \boldsymbol{\cdot} \nabla v}{ \sqrt{1+|\nabla \varphi|^2}} \ d\sigma= \int_D  \frac{\nabla \varphi  \boldsymbol{\cdot} \nabla (v e^{(N-1)\varphi})}{ \sqrt{1+|\nabla \varphi|^2}} \ d\sigma - \int_D (N-1) e^{(N-1)\varphi} \frac{|\nabla \varphi|^2}{ \sqrt{1+|\nabla \varphi|^2}} v \ d\sigma. $$
Hence we can rewrite \eqref{eq1:sect5} as
$$ \int_D  \frac{\nabla \varphi  \boldsymbol{\cdot} \nabla (v e^{(N-1)\varphi})}{ \sqrt{1+|\nabla \varphi|^2}} \ d\sigma + \int_D e^{(N-1)\varphi} \frac{N-1}{\sqrt{1+|\nabla \varphi|^2}} v \ d\sigma = \lambda  \int_D e^{N\varphi}v \ d\sigma. $$
Now, since $\varphi$ is smooth and $v$ has compact support in $D$, integrating by parts, we get
\beq\label{eq1:compactsupp}
\int_D  \frac{\nabla \varphi  \boldsymbol{\cdot} \nabla (v e^{N\varphi})}{ \sqrt{1+|\nabla \varphi|^2}} \ d\sigma= - \int_D  \mathrm{div}_{\S^{N-1}}\left(\frac{\nabla \varphi}{ \sqrt{1+|\nabla \varphi|^2}}\right) e^{N\varphi} v \ d\sigma, 
\eeq
and thus we deduce
$$ \int_D  \left(-\mathrm{div}_{\S^{N-1}}\left(\frac{\nabla \varphi}{ \sqrt{1+|\nabla \varphi|^2}}\right) e^{(N-1)\varphi} +  \frac{N}{\sqrt{1+|\nabla \varphi|^2}} e^{(N-1)\varphi}\right) v \ d\sigma =   \int_D \lambda e^{N\varphi} v \ d\sigma.$$
Therefore, as $v$ is arbitrary, we obtain
$$ -\mathrm{div}_{\S^{N-1}}\left(\frac{\nabla \varphi}{ \sqrt{1+|\nabla \varphi|^2}}\right) e^{(N-1)\varphi} +  \frac{N-1}{\sqrt{1+|\nabla \varphi|^2}} e^{(N-1)\varphi} = \lambda e^{N\varphi}\ \ \ \hbox{in $D$},$$
i.e.
\beq\label{eqradgraphlambda}
 -\mathrm{div}_{\S^{N-1}}\left(\frac{\nabla \varphi}{ \sqrt{1+|\nabla \varphi|^2}}\right) +  \frac{N-1}{\sqrt{1+|\nabla \varphi|^2}} = \lambda e^{\varphi} \ \ \hbox{in $D$}.
\eeq
Comparing \eqref{eqradgraphlambda} with \eqref{eq:divformt} it follows that the mean curvature of $\Gamma_\varphi$ is constant and it is equal to $\frac{\lambda}{N-1}$.
\end{proof}
\begin{proposition}
If $\varphi$ is as in the statement of Proposition \ref{prop1sect5} then $\Gamma_\varphi$ intersects orthogonally $\partial\Sigma_D\setminus\{0\}$.
\end{proposition}

\begin{proof}
Let $\nu_{_{\partial\Sigma_D}}$ be the exterior unit normal to $\partial \Sigma_D\setminus\{0\}$, and let $\nu_{_{\Gamma_\varphi}}$ be the exterior unit normal to $\Gamma_\varphi$. By \eqref{eq:Gauss} we have
\beq\label{eq:expressnormalgammaphi}
\nu_{_{\Gamma_\varphi}}(\mathcal{Y}(q))=  \frac{q- \nabla \varphi }{(1+|\nabla \varphi|^2)^{1/2}},
\eeq
where $\mathcal{Y}$ is the standard parametrization of $\Gamma_\varphi$ defined by \eqref{eq:defparstand}. Notice that, since by assumption $\varphi$ is smooth up to the boundary, then \eqref{eq:expressnormalgammaphi} is well defined on $\overline D$. If $p \in (\partial \Sigma_D\setminus\{0\}) \cap \overline{\Gamma_\varphi}$, then by definition the intersection is orthogonal at $p$ if and only if $\nu_{_{\partial\Sigma_D}}(p)  \boldsymbol{\cdot} \nu_{_{\Gamma_\varphi}}(p) = 0$. Therefore, writing $p=\mathcal{Y}(q)$ this is equivalent to
$$\nu_{_{\partial\Sigma_D}}(\mathcal{Y}(q))  \boldsymbol{\cdot} (q- \nabla \varphi(q)) = 0.$$
Since $p=\mathcal{Y}(q) \in \partial \Sigma_D\setminus\{0\}$ and $\partial \Sigma_D\setminus\{0\}$ is the boundary of a cone we have $\nu_{_{\partial\Sigma_D}}(\mathcal{Y}(q)) \boldsymbol{\cdot} q = 0$ and thus the intersection between $\partial \Sigma_D\setminus\{0\}$ and $\Gamma_\varphi$ is orthogonal if and only if
\beq\label{eq:previouseq}
\nu_{_{\partial\Sigma_D}}(\mathcal{Y}(q))  \boldsymbol{\cdot} \nabla \varphi(q)= 0  \ \ \ \forall q \in \partial D.
\eeq
Exploiting again that $\partial\Sigma_D$ is a cone, we have $\nu_{_{\partial\Sigma_D}}(p)=\nu_{_{\partial\Sigma_D}}(t p)$ for any $p\in \partial \Sigma_D\setminus\{0\}$, $t>0$. Hence, since $\mathcal{Y}(q) \in \partial \Sigma_D\setminus\{0\}$, we have $\nu_{_{\partial\Sigma_D}}(\mathcal{Y}(q))=\nu_{_{\partial\Sigma_D}}(q)=\nu_{_{\partial D}}(q)$ for any $q \in \partial D$, where $\nu_{_{\partial D}}$ is the exterior unit co-normal to $\partial D$, and thus \eqref{eq:previouseq} is equivalent to
\beq\label{eq:thesiscorol}
\frac{\partial \varphi}{\partial \nu_{_{\partial D}}} = 0 \ \ \hbox{on $\partial D$}.
\eeq
To prove this, we argue as in the proof of Proposition \ref{prop1sect5}. Taking a variation $v \in C^1(\overline D, \R)$ and integrating by parts we have
$$
\int_D  \frac{\nabla \varphi  \boldsymbol{\cdot} \nabla (v e^{(N-1)\varphi})}{ \sqrt{1+|\nabla \varphi|^2}} \ d\sigma=\int_{\partial D} e^{(N-1)\varphi} v \ \Big\langle \frac{\nabla \varphi}{ \sqrt{1+|\nabla \varphi|^2}},  \nu_{_{\partial D}} \Big{\rangle} \ d\hat\sigma - \int_D  \mathrm{div}_{\S^{N-1}}\left(\frac{\nabla \varphi}{ \sqrt{1+|\nabla \varphi|^2}}\right) e^{(N-1)\varphi} v \ d\sigma.
$$
Using this and arguing as in the proof Proposition \ref{prop1sect5}, since $\varphi$ satisfies the equation \eqref{eqradgraphlambda} we obtain 
$$\int_{\partial D} e^{(N-1)\varphi} v \ \Big\langle \frac{\nabla \varphi}{ \sqrt{1+|\nabla \varphi|^2}},  \nu_{_{\partial D}} \Big{\rangle} \ d\hat\sigma=0.$$
Since $v \in C^1(\overline D, \R)$ is arbitrary we can choose $v$ such that $v= \ \Big\langle \frac{\nabla \varphi}{ \sqrt{1+|\nabla \varphi|^2}},  \nu_{_{\partial D}} \Big{\rangle}$ on $\partial D$ and thus
$$\int_{\partial D} e^{(N-1)\varphi}   \left|\Big\langle \frac{\nabla \varphi}{ \sqrt{1+|\nabla \varphi|^2}},  \nu_{_{\partial D}} \Big{\rangle}\right|^2 \ d\hat\sigma=0,$$
which gives $\Big\langle \frac{\nabla \varphi}{ \sqrt{1+|\nabla \varphi|^2}},  \nu_{_{\partial D}} \Big{\rangle}\equiv 0$ on $\partial D$, and thus \eqref{eq:thesiscorol} is proved.
\end{proof}

Analogously to Lemma \ref{lem:secvarI}, if $\varphi \in M$ is a critical point for $\mathcal{I}$ then
$$\mathcal{I}^{\prime\prime}(\varphi)=\P^{\prime\prime}(\varphi)-\lambda \V^{\prime\prime}(\varphi).$$
Choosing $c=|\Omega_0|=|\Sigma_D\cap B_1(0)|$ in \eqref{def:M} we observe that the function $\varphi\equiv 0$ belongs to $M$ and it is a critical point for $\mathcal{I}$. In particular \eqref{eq:PprimelambdaVprime} yields $\lambda=N-1$. Moreover, for any $v,w \in T_0M$, since
 \begin{eqnarray*}
\P^{\prime\prime}(0)[v,w]&=& \int_D \left((N-1)^2 vw + {\nabla \varphi  \boldsymbol{\cdot} \nabla v}\right) \ d\sigma,
\end{eqnarray*}
and recalling \eqref{eq:secvarvol}, it follows that
 \begin{equation}\label{exprIvw0}
\mathcal{I}^{\prime\prime}(0)[v,w]=\P^{\prime\prime}(0)[v,w]-(N-1) \V^{\prime\prime}[v,w]= \int_D  \left(\nabla v \boldsymbol{\cdot} \nabla w - (N-1) vw\right) \ d\sigma.
\end{equation}
From \eqref{exprIvw0} we easily have the analogue of Theorem \ref{prop:critfirsteigenminim} for the perimeter functional $\mathcal{I}$.

\begin{theorem}\label{prop:nonisopercrit}
Let $\lambda_1(D)$ be the first nontrivial eigenvalue of $-\Delta_{\S^{N-1}}$ on the domain $D$ with zero Neumann condition on $\partial D$. Then:
\begin{itemize}
\item[(i)] if $\lambda_1(D) < N-1$ then $\varphi\equiv 0$ is not a local minimizer for $\mathcal{I}$;
\item[(ii)] if $\lambda_1(D) > N-1$ then $\varphi\equiv 0$ is a local minimizer for $\mathcal{I}$.
\end{itemize}
\end{theorem}
\begin{proof}
Since $T_0M$ is made by functions with zero mean value (see \eqref{def:TvarphiM}), considering the $L^2$-normalized eigenfunction $w_1$ corresponding to the eigenvalue $\lambda_1(D)$, from \eqref{exprIvw0} we get $\mathcal{I}^{\prime\prime}(0)[w_1,w_1]<0$ whenever $\lambda_1(D)<N-1$. This proves $(i)$.

Viceversa, if $\lambda_1(D)>N-1$, from \eqref{exprIvw0} and the variational characterization of $\lambda_1(D)$ we get that $\mathcal{I}^{\prime\prime}[v,v]>0$ for all $v\in T_0M$ with $v\neq 0$, and hence $(ii)$ holds.
\end{proof}
To find examples of domains $D\subset\S^{N-1}$ satisfying $\lambda_1(D)<N-1$ we can use the function $u_e\in C^\infty(\S^{N-1})$ introduced in \eqref{eq:defue} and Proposition \ref{lemma:unstablcrit}. Hence for the nonconvex domains constructed in the Appendix, the spherical sectors are not the minimizers of $\mathcal{I}$.\\

Concerning the existence of a minimizer for the relative perimeter $\P(E;\Sigma_D)$ in the whole class of finite perimeter subsets of $\Sigma_D$, with a fixed volume, we summarise in the following the results stated in \cite{RR}.
\begin{theorem}\label{teo:RR}
Let $D\subset \S^{N-1}$ be a domain such that $\mathcal{H}_{N-1}(D)\leq \mathcal{H}_{N-1}(\S^{N-1}_+)$. Then, there exists a set of finite perimeter $E^*$ inside $\Sigma_D$ which minimizes the relative perimeter under a volume constraint, for any value of the volume. Moreover any minimizer of the relative perimeter, with fixed volume, is a bounded set.
 \end{theorem}
\begin{proof}
It follows from Proposition 3.5 and Proposition 3.7 in \cite{RR}.
\end{proof}
We conclude this section with the following
\begin{proof}[Proof of Theorem \ref{mainteo2}]
The existence of a set of finite perimeter $E^*$ inside $\Sigma_D$ which minimizes the relative perimeter under a volume constraint, and its boundedness, follows from Theorem \ref{teo:RR}. From Theorem \ref{prop:nonisopercrit} we infer that $E^*$ cannot be a spherical sector, while the properties (i)-(iii) of $\Gamma_{E^*}$ derive from classical results for isoperimetric problems (see e.g. \cite[Sect. 2]{RR} and the references therein).
\end{proof}
\appendix

\section{Examples of non convex domains satisfying condition \eqref{eq:hpmainteo}}
In this section we construct two classes of non-convex domains of $\S^{N-1}$ satisfying hypothesis \eqref{eq:hpmainteo} of Theorem \ref{mainteo}. In particular to show the instability condition $\lambda_1(D)<N-1$ we will prove the inequality (i) of Proposition \ref{lemma:unstablcrit}. We begin with some preliminary technical results.\\

Let $N\geq 3$, let $\{{e}_1,\ldots, {e}_N\}$ be the standard basis in $\R^N$, fix $\theta\in(-\frac{\pi}{2},\frac{\pi}{2})$ and consider the vector ${e}_{\theta}\in \S^{N-1}$ defined by $${e}_{\theta}:={e}_{1} \sin\theta + {e}_{N} \cos\theta.$$
By construction the vector ${e}_{\theta}$ lies in the sector $\{ x_1>0, \ x_N>0\}$ if $\theta\in(0,\frac{\pi}{2})$, in $\{x_1<0, \ x_N>0\}$ if $\theta\in(-\frac{\pi}{2},0)$, and coincides with $ {e}_{N}$ if $\theta=0$. For any given $r\in(0,1)$ we consider the hyperplane $H_{\theta,r}$ orthogonal to ${e}_{\theta}$ and passing through $r {e}_{\theta}$, i.e. $$H_{\theta,r}=\{x\in\RN; \ x \boldsymbol{\cdot} {e}_{\theta}= r\}.$$ Let $D_{\theta,r}$ be the region of $\S^{N-1}$ above $H_{\theta,r}$, namely $D_{\theta,r}$ is the spherical cap given by 
\beq\label{eq:defDthetar}
D_{\theta,r}=\{x\in \S^{N-1}; \ x \boldsymbol{\cdot} {e}_{\theta}> r \}.
\eeq
By definition it is easy to check that $D_{\theta,r}$ has angular radius $\arctan\left(\frac{\sqrt{1-r^2}}{r}\right)$.
We consider the function $u_{{e_1}}$ defined in \eqref{eq:defue}, namely $u_{{e_1}}(x)=x \boldsymbol{\cdot} {e}_1=x_1$.
The first result we prove is an explicit formula for the boundary integral appearing in Proposition \ref{lemma:unstablcrit} applied to the function $u_{{e_1}}$ and to the domain $D_{\theta,r}$.

\begin{lemma}\label{lem:signintegralboundary}
For $\theta\in (-\frac{\pi}{2},\frac{\pi}{2})$ and $r\in(0,1)$  it holds
\beq\label{eq:tesiintegralsign}
\int_{\partial D_{\theta,r}} u_{{e}_1} \frac{\partial u_{{e}_1}}{\partial \nu} \ d\hat\sigma=r(1-r^2)^{\frac{N-1}{2}}c_N(1-N\sin^2\theta),
\eeq
where $\nu$ is the exterior unit co-normal (i.e., for any $x\in \partial D_{\theta,r}$, $\nu(x)$ is the unique unit vector in $T_x\S^{N-1}$ which is orthogonal to $T_x\partial D_{\theta,r}$ and pointing outward $D_{\theta,r}$), $d\hat\sigma$ is the $(N-2)$-dimensional area element of ${\partial D_{\theta,r}}$, $c_N$ is a positive constant depending only on $N$ which is explicit (see \eqref{eq:CN}). 
\end{lemma}

\begin{proof}
Let $\theta\in(-\frac{\pi}{2},\frac{\pi}{2})$, $r\in(0,1)$. We begin observing that for any $x\in \partial D_{\theta,r}$ it holds 
\beq\label{eq:normext}
\nu(x)=\frac{1}{\sqrt{1-r^2}} (rx-{e}_\theta).
\eeq
Indeed, $x \boldsymbol{\cdot} {e}_\theta=r$ and thus $(rx-{e}_\theta) \boldsymbol{\cdot} (rx-{e}_\theta)=1-r^2$, namely $\left|\frac{1}{\sqrt{1-r^2}} (rx-{e}_\theta)\right|=1$, and we readily check that $\frac{1}{\sqrt{1-r^2}} (rx-{e}_\theta) \boldsymbol{\cdot} x= \frac{1}{\sqrt{1-r^2}} (r-x \boldsymbol{\cdot} {e}_\theta)=0$, which means that $\frac{1}{\sqrt{1-r^2}} (rx-{e}_\theta)\in T_x\S^{N-1}$. In addition, if $v\in T_x\partial D_{\theta,r}$ and $\gamma:(-\delta,\delta)\to \partial D_{\theta,r}$ is curve such that $\gamma(0)=x$ and $\gamma^\prime(0)=v$, for some small $\delta>0$, then as $ \partial D_{\theta,r}\subset\S^{N-1}$, differentiating the identity $|\gamma|^2\equiv1$ we get $x \boldsymbol{\cdot} v=0$. Similarly, differentiating the identity $\gamma \boldsymbol{\cdot} {e}_\theta\equiv r$ we obtain $v \boldsymbol{\cdot}  {e}_\theta=0$. Hence, $\frac{1}{\sqrt{1-r^2}} (rx-{e}_\theta) \boldsymbol{\cdot} v=0$ and thus from the arbitrariness of $v\in T_x\partial D_{\theta,r}$ we infer that $\frac{1}{\sqrt{1-r^2}} (rx-{e}_\theta) \perp T_x\partial D_{\theta,r}$, so that $\nu(x)=\pm \frac{1}{\sqrt{1-r^2}} (rx-{e}_\theta)$. To choose the right sign we now show that $\frac{1}{\sqrt{1-r^2}} (rx-{e}_\theta)$ points outward $D_{\theta,r}$. To this end, let $x\in \partial D_{\theta,r}$ and consider the vector $y_s:=x+ \frac{s}{\sqrt{1-r^2}} (rx-{e}_\theta)$, with $|s|$ small. As $r\in(0,1)$ and since $y_s \boldsymbol{\cdot} {e}_\theta = r + \frac{s(r^2-1)}{\sqrt{1-r^2}}$, then, using the definition \eqref{eq:defDthetar} and by elementary computations we readily check that $\frac{y_s}{|y_s|} \in D_{\theta,r}$ if and only if $s<0$ with $|s|$ is small enough, and this proves that $\frac{1}{\sqrt{1-r^2}} (rx-{e}_\theta)$ points outward $D_{\theta,r}$.\\

With \eqref{eq:normext} at hand, and recalling the definition of $u_{{e}_1}$ we easily obtain that
\beq\label{eq:integrandfunct}
u_{{e}_1}(x) \frac{\partial u_{{e}_1}}{\partial \nu}(x) =  \frac{1}{\sqrt{1-r^2}} (r x_1^2-\sin\theta x_1),
\eeq
for any $x=(x_1,\ldots,x_N)\in \partial D_{\theta,r}$. In order to compute the integral in \eqref{eq:tesiintegralsign} we determine a suitable parametrization for $\partial D_{\theta,r}$. To this aim, we recall that $x=(x_1,\ldots,x_N)\in \partial D_{\theta,r}$ if and only if $x$ satisfies
$$\begin{cases} x_1^2+\ldots+x_N^2=1,\\
x_1\sin\theta+x_N\cos\theta=r,\end{cases}$$
and by elementary computations we find that
\beq\label{eq2:parampartialD}
 \frac{(x_1-r\sin\theta)^2}{\cos^2\theta}+x_2^2\ldots+x_{N-1}^2=1-r^2, \ \ x_N=\frac{r}{\cos\theta} - x_1\frac{\sin\theta}{\cos\theta}.
\eeq
Using the spherical coordinates for $\S^{N-2}\subset \R^{N-1}$ and the second equation in \eqref{eq2:parampartialD} we find a parametrization for $\partial D_{\theta,r}$. Indeed, first assuming that $N\geq 4$ and denoting by $\phi_1,\ldots,\phi_{N-2}$ the angular coordinates, where $\phi_1,\ldots,\phi_{N-3}$ have the range $(0,\pi)$  and $\phi_{N-2}$ ranges over $(0,2\pi)$, we easily check that

\begin{equation}
\left[\begin{array}{c}
x_{1}\\
x_{2}\\
\vdots\\
x_{N-2}\\
x_{N-1}\end{array}\right]
=\sqrt{1-r^2}\left[\begin{array}{c}
\cos\theta \cos \phi_{1} \\
\sin \phi_1 \cos\phi_{2} \\
\vdots\\
\sin\phi_1\ldots\sin\phi_{N-3} \cos\phi_{N-2}\\
\sin\phi_1\ldots\sin\phi_{N-3} \sin\phi_{N-2}
\end{array}\right] + \left[\begin{array}{c}
r\sin\theta\\
0\\
\vdots\\
0\\
0
\end{array}\right]
\end{equation}
is a parametrization of the ellipsoid in $\R^{N-1}$ (or the sphere if $\theta=0$) described by the first equation in \eqref{eq2:parampartialD} (up to a zero measure set with respect to the $(N-2)$-dimensional Haussdorf measure). Hence, exploiting the second equation in \eqref{eq2:parampartialD}, we deduce that $\psi:(0,\pi)\times\ldots\times(0,\pi)\times(0,2\pi)\to \R^N$ defined by 
\beq\label{eq:parampsi}
\psi(\phi_1,\ldots,\phi_{N-1}):=\sqrt{1-r^2}\left[\begin{array}{c}
\cos\theta \cos\phi_{1} \\
\sin\phi_1\cos\phi_{2} \\
\vdots\\
\sin\phi_1\ldots \cos\phi_{N-2}\\
\sin\phi_1\ldots \sin\phi_{N-2}\\
- \sin\theta \cos\phi_{1} \ \ \ 
\end{array}\right] + \left[\begin{array}{c}
r\sin\theta\\
0\\
\vdots\\
0\\
0\\
r\cos\theta
\end{array}\right]
\eeq
is a parametrization of the $(N-2)$-dimensional sphere $\partial D_{\theta,r}$. Then, arguing by induction and after a straightforward computation we see that the coefficients $g_{ij}=\left\langle\frac{\partial \psi}{\partial \phi_i}, \frac{\partial \psi}{\partial \phi_j}\right\rangle$, $i,j=1,\ldots,N-2$, of the induced metric on $\partial D_{\theta,r}$, are given by $$g_{ij}=\begin{cases} 1-r^2 & \hbox{if $i=j=1$,}\\
(1-r^2)\sin^2\phi_1 &  \hbox{if $i=j=2$,}\\
 (1-r^2)\sin^2\phi_1\cdots\sin^2\phi_{i-1} &  \hbox{if $i=j$ and $i\neq 1,2$,}\\ 0 &\hbox{if $i\neq j$}.\end{cases}$$
In particular, the matrix $(g_{ij})_{i,j=1,\ldots,N-2}$ is diagonal, positive definite and the square root of its determinant is given by 
$$\sqrt{Det(g_{ij})}=\begin{cases} (1-r^2) \sin\phi_1 & \hbox{if $N=4$}, \\
 (1-r^2)^{\frac{3}{2}}(\sin\phi_1)^{2}(\sin\phi_2) &  \hbox{if $N=5$,}\\
(1-r^2)^{\frac{N-2}{2}}(\sin\phi_1)^{N-3}(\sin\phi_2)^{N-4}\cdots\sin\phi_{N-3} &  \hbox{if $N\geq 6$.}\end{cases}$$
Therefore, when $N\geq 6$ (the other cases $N=4,5$ being similar and easier) the $(N-2)$-dimensional area element of $\partial D_{\theta,r}$ is expressed in these local coordinates by $$d\hat\sigma=(1-r^2)^{\frac{N-2}{2}}(\sin\phi_1)^{N-3}(\sin\phi_2)^{N-4}\cdots\sin\phi_{N-3}\ d\phi_1\cdots d\phi_{N-3} d\phi_{N-2}.$$
Setting for brevity $G(\phi_2,\ldots,\phi_{N-2}):=(\sin\phi_2)^{N-4}\cdots\sin\phi_{N-3}$, $d\Phi:=d\phi_{2}\cdots d\phi_{N-3} d\phi_{N-2}$, and observing that
\beq\label{eq:intGPhi}
\int_{(0,\pi)^{N-4}\times(0,2\pi)} G(\phi_2,\ldots,\phi_{N-2})\ d\Phi=|\S^{N-3}|=(N-2)\omega_{N-2},
\eeq
where $\omega_{N-2}$ is the volume of the unit ball in $\R^{N-2}$, then, recalling \eqref{eq:integrandfunct}, exploiting Fubini's theorem and taking into account \eqref{eq:parampsi}, \eqref{eq:intGPhi} we get that
\beq\label{eq:boundaryintegralcomp1}
\begin{array}{lll}
&&\displaystyle \int_{\partial D_{\theta,r}} u_{{e}_1} \frac{\partial u_{{e}_1}}{\partial \nu} d\hat\sigma\\[12pt]
&=&\displaystyle (1-r^2)^{\frac{N-3}{2}}\left[ \int_{(0,\pi)^{N-3}\times(0,2\pi)}    r\left(\sqrt{1-r^2} \cos\theta\cos\phi_1+r\sin\theta\right)^2 (\sin\phi_1)^{N-3}G\ d\phi_{1} d\Phi\right.\\[12pt]
&&\displaystyle\left.\quad\quad\ \ \quad\quad\quad-\int_{(0,\pi)^{N-3}\times(0,2\pi)}    \sin\theta\left(\sqrt{1-r^2} \cos\theta\cos\phi_1+r\sin\theta\right) (\sin\phi_1)^{N-3}G\ d\phi_{1} d\Phi\right]\\[12pt]
&=&\displaystyle (N-2)\omega_{N-2}(1-r^2)^{\frac{N-3}{2}} \int_0^{\pi} \left[ r(1-r^2) \cos^2\theta\cos^2\phi_1+r^3\sin^2\theta-r\sin^2\theta\right](\sin\phi_1)^{N-3} d\phi_{1},
\end{array}
\eeq
where in the last integral we have discarded the linear terms in $\cos\phi_1$ because
$$\int_0^\pi \cos\phi_1 (\sin\phi_1)^{N-3}\ d\phi_{1}=\frac{1}{N-2}(\sin \phi_1)^{N-2}\big|_0^{\pi}=0. $$
Rearranging the terms in the last integral of \eqref{eq:boundaryintegralcomp1} we obtain
\beq\label{eq:boundaryintegralcomp2}
\begin{array}{lll}
\displaystyle \int_{\partial D_{\theta,r}} u_{{e}_1} \frac{\partial u_{{e}_1}}{\partial \nu} d\hat\sigma
&=&\displaystyle (N-2)\omega_{N-2}r(1-r^2)^{\frac{N-1}{2}} \int_0^{\pi} \left(\cos^2\theta\cos^2\phi_1-\sin^2\theta\right)(\sin\phi_1)^{N-3} d\phi_{1}.
\end{array}
\eeq
Now, let us observe that for any $k\in \N$, $k\geq 2$, integrating by parts we have
\beq\label{eq:iterationsink}
\begin{array}{lll}
\displaystyle \int_0^\pi  \cos^2\phi_1(\sin\phi_1)^{k} \ d\phi_1&=&\displaystyle  \int_0^\pi  (\sin\phi_1)^{k} - (\sin\phi_1)^{k+2} \ d\phi_1\\[12pt]
&=&\displaystyle \int_0^\pi  (\sin\phi_1)^{k}  \ d\phi_1- \frac{k+1}{k+2} \int_0^\pi  (\sin\phi_1)^{k}  \ d\phi_1\\[8pt]
&=&\displaystyle  \frac{1}{k+2} \int_0^\pi (\sin\phi_1)^{k} \ d\phi_1.
\end{array}
\eeq
Combining \eqref{eq:boundaryintegralcomp2} and \eqref{eq:iterationsink} (with $k=N-3$) we get that
\beq\label{eq:boundaryintegralcomp3}
\begin{array}{lll}
\displaystyle \int_{\partial D_{\theta,r}} u_{{e}_1} \frac{\partial u_{{e}_1}}{\partial \nu} d\hat\sigma
&=&\displaystyle (N-2)\omega_{N-2}r(1-r^2)^{\frac{N-1}{2}} \int_0^\pi (\sin\phi_1)^{N-3} \ d\phi_1 \left(\frac{1}{N-1}\cos^2\theta-\sin^2\theta\right)\\[12pt]
&=&\displaystyle (N-2)\omega_{N-2}r(1-r^2)^{\frac{N-1}{2}} \int_0^\pi (\sin\phi_1)^{N-3} \ d\phi_1 \frac{\left(1-N\sin^2\theta\right)}{N-1},
\end{array}
\eeq
and this proves \eqref{eq:tesiintegralsign}, with $c_N$ given by
\beq\label{eq:CN}
c_N:= \frac{N-2}{N-1}\omega_{N-2} \int_0^\pi (\sin\phi_1)^{N-3} \ d\phi_1.
\eeq
The proof of the Lemma is then complete for $N\geq 4$.\\

When $N=3$ the proof is much simpler. Indeed, in this case, a parametrization of the circle $\partial D_{\theta,r}$ is given by the map $\psi:(0,2\pi)\to \R^3$ defined by
$$
\psi(\phi):=\sqrt{1-r^2}\left[\begin{array}{c}
\cos\theta \cos\phi \\
\sin\phi \\
- \sin\theta \cos\phi \ \ \ 
\end{array}\right] + \left[\begin{array}{c}
r\sin\theta\\
0\\
r\cos\theta
\end{array}\right].
$$
Hence, using the definition, \eqref{eq:integrandfunct}, and taking into account that $|\psi^\prime(\phi)|=\sqrt{1-r^2}$, we easily check that
\begin{eqnarray*}
 &&\int_{\partial D_{\theta,r}} u_{{e}_1} \frac{\partial u_{{e}_1}}{\partial \nu} d\hat\sigma\\[10pt]
 &=&\frac{1}{\sqrt{1-r^2}}\int_0^{2\pi}\left[r\left(\sqrt{1-r^2}\cos\theta \cos\phi +r\sin\theta\right)^2 - \sin\theta \left(\sqrt{1-r^2}\cos\theta \cos\phi+r\sin\theta\right)\right]|\psi^\prime(\phi)| \ d\phi\\[10pt]
 &=&r(1-r^2)\int_0^{2\pi}\cos^2\theta \cos^2\phi  - \sin^2\theta \ d\phi=\pi r(1-r^2)\left(\cos^2\theta   - 2\sin^2\theta \right)=\pi r(1-r^2)\left(1   - 3\sin^2\theta \right).
\end{eqnarray*}
In particular, as $\omega_1=2$, \eqref{eq:boundaryintegralcomp3} holds true with $c_N$ given by \eqref{eq:CN} even for $N=3$. The proof is complete.
\end{proof}
\begin{remark}
Let $\theta\in(0,\frac{\pi}{2})$. From the analytic expression of $\partial D_{\theta,r}$ given by \eqref{eq2:parampartialD} or, equivalently, by using \eqref{eq:parampsi}, but with $\phi_1$ varying in $[0,2\pi)$, it is easy to check that $\partial D_{\theta,r}$ is contained in the sector $\{x\in \S^{N-1}; \ x_1>0, \ x_N>0\}$ if $-\sqrt{1-r^2} \cos\theta + r \sin\theta>0$ and $-\sqrt{1-r^2} \sin\theta + r \cos\theta>0$, which both hold true if $r>\sqrt{\max\{\cos^2\theta, \sin^2\theta\}}$. By a similar argument, we take $-\theta\in (-\frac{\pi}{2},0)$, and under the same condition on $r$, then $D_{-\theta,r}$ is contained in$\{x\in \S^{N-1}; \ x_1<0, \ x_N>0\}$.
\end{remark}

An immediate consequence of the previous remark and Lemma \ref{lem:signintegralboundary} is the following.
\begin{corollary}\label{cor:boundintegsign}
Let $\theta\in(\arcsin(\frac{1}{\sqrt{N}}), \frac{\pi}{2})$ and set $r_\theta:=\sqrt{\max\{\cos^2\theta, \sin^2\theta\}}$. Then, for any $r\in(r_\theta,1)$, the spherical cap $D_{\theta,r}$ is contained in $\{x\in \S^{N-1}; \ x_1>0, \ x_N>0\}$, while the symmetrical domain (with respect to the hyperplane $\{x_1=0\}$), namely $D_{-\theta,r}$, is contained in $\{x\in \S^{N-1}; \ x_1<0, \ x_N>0\}$. Moreover it holds that
\beq\label{eq1:corboundaryintegral}
\int_{\partial D_{\theta,r}} u_{{e}_1} \frac{\partial u_{{e}_1}}{\partial \nu} \ d\hat\sigma + \int_{\partial D_{-\theta,r}} u_{{e}_1} \frac{\partial u_{{e}_1}}{\partial \nu} \ d\hat\sigma<0
\eeq
and
\beq\label{eq2:corboundaryintegral}
\int_{D_{\theta,r}\cup D_{-\theta,r}} u_{{e}_1} \ d\sigma= 0.
\eeq
\end{corollary}
Notice that \eqref{eq2:corboundaryintegral} follows from the symmetry of $D_{\theta,r} \cup D_{-\theta,r}$, as $u_{{e}_1}$ is odd. Since $D_{\theta,r}\cup D_{-\theta,r}$ is not connected, in order to apply the instability criterion given by Proposition \ref{lemma:unstablcrit} our idea is to join the two domains $D_{\theta,r}$, $D_{-\theta,r}$ by a suitably ``small'' tunnel-like domain which is symmetric with respect to the hyperplane $\{x_1=0\}$. More precisely, we have the following.
\begin{example}
\label{example1}
 Let $\e>0$ be a small number to be determined later and consider the open region $A_\e\subset \S^{N-1}$ between the two symmetric hyperplanes $\{x_{N-1}=-\e\}$ and $\{x_{N-1}=+\e\}$, namely
$$A_{\e}:=\{x\in \S^{N-1}; \ -\e< x \boldsymbol{\cdot}{e_{N-1}} < \e\}.$$
Setting $\partial A_\e^+:=\{x\in \S^{N-1}; \  x \boldsymbol{\cdot}{e_{N-1}} =\e\}$, $\partial A_\e^-:=\{x\in \S^{N-1}; \  x \boldsymbol{\cdot}(-{e_{N-1}}) =\e\}$ and arguing as in the proof of Lemma \ref{lem:signintegralboundary} (see \eqref{eq:integrandfunct}) we can check that
the exterior unit co-normal to $\partial A_{\e}=\partial A_{\e}^-\cup A_\e^+$, pointing outwards $A_{\e}$, is given by 
\beq\label{eq:normalAeps}
\nu(x)=\begin{cases} \frac{1}{\sqrt{1-\e^2}} ({e}_{N-1} - \e x) & \hbox{if $x\in \partial A_\e^+$},\\
\frac{1}{\sqrt{1-\e^2}} (-{e}_{N-1}-  \e x) & \hbox{if $x\in \partial A_\e^-$}.
 \end{cases} 
\eeq
In view of \eqref{eq:normalAeps}, and as ${e}_1 \boldsymbol{\cdot} {e}_{N-1}=0$, it follows that $\nu(x) \boldsymbol{\cdot} {e}_1= -  \e x_1$ for all $x\in  \partial A_{\e}=\partial A_{\e}^-\cup A_\e^+$, and thus
\beq\label{eq:integrandpartialAeps}
u_{{e}_1}(x) \frac{\partial u_{{e}_1}}{\partial \nu}(x)=-\e x_1^2\ \ \ \hbox{for all $x\in \partial A_{\e}$.}
\eeq

Now, fixing $\theta$ and $r\in(r_\theta,1)$ as in Corollary \ref{cor:boundintegsign} we can choose $\eps>0$ sufficiently small (depending on $r$ and $\theta$) so that $\partial A_{\e}$ intersects $D_{\theta,r}\cup D_{-\theta,r}$. We then take as tunnel-like domain the connected subset of
 $$A_{\e}\setminus\left({D_{\theta,r}}\cup {D_{-\theta,r}}\right)=\{x\in \S^{N-1}; \ -\e< x \boldsymbol{\cdot}{e_{N-1}} < \e, \ x \boldsymbol{\cdot}{e_\theta}\leq r, \ x \boldsymbol{\cdot}{e_{-\theta}}\leq r\}$$
containing ${e}_N$, and we denote it by $T_{\e,\theta,r}$. We set $$D_{\e,\theta,r}:= D_{\theta,r}\cup T_{\e,\theta,r}\cup D_{-\theta,r}.$$
By definition it is easy to check that $D_{\e,\theta,r}$ is a domain symmetric with respect to the hyperplane $\{x_1=0\}$ and thus, as $u_{{e}_1}$ is odd, we have
$$ \int_{D_{\eps,\theta,r}} u_{{e}_1}\ d\sigma= 0.$$
Then by \eqref{eq1:corboundaryintegral} and \eqref{eq:integrandpartialAeps} we easily check that
$$
\int_{\partial D_{\e,\theta,r}} u_{{e}_1} \frac{\partial u_{{e}_1}}{\partial \nu} \ d\hat\sigma < 0
$$
if $\e>0$ is sufficiently small, so that $D_{\eps, \theta,r}$ satisfies (i) of Proposition \ref{lemma:unstablcrit} and hence $\lambda_1(D_{\eps, \theta,r})<N-1$. Moreover, by construction $D_{\eps, \theta,r} \subset \S_+^{N-1}$ and also the inequality $\mathcal{H}_{N-1}(D_{\eps, \theta,r})<\mathcal{H}_{N-1}(\S_+^{N-1})$ holds. Note that the domain $D_{\eps, \theta,r}$ is not smooth but we can take a smooth domain close to $D_{\eps, \theta,r}$ for which the same properties hold.
\end{example}

Next we exhibit another class of non-convex domains satisfying condition (i) of Proposition \ref{lemma:unstablcrit} which are not contained in a hemisphere.
\begin{example}\label{example2}
Let us fix $k\in\{1,\ldots,N-2\}$ and let
\beq\label{defSk}
\S^k:=\left\{(x_1,\ldots,x_{k+1},0,\ldots,0)\in \RN; \ \sum_{i=1}^{k+1} x_i^2=1 \right\}\subset \S^{N-1}.
\eeq
Moreover, we fix $r\in (0,\frac{\pi}{2})$ and consider
$$D_r:=\{x\in \S^{N-1}; \ \mathrm{dist}_{\S^{N-1}}(x,\S^k)<r\},$$
where $\mathrm{dist}_{\S^{N-1}}$ denotes the geodesic distance in $\S^{N-1}$. If $e\in \S^k$, we have $\int_{D_r} u_e \ d\sigma=0$ since $D_r$ is symmetric with respect to reflection at the hyperplane 
$$H_e:=\{x\in \RN; \ x \boldsymbol{\cdot} e=0\}$$
and $u_e$ is odd with respect to this reflection. We write points in $\S^{N-1}$ as
$$x= y\cos\theta + z \sin\theta,$$
with $y\in \S^k$ (see \eqref{defSk}), $\theta=\mathrm{dist}_{\S^{N-1}}(x,\S^k)\in(0,\frac{\pi}{2})$ and
$$z\in \S^{N-2-k}:=\left\{(0,\ldots,0,x_{k+2},\ldots,x_{N})\in \RN; \ \sum_{k+2}^N x_i^2=1 \right\}\subset \S^{N-1}.$$
In these coordinates, and since $e\in \S^k$, we have $u_e(x)= (e \boldsymbol{\cdot} y)\cos\theta.$
In addition we check that
$$\partial D_r= \{ x \in \S^{N-1}; \ x=y\cos r + z \sin r,\  y\in \S^k, \ z\in\S^{N-2-k} \},$$
and the exterior unit co-normal in a point $x=y\cos r + z \sin r \in \partial D_r$ is given by
$\nu(x)=-y \sin r + z\cos r$. Consequently, for any $x\in \partial D_r$ we have
$$\frac{\partial u_e}{\partial \nu}(x)= (\cos r) e \boldsymbol{\cdot} (-y \sin r + z\cos r)= - \cos r \sin r\ (e \boldsymbol{\cdot} y)=- \sin r \ u_e(x).$$
Hence it follows that
$$ u_e \frac{\partial u_e}{\partial \nu} <0 \ \ \hbox{on $\partial D_r\setminus H_e$,}$$
and thus
\beq\label{eq:lastineqappendix}
\int_{\partial D_r} u_e \frac{\partial u_e}{\partial \nu} \ d\hat\sigma <0.
\eeq
By Proposition \ref{lemma:unstablcrit}, the inequality \eqref{eq:lastineqappendix} implies that $\lambda_1(D_r)<N-1$. Finally if $r$ is small also the condition $\mathcal{H}_{N-1}(D_r)<\mathcal{H}_{N-1}(\S_+^{N-1})$ holds.
\end{example}


\begin{thebibliography}{10}

\bibitem{AM} 
Yu. Alkhutov, V. G. Maz'ya, \emph{$L^{1,p}$-coercitivity and estimates of the Green function of the Neumann problem in a convex domain}, Journal of Mathematical Sciences \textbf{196}, no. 3, 245--261 (2014)

\bibitem{BF} E. Baer, A. Figalli, \emph{Characterization of isoperimetric sets inside almost-convex cones}, Discrete Contin. Dyn. Syst. \textbf{37}(1), 1--14 (2017).

\bibitem{BHS} C. Bianchini, A. Henrot,  P. Salani, \emph{An overdetermined problem with non-constant boundary condition}, Interfaces and Free Boundaries \textbf{16}, 215--241 (2014).

\bibitem{BI2} D. Bonheure, A. Iacopetti, \emph{Spacelike radial graphs of prescribed mean curvature in the Lorentz-Minkowski space}, Analysis  \& PDE \textbf{12}, n. 7, 1805--1842 (2019).


\bibitem{BRP} T. Brian\c{c}on, M. Hayouni, M. Pierre, \emph{Lipschitz continuity of state functions
in some optimal shaping}, Calculus of Var. and PDE's \textbf{23}, 13--32 (2005).

\bibitem{BU} D. Bucur, \emph{Uniform concentration-compactness for Sobolev spaces on variable domains}. Journal of Differential Equations \textbf{162}, 427--450 (2000).

\bibitem{BU2}  D. Bucur, \emph{Minimization of the k-th eigenvalue of the Dirichlet Laplacian}, Arch. Rational Mech. Anal. \textbf{206}, 1073--1083 (2012). 

\bibitem{BV1} G. Buttazzo, B. Velichkov, \emph{Shape optimization problems on metric measure spaces}, Journal of Functional Analysis \textbf{264}, 1--33 (2013).

\bibitem{BV} G. Buttazzo, B. Velichkov, \emph{The spectral drop problem},  Contemporary mathematics \textbf{666}, 111--135 (2016).

\bibitem{CRS} X. Cabr\'e, X. Ros-Oton, J. Serra, \emph{Sharp isoperimetric inequalities via the ABP method}, J. Eur. Math. Soc. \textbf{18}, 2971--2998 (2016).

\bibitem{DJ} D. De Silva, D. Jerison, \emph{A singular energy minimizing free boundary}, J. Reine Angew. Math. \textbf{635}, 1--21 (2009).

\bibitem{ES} J. F. Escobar, \emph{Uniqueness theorems on conformal deformation of metrics, sobolev inequalities, and an eigenvalue estimate}, Comm. Pure Appl. Math. \textbf{43}(7), 857--883 (1990).

\bibitem{FI} A. Figalli, E. Indrei,  \emph{A Sharp Stability Result for the Relative Isoperimetric Inequality Inside Convex Cones}, J. Geom. Anal. \textbf{23}, 938--969 (2013).

 
 \bibitem{HP} A. Henrot, M. Pierre, \emph{Shape Variation and Optimization. A Geometrical Analysis}. EMS Tracts in Mathematics \textbf{28}, European Mathematical Society (2018).
 
 \bibitem{JS} D.S. Jerison, O. Savin, \emph{Some remarks on stability of cones for the one phase free boundary problem}, Geom. Funct. Anal. \textbf{25}, 1240--1257 (2015).

\bibitem{LS}  J. Lamboley, P. Sicbaldi, \emph{New Examples of Extremal Domains for the First Eigenvalue of the Laplace-Beltrami Operator in a Riemannian Manifold with Boundary}, International Mathematics Research Notices \textbf{18}, 8752--8798 (2015).


\bibitem{LIO} P. L. Lions, \emph{The concentration-compactness principle in the calculus of variations. The locally compact case, part 1}, Annales de l'I.H.P., section C \textbf{1}, n. 2, 109--145 (1984).

\bibitem{LP} P. L. Lions, F. Pacella, \emph{Isoperimetric inequalities for convex cones}, Proc.
Amer. Math. Soc. \textbf{109}, n. 2, 477--485 (1990)

\bibitem{LPT} P. L. Lions, F. Pacella, M. Tricarico, \emph{Best Constants in Sobolev Inequalities for Functions Vanishing on Some Part of the Boundary and Related Questions}, Indiana Univ. Math. Journal \textbf{37}, n. 2 (1988).

\bibitem{Lo03} 
{R. L\'opez}, \emph{A note on radial graphs with constant mean curvature}, 
{manuscripta math.} \textbf{110}, 45--54  (2003).

\bibitem{MPV} D. Mazzoleni, B. Pellacci, G. Verzini, \emph{Asymptotic spherical shapes in some spectral optimization problems}, Journal de Math\'ematiques Pures et Appliqu\'ees
\textbf{135}, 256--283 (2020).

\bibitem{PT}  F. Pacella, G. Tralli, \emph{Overdetermined problems and constant mean curvature surfaces in cones}, Rev. Mat. Iberoam. (2018) doi 10.4171/rmi/1151.

\bibitem{PT2}  F. Pacella, G. Tralli, \emph{Isoperimetric cones and minimal solutions of partial overdetermined problems},  Publ. Mat. \textbf{65}(1): 61--81 (2021). 

\bibitem{PTR} F. Pacella, M. Tricarico, \emph{Symmetrization for a Class of Elliptic Equations with Mixed Boundary Bonditions}, Atti Sem. Mat. Fis. Univ. Modena XXXIV, 75--94 (1985-86).

\bibitem{RR} M. Ritor\'e, C. Rosales, \emph{Existence and characterization of regions minimizing perimeter under a volume constraint inside euclidean cones}, Trans. Amer. Math. Soc. \textbf{356}(11), 4601--4622  (2004).

\bibitem{TW} A. Treibergs, W. Wei, \emph{Embedded hyperspheres with prescribed mean curvature}, J. Differential Geom. \textbf{18}, 513--521 (1983).

\bibitem{W} G. S. Weiss, \emph{A homogeneity improvement approach to the obstacle problem}, Invent. Math. \textbf{138}(1), 23--50 (1999).

\bibitem{WI} M. Willem,  \emph{Minimax Theorems}, Progress in Nonlinear Differential Equations and their Applications \textbf{24}, Birkh\"auser Basel (1996).
\end{thebibliography}
\end{document}